# BAYESIAN FREQUENTIST HYBRID INFERENCE[1]

By Ao Yuan

*Howard University*

Bayesian and frequentist methods differ in many aspects, but share some basic optimality properties. In practice, there are situations in which one of the methods is more preferred by some criteria. We consider the case of inference about a set of multiple parameters, which can be divided into two disjoint subsets. On one set, a frequentist method may be favored and on the other, the Bayesian. This motivates a joint estimation procedure in which some of the parameters are estimated Bayesian, and the rest by the maximum-likelihood estimator in the same parametric model, and thus keep the strengths of both the methods and avoid their weaknesses. Such a hybrid procedure gives us more flexibility in achieving overall inference advantages. We study the consistency and high-order asymptotic behavior of the proposed estimator, and illustrate its application. Also, the results imply a new method for constructing objective prior.

**1. Introduction.** In statistical practice, usually either the frequentist or the Bayesian method is used in parametric inference. Often the choice of methods is subjective. The two methods share some common basic asymptotical properties, which have been studied extensively. The Bernstein–von Mises theorem, for example [23, 30], states that in broad cases the Bayes and frequentist inferences are equivalent: the two estimators are close, and the posterior distribution around its mean is close to the distribution of the maximum likelihood estimate (MLE) around the true parameter—both are asymptotically normal with mean zero and the same asymptotic covariance matrix. However, the two methods are different in many other aspects, each of them has its own advantage(s) in some situations. In application they have received different appreciations for various reasons. Although all admissible

Received January 2008; revised August 2008.

[1]Supported in part by the National Center for Research Resources NIH Grant 212RR003048.

*AMS 2000 subject classifications.* Primary 62F10; secondary 62F15.

*Key words and phrases.* Bayesian inference, consistency, frequentist inference, high-order expansion, hybrid inference.

solutions, including the MLE, to a decision problem can be formulated as Bayesian [37], the two methods are regarded as different in concept, theory and history; so are they in this paper.

The Bayesian has appreciated steady growth partially due to the development of computation facilities, but in practice the main stream statistical tool is still frequentist. Efron [12] summarized the main reasons for this as the ease of use, modeling and objectivity. Lindley [25] gave a broad review of the present position of Bayesian statistics. The two schools mainly favor their own method, and the practitioners often have to choose one of the methods and ignore the other. In practice, there are situations in which one of the methods is more favorable than the other by some criteria. Thus in inferring multiparameters, it may happen that on part of the parameters, the frequentist method is preferable, while on the other part a Bayesian procedure is more appropriate. A practical example comes from our analysis of genetic data, in which the means of traits underlining each genotype are well studied in the literature. The sound prior knowledge prefers a Bayesian on this subset of parameters, while the mixing proportions and variances for the subdistributions are new in the investigation, and the MLE is favored on these parameters. This motivates a joint operation of the two methods on different parts of the parameters in the same model. Such hybrid inference will give us more flexibility than using either methods alone in achieving overall advantage. In this paper, we propose a hybrid estimator and study its consistency and asymptotic high-order behavior, and we illustrate its application. Also, using the high-order expansions, we considered a new type of second-order matching prior in the objective Bayes context.

There are some combined Bayesian and frequentist methods [1] and compromises between the two [16]. Our method here is not such combination nor compromise, not the quasi-Bayesian (pseudo- or semi-Bayesian), nor the empirical Bayes in the literature. The profile likelihood approach [28, 32] also divides the parameters into two parts, one of interest and the other nuisance, often of infinite dimensional. Fixing the parameters of interest, the nuisance parameters are eliminated by maximization, then the parameters of interest are estimated based on the resulting profile likelihood. This approach and ours have some operational similarity, but are different in nature. [7, 8, 24] study Bayesian method based on profile likelihood. They obtain the frequentist inference about parameters of interest by sampling from the posterior distribution of the semiparametric profile likelihood. They show that the estimator is of high-order accuracy to the corresponding frequentist's, and can have advantages on small samples. They further studied the case that the nuisance parameter may not have root-$n$ convergence rate. Their method and ours share more in common than the others. Although not identical, the two can be the same in some cases. This point will be made clear after we introduce the notations in the next section. Shen [33] studied the inference of parameters of interest by marginalizing the posterior distribution. Berger, Liseo and Wolpert [2] studied a method of eliminating nuisance parameters by integrating the likelihood over them with a prior, and the parameters of interest is estimated by maximizing the resulting likelihood.

In Section 2, we describe our method and study its consistency and asymptotic high-order behavior. For the latter, we first extend the existing high-order results of Bayes estimate and MLE to multivariate case, then based on the results, obtain the expansion of the hybrid estimator, which depicts the interplay of the Bayes and MLE components in each order. We show that the Bayes estimator, MLE and the hybrid estimator are first-order equivalent, asymptotic normal and efficient. In Section 3, we discuss implications from the results obtained, the evaluation of high-order terms, and advantages and weaknesses of MLE and Bayes, so as to have some references in choosing which method or how to hybridize them in practice. In particular, we derive a new method for constructing objective priors by matching high-order terms in expansions of the Bayes estimator and the MLE. We then illustrate the application of this method and the construction of the second-order objective prior in this sense by some examples. Also, a simple example is given in which neither the MLE nor Bayes estimator is consistent while the hybrid estimator is. The used regularity conditions and technical proofs of theorems are given in the Appendix.

**2. The method.** Let $f(\cdot|\boldsymbol{\theta})$ be a given density function of the data distribution and $\boldsymbol{\theta} \in \boldsymbol{\Theta} \subset R^d$ ($d > 1$) be the parameter of interest. Partition the parameter as $\boldsymbol{\theta} = (\boldsymbol{\alpha}, \boldsymbol{\beta}) \in (\boldsymbol{\Lambda}, \boldsymbol{\Omega}) \subset (R^{d_1}, R^{d_2})$. Assume that according to some criteria, on part of the parameters $\boldsymbol{\alpha}$, the Bayesian method is preferred, and on the other part $\boldsymbol{\beta}$, the MLE. This motivates the operation of the two methods on different parts of parameters in the same model simultaneously. We call such joint procedure of the two methods on the parameters in the same model a hybrid estimator, which is the goal of this study.

Specifically, let $\mathbf{x}^n = (\mathbf{x}_1, \ldots, \mathbf{x}_n)$ be an i.i.d. sample with likelihood $f(\mathbf{x}^n|\boldsymbol{\alpha}, \boldsymbol{\beta}) = \prod_{i=1}^n f(\mathbf{x}_i|\boldsymbol{\alpha}, \boldsymbol{\beta})$, $\pi(\boldsymbol{\alpha})$ be the prior density for $\boldsymbol{\alpha}$ and $\pi(\boldsymbol{\alpha}|\mathbf{x}^n, \boldsymbol{\beta}) \propto f(\mathbf{x}^n|\boldsymbol{\alpha}, \boldsymbol{\beta})\pi(\boldsymbol{\alpha})$ the posterior density of $\boldsymbol{\alpha}$ given the data and $\boldsymbol{\beta}$. Let $\mathcal{D}$ be the decision space for inferring $\boldsymbol{\alpha}$ ($\mathcal{D} = \boldsymbol{\Lambda}$ for estimation of $\boldsymbol{\alpha}$), $\mathbf{d}(\mathbf{x}^n) \in \mathcal{D}$ a decision rule, $W(\mathbf{d}(\mathbf{x}^n), \boldsymbol{\alpha})$ the loss function, $R(\mathbf{d}, \boldsymbol{\alpha}|\boldsymbol{\beta}) = E_{(\boldsymbol{\alpha}, \boldsymbol{\beta})} W(\mathbf{d}(\mathbf{x}^n), \boldsymbol{\alpha})$ the risk on $\boldsymbol{\alpha}$ for given $\boldsymbol{\beta}$, $R(\mathbf{d}|\boldsymbol{\beta}) = \int R(\mathbf{d}, \boldsymbol{\alpha}|\boldsymbol{\beta}) \pi(\boldsymbol{\alpha}) d\boldsymbol{\alpha}$ the Bayes risk on $\boldsymbol{\alpha}$ for given $\boldsymbol{\beta}$ and

$$R(\mathbf{d}|\mathbf{x}^n, \boldsymbol{\beta}) = \int W(\mathbf{d}(\mathbf{x}^n), \boldsymbol{\alpha}) \pi(\boldsymbol{\alpha}|\mathbf{x}^n, \boldsymbol{\beta}) d\boldsymbol{\alpha}$$

the posterior risk for inferring $\boldsymbol{\alpha}$ given $\boldsymbol{\beta}$. Then for fixed $\boldsymbol{\beta}$, the Bayes decision for $\boldsymbol{\alpha}$ is $\mathbf{d}^*(\cdot) = \mathbf{d}^*(\cdot|\boldsymbol{\beta}) = \arg\inf_{\mathbf{d} \in \boldsymbol{\Lambda}} R(\mathbf{d}|\boldsymbol{\beta})$, and from Bayes inference

theory

$$\mathbf{d}^*(\mathbf{x}^n) \stackrel{\text{a.s.}}{=} \arg\inf_{\mathbf{d} \in \mathbf{\Lambda}} R(\mathbf{d}|\mathbf{x}^n, \boldsymbol{\beta}) = \arg\inf_{\mathbf{d} \in \mathbf{\Lambda}} \int W(\mathbf{d}(\mathbf{x}^n), \boldsymbol{\alpha}) f(\mathbf{x}^n|\boldsymbol{\alpha}, \boldsymbol{\beta}) \pi(\boldsymbol{\alpha}) \, d\boldsymbol{\alpha}.$$

The right-hand side above is the generalized Bayesian estimator of $\boldsymbol{\alpha}$ for fixed $\boldsymbol{\beta}$.

In this hybrid inference, we infer $\boldsymbol{\alpha}$ by the generalized Bayesian rule for each fixed $\boldsymbol{\beta}$ and at the same time infer $\boldsymbol{\beta}$ by the frequentist MLE, that is, we are to find $\boldsymbol{\theta}_n = (\check{\boldsymbol{\alpha}}_n, \hat{\boldsymbol{\beta}}_n) = (\check{\boldsymbol{\alpha}}(\mathbf{x}^n), \hat{\boldsymbol{\beta}}(\mathbf{x}^n))$ such that

$$(2.1) \qquad (\check{\boldsymbol{\alpha}}_n, \hat{\boldsymbol{\beta}}_n) = \arg\inf\sup_{(\mathbf{d}, \boldsymbol{\beta})} \int W(\mathbf{d}(\mathbf{x}^n), \boldsymbol{\alpha}) f(\mathbf{x}^n|\boldsymbol{\alpha}, \boldsymbol{\beta}) \pi(\boldsymbol{\alpha}) \, d\boldsymbol{\alpha}$$

is the joint optimizer over $(\mathbf{d}, \boldsymbol{\beta}) \in (\mathbf{\Lambda}, \mathbf{\Omega})$.

REMARK 1. By imposing a 0–1 loss and constant prior on $\boldsymbol{\beta}$, (2.1) can be formulated as a full Bayesian solution (as in the proof of Theorem 2.1). Thus, $(\check{\boldsymbol{\alpha}}_n, \hat{\boldsymbol{\beta}}_n)$ generally exists and is locally unique.

In the above $(\check{\boldsymbol{\alpha}}_n, \hat{\boldsymbol{\beta}}_n)$ is jointly a generalized Bayes estimator and a MLE of $(\boldsymbol{\alpha}, \boldsymbol{\beta})$. Here $\check{\boldsymbol{\alpha}}_n$ is not identical to the Bayes estimator based on profile likelihood such as in [8]. The latter is first eliminating the nuisance parameter $\boldsymbol{\beta}$ by maximizing the likelihood over it, along some least favorable curve, to get $\tilde{\boldsymbol{\beta}} = \tilde{\boldsymbol{\beta}}(\boldsymbol{\alpha}) = \arg\sup_{\boldsymbol{\beta}} f(\mathbf{x}^n|\boldsymbol{\alpha}, \boldsymbol{\beta})$, then computing the Bayes estimator $\tilde{\boldsymbol{\alpha}}_n$ of $\boldsymbol{\alpha}$ based on the profile likelihood $\tilde{f}(\mathbf{x}^n|\boldsymbol{\alpha}) = f(\mathbf{x}^n|\boldsymbol{\alpha}, \tilde{\boldsymbol{\beta}}(\boldsymbol{\alpha}))$, that is, $\tilde{\boldsymbol{\alpha}}_n = \arg\inf_{\mathbf{d}} \int W(\mathbf{d}(\mathbf{x}^n), \boldsymbol{\alpha}) \tilde{f}(\mathbf{x}^n|\boldsymbol{\alpha}) \pi(\boldsymbol{\alpha}) \, d\boldsymbol{\alpha} = \arg\inf_{\mathbf{d}} \int W(\mathbf{d}(\mathbf{x}^n), \boldsymbol{\alpha}) \sup_{\boldsymbol{\beta}} f(\mathbf{x}^n|\boldsymbol{\alpha}, \boldsymbol{\beta}) \pi(\boldsymbol{\alpha}) \, d\boldsymbol{\alpha}$. It is seen that generally $\check{\boldsymbol{\alpha}}_n \neq \tilde{\boldsymbol{\alpha}}_n$, and they can be equal under some fair conditions, such as that the integrand in (2.1) can be dominated, with respect to $(\mathbf{x}^n, \boldsymbol{\beta})$, by some integrable function in $\boldsymbol{\alpha}$.

In the following, we study the consistency of the hybrid estimator. As Bayes estimator and the MLE are generally first-order equivalent, their competition goes into the asymptotic high-order terms. We investigate the high-order asymptotic behavior of (2.1); this will give us flexibility in choosing which method to use on which parameter component(s) to achieve high-order advantage.

*Consistency of the estimator.* The study of the consistencies of Bayes estimates, MLE and their relationships has a relatively long history [4, 11, 22, 34, 36], among others. Doob [11] established strong consistency of Bayes estimators under very general conditions, and there is some speculation that conditions for Bayesian consistency might be found which are weaker than those for the MLE. Under some basic assumptions, Strasser [34] showed that any conditions for the convergence (a.s.) of MLE assert the concentration (a.s.) of the posterior distribution to the true parametric value. This does not directly imply that conditions for Bayesian consistency are weaker since posterior concentration to the true parameter is not equivalent to the consistency of Bayes estimate. The latter also depends on the loss. There are examples in which one of the estimator is consistent while the other is not. However, for multiple parameters using a hybrid estimator may overcome possible difficulty in using one method alone. We will give such an example later. Let $\boldsymbol{\theta}_0 = (\boldsymbol{\alpha}_0, \boldsymbol{\beta}_0)$ be the "true" parameters generating the data under the specified model. In the following, we give the consistency of the hybrid estimate using a method similar to that in [5]. Generally, the loss $W(\mathbf{d}, \boldsymbol{\alpha})$ has the form $W(\|\mathbf{d} - \boldsymbol{\alpha}\|) = W(\mathbf{d} - \boldsymbol{\alpha})$. To avoid confusion, we will use $W$ for any of these functional forms. Let $d = dim(\boldsymbol{\theta}) = dim(\boldsymbol{\alpha}) + dim(\boldsymbol{\beta}) = d_1 + d_2$.

THEOREM 2.1. *Assume conditions* (A1)–(A9) *in the Appendix*, $W(\cdot)$ *satisfies* $W(\mathbf{0}) = 0$, *is strictly increasing and continuous in a neighborhood of* $\mathbf{0}$. *Then, as* $n \to \infty$ *we have*

$$(\check{\boldsymbol{\alpha}}_n, \hat{\boldsymbol{\beta}}_n) \to (\boldsymbol{\alpha}_0, \boldsymbol{\beta}_0) \qquad (a.s.).$$

*High-order asymptotic behavior.* High-order asymptotic expansions are used to assess estimators when they have similar lower-order behavior. In [9, 17, 20, 26] (among others) such expansions of Bayes estimate and MLE were obtained, so were their densities and related quantities in the one-dimensional case. The results in [17] are more suitable to our case. Here, we first generalize the results there to the multidimensional case, then use them to get the expansion of the hybrid estimator.

We first give a multivariate generalization of the asymptotic expansion of maximum posterior density estimator, which is given by

$$\hat{\boldsymbol{\theta}}_n = \arg\sup_{\boldsymbol{\theta}} \pi(\boldsymbol{\theta}|\mathbf{x}^n) = \arg\sup_{\boldsymbol{\theta}} f(\mathbf{x}^n|\boldsymbol{\theta})\pi(\boldsymbol{\theta}).$$

We introduce the following notation. For an integer vector $\mathbf{i} = (i_1, \ldots, i_d)$ with $i_j \geq 0$ $(j = 1, \ldots, d)$, denote $|\mathbf{i}| = \sum_{j=1}^d i_j$, $\frac{\partial^{|\mathbf{i}|}}{\partial \boldsymbol{\theta}^\mathbf{i}} = \frac{\partial^{|\mathbf{i}|}}{\partial \theta_1^{i_1} \cdots \partial \theta_d^{i_d}}$, and for any $g(\cdot) \geq 0$, define $\log g(\cdot) = 0$ if $g(\cdot) = 0$. Denote $\mathbf{1} = (1, \ldots, 1)'$ of length $d$, $\mathbf{0} = (0, \ldots, 0)'$ of length $d$,

$$\mathbf{L}(\mathbf{x}|\boldsymbol{\theta}) := (L_1(\mathbf{x}|\boldsymbol{\theta}), \ldots, L_d(\mathbf{x}|\boldsymbol{\theta}))' = \left(\frac{\partial}{\partial \theta_1} \log f(\mathbf{x}|\boldsymbol{\theta}), \ldots, \frac{\partial}{\partial \theta_d} \log f(\mathbf{x}|\boldsymbol{\theta})\right)',$$

$$\mathbf{L_i}(\mathbf{x}|\boldsymbol{\theta}) = \left(\frac{\partial^{|\mathbf{i}|}}{\partial \boldsymbol{\theta}^\mathbf{i}} L_1(\mathbf{x}|\boldsymbol{\theta}), \ldots, \frac{\partial^{|\mathbf{i}|}}{\partial \boldsymbol{\theta}^\mathbf{i}} L_d(\mathbf{x}|\boldsymbol{\theta})\right)',$$

$$\boldsymbol{\rho}(\boldsymbol{\theta}) := (\rho_1(\boldsymbol{\theta}), \ldots, \rho_d(\boldsymbol{\theta}))' = \left(\frac{\partial}{\partial \theta_1} \log \pi(\boldsymbol{\theta}), \ldots, \frac{\partial}{\partial \theta_d} \log \pi(\boldsymbol{\theta})\right)',$$

$$\boldsymbol{\rho}_{\mathbf{i}}(\boldsymbol{\theta}) = \left(\frac{\partial^{|\mathbf{i}|}}{\partial \boldsymbol{\theta}^{\mathbf{i}}}\rho_1(\boldsymbol{\theta}), \ldots, \frac{\partial^{|\mathbf{i}|}}{\partial \boldsymbol{\theta}^{\mathbf{i}}}\rho_d(\boldsymbol{\theta})\right)',$$

$$\mathbf{S}_{\mathbf{i}}(\boldsymbol{\theta}) = \frac{1}{\sqrt{n}}\sum_{j=1}^n \mathbf{L}_{\mathbf{i}}(\mathbf{x}_j|\boldsymbol{\theta}),$$

$$\boldsymbol{\Delta}_{\mathbf{i}}(\boldsymbol{\theta}) = \frac{1}{\sqrt{n}}\sum_{j=1}^n (\mathbf{L}_{\mathbf{i}}(\mathbf{x}_j|\boldsymbol{\theta}) - E_{\boldsymbol{\theta}}\mathbf{L}_{\mathbf{i}}(\mathbf{x}_1|\boldsymbol{\theta})),$$

and set $\mathbf{S}_{\mathbf{i}} = \mathbf{S}_{\mathbf{i}}(\boldsymbol{\theta}_0)$, $\boldsymbol{\Delta}_{\mathbf{i}} = \boldsymbol{\Delta}_{\mathbf{i}}(\boldsymbol{\theta}_0)$ and $\mathbf{E}_{\mathbf{i}} = E_{\boldsymbol{\theta}_0}\mathbf{L}_{\mathbf{i}}(\mathbf{x}_1|\boldsymbol{\theta}_0)$. For vector $\mathbf{H} = (H_1, \ldots, H_d)'$ and integer vector $\mathbf{i} = (i_1, \ldots, i_d)$, define $\mathbf{H}^{\mathbf{i}} = (H_1^{i_1}, \ldots, H_d^{i_d})'$, $\langle \mathbf{H}^{\mathbf{i}} \rangle = \prod_{j=1}^d H_j^{i_j}$ and $\mathbf{i}! = \prod_{j=1}^d i_j!$. For $\mathbf{a} = (a_1, \ldots, a_d)'$ and $\mathbf{b} = (b_1, \ldots, b_d)'$, define $\mathbf{a} + \mathbf{b} = (a_1 + b_1, \ldots, a_d + b_d)'$, $\mathbf{ab} = (a_1 b_1, \ldots, a_d b_d)'$ and $\langle \mathbf{ab} \rangle = \prod_{i=1}^d a_i b_i$. Denote $\mathbf{e}_j = (0, \ldots, 0, 1, 0, \ldots, 0)'$, the $d$-vector with $j$th element be 1 and the others be zeros. For nonnegative integers $r \leq s$ and nonnegative integer $d$-vectors $\mathbf{l}$ and $\mathbf{i}$, the notation $(r, s, \mathbf{l}, \mathbf{i})$ stands for the collection of all nonnegative integer $d$-vector sets $\{(\mathbf{i}_r, \ldots, \mathbf{i}_s)\}$,

$$(r, s, \mathbf{l}, \mathbf{i}) = \left\{(\mathbf{i}_r, \ldots, \mathbf{i}_s) : \sum_{v=r}^s v\mathbf{i}_v = \mathbf{l}, \sum_{v=r}^s \mathbf{i}_v = \mathbf{i}\right\}.$$

Denote $\mathbf{I}$ the Fisher information, and $\mathbf{I}^{-1}$ its inverse, evaluated at $\boldsymbol{\theta}_0$.

THEOREM 2.2. *Under conditions* (B1)–(B7) *in the [Appendix](), we have*

$$\sqrt{n}(\hat{\boldsymbol{\theta}}_n - \boldsymbol{\theta}_0) = \sum_{r=0}^{k-1} n^{-r/2}\mathbf{H}_r + O_p(n^{-k/2}),$$

*where the term $O_p(n^{-k/2})$ is in the sense of* [17], *and the $\mathbf{H}_r$'s are $d$-vectors of polynomials in the $\boldsymbol{\Delta}_{\mathbf{i}}$'s of degree $r+1$, their coefficients are polynomials in the $\mathbf{E}_{\mathbf{i}}$'s, $|\mathbf{i}| = 2, \ldots, r+1$, in $\mathbf{I}^{-1}$ and the $\boldsymbol{\rho}_{\mathbf{i}}$'s, given by* $(0 \leq r \leq k-1)$

$$\mathbf{H}_0 = \mathbf{I}^{-1}\boldsymbol{\Delta}_0,$$

$$\mathbf{H}_r = \mathbf{I}^{-1}\sum_{s+t=r}\left(\sum_{|\mathbf{i}|=t-1}\boldsymbol{\rho}_{\mathbf{i}}\sum_{|\mathbf{l}|=s}\sum_{(0,s,\mathbf{l},\mathbf{i})}\prod_{v=0}^s \frac{\langle \mathbf{H}_v^{\mathbf{i}_v} \rangle}{\mathbf{i}_v!} + \sum_{|\mathbf{i}|=t}\boldsymbol{\Delta}_{\mathbf{i}}\sum_{|\mathbf{l}|=s}\sum_{(0,s,\mathbf{l},\mathbf{i})}\prod_{v=0}^s \frac{\langle \mathbf{H}_v^{\mathbf{i}_v} \rangle}{\mathbf{i}_v!}\right.$$

$$\left. + \sum_{|\mathbf{i}|=t+1, t>0}\mathbf{E}_{\mathbf{i}}\sum_{|\mathbf{l}|=s}\sum_{(0,s,\mathbf{l},\mathbf{i})}\prod_{v=0}^s \frac{\langle \mathbf{H}_v^{\mathbf{i}_v} \rangle}{\mathbf{i}_v!}\right).$$

In the above set $\pi(\cdot)$ to constant, then $\hat{\boldsymbol{\theta}}_n$ is the MLE, and we have

$$\sqrt{n}(\hat{\boldsymbol{\theta}}_n - \boldsymbol{\theta}_0) = \sum_{r=0}^{k-1} n^{-r/2}\mathbf{H}_r^{\circ} + O_p(n^{-k/2}),$$

where $\mathbf{H}_0^\circ = \mathbf{H}_0$, and for $1 \leq r \leq k-1$,

$$\mathbf{H}_r^\circ = \mathbf{I}^{-1} \sum_{s+t=r} \left( \sum_{|\mathbf{i}|=t} \boldsymbol{\Delta}_{\mathbf{i}} \sum_{|\mathbf{l}|=s} \sum_{(0,s,\mathbf{l},\mathbf{i})} \prod_{v=0}^{s} \frac{\langle \mathbf{H}_v^{\circ \mathbf{i}_v} \rangle}{\mathbf{i}_v!} \right.$$

$$\left. + \sum_{|\mathbf{i}|=t+1, t>0} \mathbf{E}_{\mathbf{i}} \sum_{|\mathbf{l}|=s} \sum_{(0,s,\mathbf{l},\mathbf{i})} \prod_{v=0}^{s} \frac{\langle \mathbf{H}_v^{\circ \mathbf{i}_v} \rangle}{\mathbf{i}_v!} \right).$$

Next, set $\boldsymbol{\theta} = \boldsymbol{\alpha}$, then $\boldsymbol{\theta}_n = \check{\boldsymbol{\theta}}_n$ is the Bayes estimate of $\boldsymbol{\theta}$. To get the corresponding expansion for $\sqrt{n}(\check{\boldsymbol{\theta}}_n - \boldsymbol{\theta}_0)$, we need the following notation. Let $l(x|\boldsymbol{\theta}) = \log f(x|\boldsymbol{\theta})$, $\varrho(\boldsymbol{\theta}) = \log \pi(\boldsymbol{\theta})$,

$$l_{\mathbf{i}}(x|\boldsymbol{\theta}) = \frac{\partial^{|\mathbf{i}|}}{\partial \boldsymbol{\theta}^{\mathbf{i}}} l(x|\boldsymbol{\theta}), \qquad \varrho_{\mathbf{i}}(\boldsymbol{\theta}) = \frac{\partial^{|\mathbf{i}|}}{\partial \boldsymbol{\theta}^{\mathbf{i}}} \log \pi(\boldsymbol{\theta}),$$

$$\mathcal{S}_{\mathbf{i}}(\boldsymbol{\theta}) = \frac{1}{\sqrt{n}} \sum_{j=1}^{n} l_{\mathbf{i}}(x_j|\boldsymbol{\theta}), \qquad \delta_{\mathbf{i}}(\boldsymbol{\theta}) = \frac{1}{\sqrt{n}} \sum_{j=1}^{n} (l_{\mathbf{i}}(x_j|\boldsymbol{\theta}) - E_{\boldsymbol{\theta}} l_{\mathbf{i}}(\mathbf{x}_1|\boldsymbol{\theta})),$$

and set $\mathcal{S}_{\mathbf{i}} = \mathcal{S}_{\mathbf{i}}(\boldsymbol{\theta}_0)$, $\delta_{\mathbf{i}} = \delta_{\mathbf{i}}(\boldsymbol{\theta}_0)$ and $\mathcal{E}_{\mathbf{i}} = E_{\boldsymbol{\theta}_0} l_{\mathbf{i}}(X|\boldsymbol{\theta}_0)$. We make the convention that for nonnegative integer vectors $\mathbf{i} = (i_1, \ldots, i_d)'$ and $\mathbf{j} = (j_1, \ldots, j_d)'$, the notation $\mathbf{j} - \mathbf{i}$ implies $\mathbf{j} \geq \mathbf{i}$, that is, $j_r \geq i_r$ $(r = 1, \ldots, d)$.

THEOREM 2.3. *Under conditions* (B1)–(B10) *in the Appendix, we have*

$$\sqrt{n}(\check{\boldsymbol{\theta}}_n - \boldsymbol{\theta}_0) = \sum_{r=0}^{k-1} n^{-r/2} \mathbf{G}_r + O_P(n^{-k/2}),$$

*where* $\mathbf{G}_0 = \mathbf{H}_0$ *and for* $1 \leq r \leq k-1$,

$\mathbf{G}_r = \mathbf{H}_r + \mathbf{Q}_r$,

$$\mathbf{Q}_r = \mathbf{M}_{\mathbf{0},r} + \sum_{m=1}^{r-1} \sum_{|\mathbf{i}|=1}^{m} \mathbf{i}! \mathbf{M}_{\mathbf{i}, r-m} \sum_{|\mathbf{l}|=m} \sum_{(1,m,\mathbf{l},\mathbf{i})} \prod_{v=1}^{m} \langle \mathbf{Q}_v^{\mathbf{i}_v} \rangle / \mathbf{i}_v!$$

$$+ \sum_{|\mathbf{i}| \in \langle 3, r \rangle} \mathbf{i}! \mathbf{M}_{\mathbf{i}, 0} \sum_{|\mathbf{l}|=r} \sum_{(1,r-1,\mathbf{l},\mathbf{i})} \prod_{v=1}^{r-1} \langle \mathbf{Q}_v^{\mathbf{i}_v} \rangle / \mathbf{i}_v! \qquad (note\ \mathbf{Q}_1 = \mathbf{M}_{\mathbf{0},1}),$$

$$\mathbf{M}_{\mathbf{j},r} = (\{\boldsymbol{\sigma}(\mathbf{a})\}\mathbf{I})^{-1} \sum_{|\mathbf{i}| \in \langle 2(1 \wedge r) + |\mathbf{j}|, 3r+|\mathbf{j}| \rangle} \frac{1}{\mathbf{j}!} N_{\mathbf{i}-\mathbf{j},r} \boldsymbol{\Psi}_{\mathbf{i}-\mathbf{j}}^{(\mathbf{j})} \qquad (0 \leq |\mathbf{j}| \leq k-1),$$

*where* $\{\boldsymbol{\sigma}(\mathbf{a})\} = \mathrm{diag}(\sigma(a_1), \ldots, \sigma(a_d))$, $\sigma(a_r)$ *is the* $a_r$*th marginal moment of* $\theta_r$ *with* $\boldsymbol{\theta} \sim N(\mathbf{0}, \mathbf{I}^{-1})$, $\boldsymbol{\Psi}_{\mathbf{i}}^{(\mathbf{j})} = \boldsymbol{\Psi}_{\mathbf{i}}^{(\mathbf{j})}(\mathbf{0})$ $[\boldsymbol{\Psi}_{\mathbf{i}} := \boldsymbol{\Psi}_{\mathbf{i}}^{(\mathbf{0})}(\mathbf{0})]$ *is d-vector of multivariate normal moments associated with the loss [with components of the form* $\sigma(\mathbf{i})$ *as given in the proof],* $\langle a, b \rangle$ *is the set of odd integers* $s$ *with* $a \leq s \leq b$ *and* $N_{\mathbf{i},r}$*'s are random variables given in Lemma 1 in the Appendix.*

*For general loss functions instead of the one in* (B9), *the* $\mathbf{Q}_r$ *'s have more terms in more involved forms, and are outlined at the end of the proof.*

Now, based on the expansions of the MLE and Bayes estimate, we give asymptotic expansion for the hybrid estimator $\sqrt{n}(\boldsymbol{\theta}_n - \boldsymbol{\theta}_0)$, which depicts the interplay of the two components in each order. Denote $\mathbf{I}^{-1} = (\mathbf{I}^{ij})_{1 \le i,j \le 2}$ and $\mathbf{a} = (\mathbf{a}_1', \mathbf{a}_2')'$ as the componentwise notations corresponding to the two sets of parameters $(\boldsymbol{\alpha}, \boldsymbol{\beta})$.

THEOREM 2.4. *Assume conditions* (B1)–(B5), *and* (B6)–(B10) *(with $\boldsymbol{\theta}$, $\boldsymbol{\Theta}$ and $\mathbf{a}$ replaced by $\boldsymbol{\alpha}$, $\boldsymbol{\Lambda}$ and $\mathbf{a}_1$), in the Appendix, then*

$$\sqrt{n}(\boldsymbol{\theta}_n - \boldsymbol{\theta}_0) = \sqrt{n}\begin{pmatrix} \check{\boldsymbol{\alpha}}_n - \boldsymbol{\alpha}_0 \\ \hat{\boldsymbol{\beta}}_n - \boldsymbol{\beta}_0 \end{pmatrix}$$

$$= \sum_{r=0}^{k-1} n^{-r/2} \begin{pmatrix} \mathbf{g}_r \\ \mathbf{h}_r \end{pmatrix} + O_p(n^{-k/2}),$$

*where* $\begin{pmatrix} \mathbf{g}_r \\ \mathbf{h}_r \end{pmatrix}$ $(0 \le r \le k-1)$ *are given by*

$$\begin{pmatrix} \mathbf{g}_0 \\ \mathbf{h}_0 \end{pmatrix} = \mathbf{I}^{-1} \boldsymbol{\Delta}_0,$$

$$\begin{pmatrix} \mathbf{g}_r \\ \mathbf{h}_r \end{pmatrix} = \mathbf{I}^{-1} \sum_{s+t=r} \left( \sum_{|\mathbf{i}|=r} \boldsymbol{\Delta}_{\mathbf{i}} \sum_{|\mathbf{l}|=s} \sum_{(0,s,\mathbf{l},\mathbf{i})} \prod_{v=0}^{s} \frac{1}{\mathbf{i}_v!} \left\langle \begin{pmatrix} \mathbf{g}_v \\ \mathbf{h}_v \end{pmatrix}^{\mathbf{i}_v} \right\rangle \right.$$

$$\left. + \sum_{|\mathbf{i}|=t+1, t>0} \mathbf{E}_{\mathbf{i}} \sum_{|\mathbf{l}|=s} \sum_{(0,s,\mathbf{l},\mathbf{i})} \prod_{v=0}^{s} \frac{1}{\mathbf{i}_v!} \left\langle \begin{pmatrix} \mathbf{g}_v \\ \mathbf{h}_v \end{pmatrix}^{\mathbf{i}_v} \right\rangle \right) + \mathbf{T}_r$$

$$(r = 1, \ldots, k-1)$$

*and*

$$\mathbf{T}_r = \mathbf{I}^{-1} \sum_{s+t=r} \sum_{|\mathbf{i}|=t-1} \begin{pmatrix} \boldsymbol{\rho}_{\mathbf{i}}^1 \\ 0 \end{pmatrix} \sum_{|\mathbf{l}|=s} \sum_{(0,s,\mathbf{l},\mathbf{i})} \prod_{v=0}^{s} \frac{\langle \mathbf{H}_v^{\mathbf{i}_v} \rangle}{\mathbf{i}_v!} + \begin{pmatrix} \mathbf{q}_r \\ 0 \end{pmatrix},$$

$\boldsymbol{\rho}_{\mathbf{i}}^1 = \boldsymbol{\rho}_{\mathbf{i}}^1(\boldsymbol{\alpha}_0)$ *is the first $d_1$ components of $\boldsymbol{\rho}_{\mathbf{i}}$, $\mathbf{H}_r$'s as in Theorem 2.2 with $\boldsymbol{\rho}_{\mathbf{i}} = (\boldsymbol{\rho}_{\mathbf{i}}^1, 0)'$ and $\mathbf{q}_r = \mathbf{q}_r(\boldsymbol{\theta}_0)$ is the $d_1$-dimensional version of $\mathbf{Q}_r$ in Theorem 2.3 with $\mathbf{I}^{-1}$ replaced by $\mathbf{I}^{11}$.*

For general loss functions rather than that given in (B9), the results are analogous. Below we give the first three terms in the expansions for each estimators.

FACT. (i) For the MLE in Theorem 2.2, the first three terms are

$$\mathbf{H}_0^\circ = \mathbf{I}^{-1}\mathbf{\Delta_0}, \qquad \mathbf{H}_1^\circ = \mathbf{I}^{-1}\left(\sum_{|\mathbf{i}|=1}\mathbf{\Delta_i}\langle\mathbf{H}_0^{\circ\mathbf{i}}\rangle + \sum_{|\mathbf{i}|=2}\mathbf{E_i}\langle\mathbf{H}_0^{\circ\mathbf{i}}\rangle/\mathbf{i}!\right),$$

$$\mathbf{H}_2^\circ = \mathbf{I}^{-1}\left(\sum_{|\mathbf{i}|=2}\mathbf{\Delta_i}\langle\mathbf{H}_0^{\circ\mathbf{i}}\rangle/\mathbf{i}! + \sum_{|\mathbf{i}|=1}\mathbf{\Delta_i}\langle\mathbf{H}_1^{\circ\mathbf{i}}\rangle + \sum_{|\mathbf{i}|=3}\mathbf{E_i}\langle\mathbf{H}_0^{\circ\mathbf{i}}\rangle/\mathbf{i}!\right.$$

$$\left.+ \sum_{|\mathbf{i}|,|\mathbf{j}|=1}\mathbf{E_{i+j}}(\langle\mathbf{H}_0^{\circ\mathbf{i}}\mathbf{H}_1^{\circ\mathbf{j}}\rangle + \langle\mathbf{H}_0^{\circ\mathbf{j}}\mathbf{H}_1^{\circ\mathbf{i}}\rangle)/2\right).$$

(ii) For the maximum posterior estimator in Theorem 2.2, the first three terms are

$$\mathbf{H}_0 = \mathbf{H}_0^\circ, \qquad \mathbf{H}_1 = \mathbf{H}_1^\circ + \mathbf{I}^{-1}\boldsymbol{\rho_0},$$

$$\mathbf{H}_2 = \mathbf{I}^{-1}\left(\sum_{|\mathbf{i}|=1}\boldsymbol{\rho_i}\langle\mathbf{H}_0^\mathbf{i}\rangle + \sum_{|\mathbf{i}|=2}\mathbf{\Delta_i}\langle\mathbf{H}_0^\mathbf{i}\rangle/\mathbf{i}! + \sum_{|\mathbf{i}|=1}\mathbf{\Delta_i}\langle\mathbf{H}_1^\mathbf{i}\rangle + \sum_{|\mathbf{i}|=3}\mathbf{E_i}\langle\mathbf{H}_0^\mathbf{i}\rangle/\mathbf{i}!\right.$$

$$\left.+ \sum_{|\mathbf{i}|,|\mathbf{j}|=1}\mathbf{E_{i+j}}(\langle\mathbf{H}_0^\mathbf{i}\mathbf{H}_1^\mathbf{j}\rangle + \langle\mathbf{H}_0^\mathbf{j}\mathbf{H}_1^\mathbf{i}\rangle)/2\right).$$

(iii) For the Bayes estimator in Theorem 2.3, the first three terms are

$$\mathbf{G}_0 = \mathbf{H}_0^\circ, \qquad \mathbf{G}_1 = \mathbf{H}_1^\circ + \mathbf{I}^{-1}\boldsymbol{\rho_0} + \mathbf{Q}_1,$$

$$\mathbf{Q}_1 = \mathbf{M_{0,1}} = (\{\boldsymbol{\sigma}(\mathbf{a})\}\mathbf{I})^{-1}\sum_{|\mathbf{i}|=3}\boldsymbol{\Psi_i}\frac{\boldsymbol{\mathcal{E}_i}}{\mathbf{i}!},$$

$$\boldsymbol{\Psi_i} = E(\boldsymbol{\theta}^{\mathbf{a}-1}\langle\boldsymbol{\theta}^\mathbf{i}\rangle), \qquad \boldsymbol{\theta} \sim N(\mathbf{0},\mathbf{I}^{-1});$$

$$\mathbf{G}_2 = \mathbf{H}_2 + \mathbf{Q}_2, \qquad \mathbf{Q}_2 = \mathbf{M_{0,2}} + \sum_{|\mathbf{i}|=1}\mathbf{M_{i,1}}\langle\mathbf{Q}_1^\mathbf{i}\rangle,$$

$$\mathbf{M_{0,2}} = (\{\boldsymbol{\sigma}(\mathbf{a})\}\mathbf{I})^{-1}\sum_{|\mathbf{i}|=3}\frac{1}{\mathbf{i}!}\left(\boldsymbol{\delta_i} + \sum_{|\mathbf{j}|=1}\boldsymbol{\mathcal{E}_{i+j}}\langle(\mathbf{I}^{-1}\mathbf{\Delta_0})^\mathbf{j}\rangle\right)\boldsymbol{\Psi_i},$$

$$\mathbf{M_{i,1}} = (\{\boldsymbol{\sigma}(\mathbf{a})\}\mathbf{I})^{-1}\sum_{\mathbf{j}>\mathbf{i},|\mathbf{j}|=3}\frac{1}{(\mathbf{j}-\mathbf{i})!}\left(\boldsymbol{\delta_{j-i}} + \sum_{|\mathbf{l}|=1}\boldsymbol{\mathcal{E}_{j-i+l}}\langle(\mathbf{I}^{-1}\mathbf{\Delta_0})^\mathbf{l}\rangle\right)\boldsymbol{\Psi_i}.$$

(iv) For the hybrid estimator in Theorem 2.4, the first three terms are

$$\begin{pmatrix}\mathbf{g}_0\\\mathbf{h}_0\end{pmatrix} = \mathbf{H}_0^\circ, \qquad \begin{pmatrix}\mathbf{g}_1\\\mathbf{h}_1\end{pmatrix} = \mathbf{H}_1^\circ + \begin{pmatrix}\mathbf{I}^{11}\boldsymbol{\rho}_0^1 + (\{\boldsymbol{\sigma}(\mathbf{a}_1)\mathbf{I}^{11}\})^{-1}\sum_{|\mathbf{i}|=3}\boldsymbol{\Psi}_\mathbf{i}^1\boldsymbol{\mathcal{E}}_\mathbf{i}^1/\mathbf{i}!\\\mathbf{I}^{21}\boldsymbol{\rho}_0^1\end{pmatrix},$$

where $\boldsymbol{\Psi}_{\mathbf{i}}^1 = E(\boldsymbol{\alpha}^{\mathbf{a}_1 - \mathbf{1}_1} \langle \boldsymbol{\alpha}^{\mathbf{i}} \rangle)$, $\boldsymbol{\alpha} \sim N(\mathbf{0}, \mathbf{I}^{11})$, and $\mathcal{E}_{\mathbf{i}}^1 = E_{\boldsymbol{\theta}_0}[\frac{\partial^{|\mathbf{i}|}}{\partial \boldsymbol{\alpha}^{\mathbf{i}}} l(\mathbf{x}_1 | \boldsymbol{\alpha}, \boldsymbol{\beta}_0)]|_{\boldsymbol{\alpha} = \boldsymbol{\alpha}_0}$,

$$\begin{pmatrix} \mathbf{g}_2 \\ \mathbf{h}_2 \end{pmatrix} = \mathbf{H}_2^\bullet + \begin{pmatrix} \mathbf{I}^{11} \sum_{|\mathbf{i}|=1} \rho_{\mathbf{i}}^1 \langle \mathbf{H}_0^{\mathbf{i}} \rangle + \mathbf{q}_2 \\ \mathbf{I}^{21} \sum_{|\mathbf{i}|=1} \rho_{\mathbf{i}}^1 \langle \mathbf{H}_0^{\mathbf{i}} \rangle \end{pmatrix},$$

where $\mathbf{H}_2^\bullet$ is $\mathbf{H}_2^\circ$ in which $\mathbf{H}_r^\circ$ is replaced by $(\mathbf{g}_r', \mathbf{h}_r')'$ ($r = 0, 1$) and $\mathbf{q}_2$ is $\mathbf{Q}_2$ with $(\mathbf{I}, \mathbf{a}, \boldsymbol{\Psi}_{\mathbf{i}}, \mathcal{E}_{\mathbf{i}}, \delta_{\mathbf{i}}, \boldsymbol{\Delta}_0)$ replaced by $(\mathbf{I}^{11}, \mathbf{a}_1, \boldsymbol{\Psi}_{\mathbf{i}}^1, \mathcal{E}_{\mathbf{i}}^1, \delta_{\mathbf{i}}^1, \boldsymbol{\Delta}_0^1)$.

REMARK 2. Since $\mathbf{G}_0 = \mathbf{H}_0^\circ = (\mathbf{g}_0', \mathbf{h}_0')' = \mathbf{I}^{-1} \boldsymbol{\Delta}_0$, the Bayes estimator, MLE and hybrid estimator are asymptotically first-order equivalent, normal and efficient.

*Computation consideration.* Although in some cases $(\check{\boldsymbol{\alpha}}_n, \hat{\boldsymbol{\beta}}_n)$ has closed form, the solution of (2.1) generally may not. Denote $G_n(\mathbf{d}, \boldsymbol{\beta}) = \int W(\mathbf{d}(\mathbf{x}_n), \boldsymbol{\alpha}) \times f(\mathbf{x}_n | \boldsymbol{\alpha}, \boldsymbol{\beta}) \pi(\boldsymbol{\alpha}) \, d\boldsymbol{\alpha}$. Since

$$\sup_{\boldsymbol{\beta}} \left[ \inf_{\mathbf{d}} G_n(\mathbf{d}, \boldsymbol{\beta}) \right] \leq \inf_{(\mathbf{d}, \boldsymbol{\beta})} \sup G_n(\mathbf{d}, \boldsymbol{\beta}) \leq \inf_{\mathbf{d}} \left[ \sup_{\boldsymbol{\beta}} G_n(\mathbf{d}, \boldsymbol{\beta}) \right],$$

and the "=" signs do not always hold, so generally

$$\arg \sup_{\boldsymbol{\beta}} \left[ \inf_{\mathbf{d}} G_n(\mathbf{d}, \boldsymbol{\beta}) \right] \neq (\check{\boldsymbol{\alpha}}_n, \hat{\boldsymbol{\beta}}_n) = \arg \inf_{(\mathbf{d}, \boldsymbol{\beta})} \sup G_n(\mathbf{d}, \boldsymbol{\beta}) \neq \arg \inf_{\mathbf{d}} \left[ \sup_{\boldsymbol{\beta}} G_n(\mathbf{d}, \boldsymbol{\beta}) \right].$$

However, if $\arg \inf_{\mathbf{d}} G_n(\mathbf{d}, \boldsymbol{\beta})$ does not depend on $\boldsymbol{\beta}$, then $(\check{\boldsymbol{\alpha}}_n, \hat{\boldsymbol{\beta}}_n) = \sup_{\boldsymbol{\beta}}[\inf_{\mathbf{d}} G_n(\mathbf{d}, \boldsymbol{\beta})]$. Similarly, if $\arg \sup_{\boldsymbol{\beta}} G_n(\mathbf{d}, \boldsymbol{\beta})$ does not depend on $\mathbf{d}$, then $(\check{\boldsymbol{\alpha}}_n, \hat{\boldsymbol{\beta}}_n) = \inf_{\mathbf{d}}[\sup_{\boldsymbol{\beta}} G_n(\mathbf{d}, \boldsymbol{\beta})]$.

When $(\check{\boldsymbol{\alpha}}_n, \hat{\boldsymbol{\beta}}_n)$ is not directly computable, an iterative procedure of it can be formulated by using the Newton–Raphson method.

**3. Implications and examples.** The results obtained in Section 2 imply a new method for the construction of objective prior. They can also can be used to assess high-order behavior of the three estimators and applied to practical problems. Below we discuss these issues with some examples.

*Implication for objective Bayes.* In the objective Bayesian context, the prior is selected by some objective rule instead of subjective choice. Such results include the uniform prior, the reference prior [3] and the noninformative prior [19]. Jeffreys' general rule for such a prior is $\pi(\boldsymbol{\theta}) \propto |\mathbf{I}(\boldsymbol{\theta})|^{1/2}$. Under some regularity conditions and without nuisance parameters, the reference prior coincides with Jeffreys' prior. The coverage matching prior is one that the posterior probability matches the corresponding frequentist probability with high accuracy. Authors [27, 35, 38], among others, studied priors with second-order probability matching. A comprehensive review of this area can be found in [21].

Here we use a similar idea to that of coverage matching to select prior such that the Bayes estimate and MLE match for high-order terms in their expansions. We say that a prior $\pi(\cdot)$ on $\boldsymbol{\theta}$ (or on $\boldsymbol{\alpha}$) is $r$th *order expansion matching prior* for Bayes estimate (or for the hybrid estimator) if the first $r$ terms in its expansion under $\pi(\cdot)$ match those of the MLE, that is, under the notation of the last section,

$$\mathbf{G}_i = \mathbf{H}_i^\circ \qquad [\text{or } (\mathbf{g}_i', \mathbf{h}_i')' = \mathbf{H}_i^\circ], \qquad i = 0, \ldots, r-1.$$

As $\mathbf{G}_0 = \mathbf{H}_0^\circ = (\mathbf{g}_0', \mathbf{h}_0')'$ under the conditions of Theorems 2.3 and 2.4, any prior is automatically 1st-order matching. In the above equations, all quantities involving $\boldsymbol{\theta}_0$ in their definition should be replaced by $\boldsymbol{\theta}$ (or $\boldsymbol{\alpha}$), so these equations are a set of differential equations of order $r-1$ for the prior density as a function of $\boldsymbol{\theta}$ (or $\boldsymbol{\alpha}$). Especially, by the expressions in the Fact, a second-order matching prior $\pi(\boldsymbol{\theta})$ in this sense should satisfy

$$\mathbf{I}^{-1}\boldsymbol{\rho_0} + \mathbf{M_{0,1}} = \mathbf{0} \quad \text{or} \quad \frac{\partial \log \pi(\boldsymbol{\theta})}{\partial \boldsymbol{\theta}} = -\{\boldsymbol{\sigma}(\mathbf{a})\}^{-1} \sum_{|\mathbf{i}|=3} \frac{\mathcal{E}_\mathbf{i}(\boldsymbol{\theta})}{\mathbf{i}!} \boldsymbol{\Psi_i}.$$

To solve $\pi(\cdot)$ from the above partial differential equations, denote $\mathbf{b}(\boldsymbol{\theta}) = -\{\boldsymbol{\sigma}(\mathbf{a})\}^{-1} \sum_{|\mathbf{i}|=3} \frac{\mathcal{E}_\mathbf{i}(\boldsymbol{\theta})}{\mathbf{i}!} \boldsymbol{\Psi_i} = (b_1(\boldsymbol{\theta}), \ldots, b_d(\boldsymbol{\theta}))'$, $\boldsymbol{\theta}_{-k} = (\theta_1, \ldots, \theta_{k-1}, \theta_{k+1}, \ldots, \theta_d)'$, then it is seen that the solution exists if and only if

$$\frac{\partial b_i(\boldsymbol{\theta})}{\partial \theta_j} = \frac{\partial b_j(\boldsymbol{\theta})}{\partial \theta_i} \qquad (i \neq j),$$

which is equivalent to that there are functions $v_k(\boldsymbol{\theta}_{-k})$ $(k=1,\ldots,d)$ such that the following set of equations of indefinite integrals hold

$$\int b_i(\boldsymbol{\theta}) \, d\theta_i + v_i(\boldsymbol{\theta}_{-i}) = \int b_j(\boldsymbol{\theta}) \, d\theta_j + v_j(\boldsymbol{\theta}_{-j}) \qquad (1 \leq i < j \leq d),$$

then for any $k$ $(1 \leq k \leq d)$, up to constant,

$$\log \pi(\boldsymbol{\theta}) = \int b_k(\boldsymbol{\theta}) \, d\theta_k + v_k(\boldsymbol{\theta}_{-k}).$$

Notationally, we denote the solution as

$$\pi(\boldsymbol{\theta}) \propto \exp\left(-\int \{\boldsymbol{\sigma}(\mathbf{a})\}^{-1} \sum_{|\mathbf{i}|=3} \frac{\mathcal{E}_\mathbf{i}(\boldsymbol{\theta})}{\mathbf{i}!} \boldsymbol{\Psi_i} \, d\boldsymbol{\theta}\right).$$

Similarly, for the hybrid estimator, a second-order matching prior $\pi(\cdot)$ on $\boldsymbol{\alpha}$ should satisfy

$$\begin{cases} \mathbf{I}^{11}\boldsymbol{\rho}_0^1 + \mathbf{m}_{0,1} = \mathbf{0}, \\ \mathbf{I}^{21}\boldsymbol{\rho}_0^1 = \mathbf{0} \end{cases} \text{ or } \begin{cases} \dfrac{\partial \log \pi(\boldsymbol{\alpha})}{\partial \boldsymbol{\alpha}} + \{\boldsymbol{\sigma}(\mathbf{a}_1)\}^{-1} \displaystyle\sum_{|\mathbf{i}|=3} \dfrac{\mathcal{E}_{\mathbf{i}}^1(\boldsymbol{\alpha},\boldsymbol{\beta}_1)}{\mathbf{i}!}\boldsymbol{\Psi}_{\mathbf{i}}^1 = \mathbf{0}, \\ \mathbf{I}^{21}(\boldsymbol{\alpha},\boldsymbol{\beta}_1)\dfrac{\partial \log \pi(\boldsymbol{\alpha})}{\partial \boldsymbol{\alpha}} = \mathbf{0}, \end{cases}$$

which is a set of equations in $\partial \pi(\boldsymbol{\alpha})/\partial \boldsymbol{\alpha}$, with $\boldsymbol{\beta}_1$ as a hyper parameter. As in Example 1 below, when $\mathbf{I}^{21} = \mathbf{0}$, $\mathbf{I}^{11}$ is independent on components of $\boldsymbol{\alpha}$, and $\mathbf{a}_1 = 2\mathbf{1}_1$, we will have $\pi(\boldsymbol{\alpha}) \propto |\mathbf{I}^{11}(\boldsymbol{\alpha},\boldsymbol{\beta}_1)|^{1/2}$.

The $\boldsymbol{\Psi}_{\mathbf{i}}$'s ($\boldsymbol{\Psi}_{\mathbf{i}}^1$'s) are $d$-vector functions of the $\sigma(\mathbf{i})$'s which can be found in various sources. Denote $\mathbf{I}^{-1} = (\sigma_{st})$, we have $\sigma((4,0,0)) = 3\sigma_{11}^2$, $\sigma((3,1,0)) = 3\sigma_{11}\sigma_{12}$, $\sigma((2,2,0)) = \sigma_{11}\sigma_{22} + 2\sigma_{12}^2$ and $\sigma((2,1,1)) = \sigma_{11}\sigma_{23} + 2\sigma_{12}\sigma_{13}$.

Like the second-order probability matching prior, the second-order expansion matching prior may not always exist nor be unique. For the former, Mukerjee and Ghosh [27] gave closed-form examples only under some special parameterizations. Below we give several examples of second-order expansion matching priors in natural parametrization.

EXAMPLE 1. When $\mathbf{I} = \{(I_{11}(\theta_1), \ldots, I_{dd}(\theta_d))\}$ is in independent parametric form and $\mathbf{a} = 2\mathbf{1}$, we have $\{\boldsymbol{\sigma}(\mathbf{a})\} = \mathbf{I}^{-1}(\boldsymbol{\theta})$. Assume the conditions for exchange of expectation and differentiation, then for $\mathbf{i} = 3\mathbf{e}_j$, some $j$, $\mathcal{E}_{\mathbf{i}} = -\partial I_{jj}(\theta_j)/\partial \theta_j$, $\boldsymbol{\Psi}_{\mathbf{i}} = 3I_{jj}^{-2}(\theta_j)\mathbf{e}_j$. For $\mathbf{i} \neq 3\mathbf{e}_j$, some $j$, $\mathcal{E}_{\mathbf{i}} = 0$. So $\{\boldsymbol{\sigma}(\mathbf{a})\}^{-1}\sum_{|\mathbf{i}|=3}\frac{\mathcal{E}_{\mathbf{i}}(\boldsymbol{\theta})}{\mathbf{i}!}\boldsymbol{\Psi}_{\mathbf{i}} = -\frac{1}{2}(\frac{\partial I(\theta_1)}{\partial \theta_1}I_{11}^{-1}(\theta_1), \ldots, \frac{\partial I(\theta_d)}{\partial \theta_1}I_{dd}^{-1}(\theta_d))'$. It is easy to see that $v_k(\boldsymbol{\theta}_{-k}) = \sum_{i \neq k}(1/2)\int (\partial I_{ii}(\theta_i)/\partial \theta_i)/I_{ii}(\theta_i)\,d\theta_i$, and the second-order expansion matching prior is

$$\pi(\boldsymbol{\theta}) \propto \prod_{i=1}^{d}\exp\left\{\frac{1}{2}\int \frac{\partial I_{ii}(\theta_i)}{\partial \theta_i}I_{ii}^{-1}(\theta_i)\,d\theta_i\right\} = |\mathbf{I}(\boldsymbol{\theta})|^{1/2},$$

which is Jeffreys' prior.

EXAMPLE 2. Consider the data model $N(\mu,\sigma^2)$ with parameter $\boldsymbol{\theta} = (\mu,\sigma^2)$. We have $\mathbf{I}^{-1} = \{(\sigma^2, 2\sigma^4)\}$. In this case, $\mathcal{E}_{(3,0)} = \mathcal{E}_{(1,2)} = 0$, $\mathcal{E}_{(2,1)} = 1/\sigma^4$, $\mathcal{E}_{(0,3)} = 2/\sigma^6$, $\boldsymbol{\Psi}_{(2,1)} = (\sigma(3,1),\sigma(2,2))' = (0, 2\sigma^6)'$, $\boldsymbol{\Psi}_{(0,3)}(\mathbf{0}) = (\sigma(1,3), \sigma(0,4))' = (0, 12\sigma^8)'$. If set $\mathbf{a} = (2,2)'$, then $\{\boldsymbol{\sigma}(\mathbf{a})\} = \mathbf{I}^{-1}$ and $-\{\boldsymbol{\sigma}(\mathbf{a})\}^{-1} \times \sum_{|\mathbf{i}|=3}\frac{\mathcal{E}_{\mathbf{i}}(\boldsymbol{\theta})}{\mathbf{i}!}\boldsymbol{\Psi}_{\mathbf{i}} = (0, -(5/2)\sigma^{-2})'$. It is easy to check that if we choose $v_1(\sigma^2) = -(5/2)\int (\sigma^2)^{-1}\,d\sigma^2$ and $v_2(\mu) = const.$, then the second-order expansion matching prior is $\pi(\boldsymbol{\theta}) = \pi(\mu,\sigma^2) \propto (\sigma^2)^{-5/2}$. In contrast, Jeffreys' prior in this case is $\pi(\mu,\sigma^2) \propto (\sigma^2)^{-3/2}$.

EXAMPLE 3. Consider the bivariate normal model with parameters $\boldsymbol{\theta} = (\mu_1, \mu_2, \sigma_1^2, \sigma_2^2, \rho)$, with $\rho$ being the correlation coefficient. Suppose we make a hybrid inference with a Bayesian component on $\boldsymbol{\alpha} = (\sigma_1^2, \sigma_2^2, \rho)$ and want a second-order matching prior $\pi(\boldsymbol{\alpha})$. Here we need to replace $(\mathbf{I}^{11}, \mathbf{I}^{21})$ by $(\mathbf{I}^{22}, \mathbf{I}^{12})$ in the partial differential equations. We have $\mathbf{I}^{12} = \mathbf{0}$ and

$$I^{22}(\boldsymbol{\theta}) = I^{22}(\boldsymbol{\alpha}) = \begin{pmatrix} 2\sigma_1^4 & 2\rho^2\sigma_1^2\sigma_2^2 & \rho(1-\rho^2)\sigma_1^2 \\ 2\rho^2\sigma_1^2\sigma_2^2 & 2\sigma_2^4 & \rho(1-\rho^2)\sigma_2^2 \\ \rho(1-\rho^2)\sigma_1^2 & \rho(1-\rho^2)\sigma_2^2 & (1-\rho^2)^2 \end{pmatrix},$$

$\mathcal{E}^2_{(3,0,0)} = \frac{16-7\rho^2}{8(1-\rho^2)\sigma_1^6}$, $\mathcal{E}^2_{(0,3,0)} = \frac{16-7\rho^2}{8(1-\rho^2)\sigma_2^6}$, $\mathcal{E}^2_{(0,0,3)} = \frac{10\rho^5+8\rho^3-12\rho}{(1-\rho^2)^4}$, $\mathcal{E}^2_{(2,1,0)} = -\frac{3\rho^2}{8(1-\rho^2)\sigma_1^4\sigma_2^2}$, $\mathcal{E}^2_{(1,2,0)} = -\frac{3\rho^2}{8(1-\rho^2)\sigma_1^2\sigma_2^4}$, $\mathcal{E}^2_{(2,0,1)} = \frac{3\rho^3-5\rho}{4(1-\rho^2)^2\sigma_1^4}$, $\mathcal{E}^2_{(0,2,1)} = \frac{3\rho^3-5\rho}{4(1-\rho^2)^2\sigma_2^4}$, $\mathcal{E}^2_{(1,0,2)} = \frac{(1+2\rho^2)\sigma_1^2-\rho^2(3+\rho^2)}{(1-\rho^2)^3\sigma_1^2}$, $\mathcal{E}^2_{(0,1,2)} = \frac{(1+2\rho^2)\sigma_2^2-\rho^2(3+\rho^2)}{(1-\rho^2)^3\sigma_2^2}$, $\mathcal{E}^2_{(1,1,1)} = \frac{\rho(1+\rho^2)}{4(1-\rho^2)^2\sigma_1^2\sigma_2^2}$; $\boldsymbol{\Psi}^2_{(3,0,0)} = (12\sigma_1^8, 12\rho^2\sigma_1^6\sigma_2^2, 6\rho(1-\rho^2)\sigma_1^6)'$, $\boldsymbol{\Psi}^2_{(0,3,0)} = (12\rho^2\sigma_1^2 \times \sigma_2^6, 12\sigma_2^8, 6\rho(1-\rho^2)\sigma_2^6)'$, $\boldsymbol{\Psi}^2_{(0,0,3)} = (3\rho(1-\rho^2)^3\sigma_1^2, 3\rho(1-\rho^2)^3\sigma_2^2, 3(1-\rho^2)^4)'$, $\boldsymbol{\Psi}^2_{(2,1,0)} = (12\rho^2\sigma_1^6\sigma_2^2, 4(1+2\rho^4)\sigma_1^4\sigma_2^4, 2\rho(1-\rho^2)(1+2\rho^2)\sigma_1^4\sigma_2^2)'$, $\boldsymbol{\Psi}^2_{(1,2,0)} = (4(1+2\rho^4)\sigma_1^4\sigma_2^4, 12\rho^2\sigma_1^2\sigma_2^6, 2\rho(1-\rho^2)\sigma_2^2(\sigma_1^4+2\rho^2\sigma_2^4))'$, $\boldsymbol{\Psi}^2_{(2,0,1)} = (6\rho(1-\rho^2)\sigma_1^6, 2\rho(1-\rho^2)(1+2\rho^2)\sigma_1^4\sigma_2^2, 2(1+\rho^2)(1-\rho^2)^2\sigma_1^4)'$, $\boldsymbol{\Psi}^2_{(0,2,1)} = (2\rho(1-\rho^2)\sigma_1^2(\sigma_1^4+2\rho^2\sigma_2^4), 6\rho(1-\rho^2)\sigma_2^6, 2(1+\rho^2)(1-\rho^2)^2\sigma_1^4)'$, $\boldsymbol{\Psi}^2_{(1,0,2)} = ((1+2\rho^2)(1-\rho^2)^2\sigma_1^4, 4\rho^2(1-\rho^2)^2\sigma_1^2\sigma_2^2, 3\rho(1-\rho^2)^3\sigma_1^4)'$, $\boldsymbol{\Psi}^2_{(0,1,2)} = (4\rho^2(1-\rho^2)^2\sigma_1^2\sigma_2^2, 2(1+\rho^2)(1-\rho^2)^2\sigma_2^4, 3\rho(1-\rho^2)^3\sigma_2^4)'$, $\boldsymbol{\Psi}^2_{(1,1,1)} = (2\rho(1-\rho^2)(1+2\rho^2)\sigma_1^4\sigma_2^2, 2\rho(1-\rho^2)(1+2\rho^2)\sigma_1^2\sigma_2^4, 4\rho^2(1-\rho^2)^2\sigma_1^2\sigma_2^2)'$. Take $\mathbf{a}_2 = (2,2,2)'$, then $\{\boldsymbol{\sigma}(\mathbf{a}_2)\}^{-1} = \{(2\sigma_1^4, 2\sigma_2^4, (1-\rho^2)^2)\}^{-1}$. In this case, the second-order expansion matching prior $\pi(\boldsymbol{\alpha})$ cannot be evaluated in closed form, and numerical method, such as in [18], is needed.

Result for second-order expansion matching prior $\pi(\boldsymbol{\theta}) = \pi(\mu_1, \mu_2, \sigma_1^2, \sigma_2^2, \rho)$ can also be obtained similarly. In this case for $|\mathbf{i}| = 3$, there are 35 $\mathcal{E}_\mathbf{i}$'s and $\boldsymbol{\Psi}_\mathbf{i}$'s each, the prior will not be evaluated in closed form and numerical method is needed.

*Evaluation of high-order behavior.* Note that in the $\mathbf{H}_i$'s the $\boldsymbol{\Delta}_\mathbf{i}$'s are asymptotic normal random vectors, the $\boldsymbol{\rho}_\mathbf{i}$'s and $\mathbf{E}_\mathbf{i}$'s are constant vectors and in the $\mathbf{G}_i$'s the $N_\mathbf{i}$'s are random variables determined by the $\mathbf{H}_i$'s. Thus, asymptotically, $\mathbf{H}_i$, $\mathbf{G}_i$, $\mathbf{h}_i$ and $\mathbf{g}_i$ converge in distribution to multivariate polynomials in normal vectors of degree $i$. For the MLE, $\mathbf{H}_i^\circ$ is an $i$th form of normal vectors. Let $\tilde{\mathbf{H}}_i$, $\tilde{\mathbf{H}}_i^\circ$, $\tilde{\mathbf{G}}_i$, $\tilde{\mathbf{g}}_i$ and $\tilde{\mathbf{h}}_i$ be the weak limits of $\mathbf{H}_i$, $\mathbf{H}_i^\circ$, $\mathbf{G}_i$, $\mathbf{g}_i$ and $\mathbf{h}_i$ $(i = 0, \ldots, k-1)$. The first-order terms in the expansions often have mean zero, so their asymptotic behaviors can be characterized by their asymptotic variances. But high-order terms generally have nonzero mean, so using asymptotic variances alone as a criterion to evaluate their behavior is inappropriate. So we consider an asymptotic mean (bias) and variance combined criterion. We say that $\hat{\boldsymbol{\theta}}_n$ is $r$th order preferred over $\check{\boldsymbol{\theta}}_n$, if $\|E_{\boldsymbol{\theta}_0}(\tilde{\mathbf{H}}_i)\| + \|\operatorname{Cov}_{\boldsymbol{\theta}_0}(\tilde{\mathbf{H}}_i)\| = \|E_{\boldsymbol{\theta}_0}(\tilde{\mathbf{G}}_i)\| + \|\operatorname{Cov}_{\boldsymbol{\theta}_0}(\tilde{\mathbf{G}}_i)\|$ $(i = 0, \ldots, r-2)$, and $\|E_{\boldsymbol{\theta}_0}(\tilde{\mathbf{H}}_{r-1})\| + \|\operatorname{Cov}_{\boldsymbol{\theta}_0}(\tilde{\mathbf{H}}_{r-1})\| < \|E_{\boldsymbol{\theta}_0}(\tilde{\mathbf{G}}_{r-1})\| + \|\operatorname{Cov}_{\boldsymbol{\theta}_0}(\tilde{\mathbf{G}}_{r-1})\|$, and vice-versa.

From the Fact, we see that $\mathbf{G}_0 = \mathbf{H}_0^\circ$ (hence $\tilde{\mathbf{G}}_0 = \tilde{\mathbf{H}}_0^\circ$). Also $\mathbf{G}_1 = \mathbf{H}_1 + \mathbf{M}_{\mathbf{0},1} = \mathbf{H}_1^\circ + \mathbf{I}^{-1}\boldsymbol{\rho}_\mathbf{0} + \mathbf{M}_{\mathbf{0},1}$ with $\mathbf{M}_{\mathbf{0},1} = (\{\boldsymbol{\sigma}(\mathbf{a})\}\mathbf{I})^{-1}\sum_{|\mathbf{i}|=3}\boldsymbol{\Psi}_\mathbf{i}\mathcal{E}_\mathbf{i}/\mathbf{i}!$. Note $\mathbf{H}_1^\circ$ is a random vector and $\mathbf{G}_1$ is $\mathbf{H}_1^\circ$ plus a constant vector $\mathbf{I}^{-1}\boldsymbol{\rho}_\mathbf{0} + \mathbf{M}_{\mathbf{0},1}$. Similarly, for the second-order term of hybrid estimator, its Bayesian components are those of the MLE plus a constant vector. Hence for the second-order behavior, we only need to consider the expected asymptotic bias (EAB), and $\hat{\boldsymbol{\theta}}_n$ is second order preferred over $\check{\boldsymbol{\theta}}_n$, if

$$\|E_{\boldsymbol{\theta}_0}(\tilde{\mathbf{H}}_1)\| < \|E_{\boldsymbol{\theta}_0}(\tilde{\mathbf{G}}_1)\|$$

and vice versa. We have:

PROPOSITION. *Let* $\mathbf{D}_j = E_{\boldsymbol{\theta}_0}[\mathbf{L}_{\mathbf{e}_j}(\mathbf{X}|\boldsymbol{\theta}_0)\mathbf{L}'_\mathbf{0}(\mathbf{X}|\boldsymbol{\theta}_0)]$, *and* $_j\mathbf{I}^{-1}$ *and* $\mathbf{I}_j^{-1}$ *be the* $j$th *row and* $j$th *column of* $\mathbf{I}^{-1}$, *then*

(i) $\quad E_{\boldsymbol{\theta}_0}(\tilde{\mathbf{H}}_1^\circ) = \mathbf{I}^{-1}\left(\sum_{j=1}^d \mathbf{D}_j \mathbf{I}_j^{-1} + \sum_{i,j=1}^d \mathbf{E}_{\mathbf{e}_i+\mathbf{e}_j} {}_i\mathbf{I}^{-1}\mathbf{I}^{-1}\mathbf{I}_j^{-1}/(\mathbf{e}_i + \mathbf{e}_j)!\right),$

(ii) $\quad E_{\boldsymbol{\theta}_0}(\tilde{\mathbf{G}}_1) = E_{\boldsymbol{\theta}_0}(\tilde{\mathbf{H}}_1^\circ) + \mathbf{I}^{-1}\boldsymbol{\rho}_\mathbf{0} + \mathbf{M}_{\mathbf{0},1},$

(iii) $\quad E_{\boldsymbol{\theta}_0}\begin{pmatrix}\tilde{\mathbf{g}}_1 \\ \tilde{\mathbf{h}}_1\end{pmatrix} = E_{\boldsymbol{\theta}_0}(\tilde{\mathbf{H}}_1^\circ) + \begin{pmatrix}\mathbf{I}^{11}\boldsymbol{\rho}_\mathbf{0}^1 + (\{\boldsymbol{\sigma}(\mathbf{a}_1)\mathbf{I}^{11}\})^{-1}\sum_{|\mathbf{i}|=3}\boldsymbol{\Psi}_\mathbf{i}^1\mathcal{E}_\mathbf{i}^1/\mathbf{i}! \\ \mathbf{I}^{21}\boldsymbol{\rho}_\mathbf{0}^1\end{pmatrix}.$

By this proposition we are able to evaluate which estimator has second-order advantage under each specification of the likelihood model, prior and the loss, by the criterion of EAB.

EXAMPLE 4. We want to evaluate the overall and small-sample advantage of an estimator. We first evaluate the overall advantage for the model in Example 2 in four cases: (a) full MLE $(\hat{\mu}_n, \hat{\sigma}_n^2)$; (b) full Bayesian $(\check{\mu}_n, \check{\sigma}_n^2)$; (c) hybrid Bayes–MLE $(\check{\mu}_n, \hat{\sigma}_n^2)$; and (d) hybrid MLE–Bayes $(\hat{\mu}_n, \check{\sigma}_n^2)$. We know that the four estimates are first order equivalent, we are to study their second-order behavior by the EAB criterion. Here,

$$\mathbf{I}^{-1} = \mathbf{I}^{-1}(\boldsymbol{\theta}_0) = \begin{pmatrix}\sigma_0^2 & 0 \\ 0 & 2\sigma_0^4\end{pmatrix},$$

$$\mathbf{D}_1 = \begin{pmatrix}0 & 0 \\ \frac{1}{\sigma_0^4} & 0\end{pmatrix}, \qquad \mathbf{D}_2 = \begin{pmatrix}-\frac{1}{\sigma_0^4} & 0 \\ 0 & -\frac{1}{\sigma_0^6}\end{pmatrix},$$

$\mathbf{E}_{\mathbf{e}_1+\mathbf{e}_1} = (0, 1/\sigma_0^4)'$, $\mathbf{E}_{\mathbf{e}_1+\mathbf{e}_2} = \mathbf{E}_{\mathbf{e}_2+\mathbf{e}_1} = (1/\sigma_0^4, 0)'$, $\mathbf{E}_{\mathbf{e}_2+\mathbf{e}_2} = (0, 2/\sigma_0^6)'$.

For (a), we have $E_{\boldsymbol{\theta}_0}(\tilde{\mathbf{H}}_1^\circ) = (0, -2\sigma_0^2 + \sigma_0^6 + 16\sigma_0^{10})'$.

For (b), if we use the prior $\pi(\boldsymbol{\theta}) = \pi(\mu)\pi(\sigma^2)$, $\mu \sim N(\mu_1, \sigma_1^2)$ with $(\mu_1, \sigma_1^2)$ known, $\pi(\sigma^2) = \lambda_1 e^{-\lambda_1 \sigma^2}$, $\lambda_1 > 0$ known. Then, $\boldsymbol{\rho_0} = -((\mu_0 - \mu_1)/\sigma_1^2, \lambda_1)'$. If we use the loss $W(\mathbf{d}, \boldsymbol{\theta}) = (d_1 - \mu)^2 + (d_2 - \sigma^2)^2$, then $\{\boldsymbol{\sigma}(\mathbf{a})\} = \mathbf{I}^{-1}$, and for $|\mathbf{i}| = 3$, $\boldsymbol{\Psi_i} = (\sigma(\mathbf{i} + \mathbf{e}_1), \sigma(\mathbf{i} + \mathbf{e}_2))'$. We have $\mathcal{E}_{(3,0)} = \mathcal{E}_{(1,2)} = 0$, $\mathcal{E}_{(2,1)} = \sigma_0^{-4}$, $\mathcal{E}_{(0,3)} = 2\sigma_0^{-6}$, $\boldsymbol{\Psi}_{(0,3)} = (\sigma(1,3), \sigma(0,4))' = (0, 12\sigma_0^8)'$. Thus, $\mathbf{M}_{\mathbf{0},1} = (0, 2\sigma_0^2)'$, and $E_{\boldsymbol{\theta}_0}(\tilde{\mathbf{G}}_1^\circ) = ((\mu_1 - \mu_0)\sigma_0^2/\sigma_1^2, -2\sigma_0^4\lambda_1 + \sigma_0^6 + 16\sigma_0^{10})'$. Note the MLE is in closed form and the Bayesian is not. If one has good prior knowledge of $\boldsymbol{\theta}$, that is, $\mu_1 \approx \mu_0$ and $\lambda_1 \approx \sigma_0^{-2}$, then $E_{\boldsymbol{\theta}_0}(\tilde{\mathbf{G}}_1^\circ) \approx E_{\boldsymbol{\theta}_0}(\tilde{\mathbf{H}}_1^\circ)$, so the Bayes estimator and MLE have similar second-order behavior, while the former has small-sample advantage and the latter is computationally preferable. One can even choose $\lambda_1 \approx (\sigma_0^2 + 16\sigma_0^4)/2$, so $E_{\boldsymbol{\theta}_0}(\tilde{\mathbf{G}}_1^\circ) \approx \mathbf{0}$, thus the Bayesian has smaller asymptotic second-order bias than the MLE, but also lost its small-sample advantage on the estimate of $\sigma^2$.

For (c), let $\pi(\mu)$ be as in (b), $W(d_1, \mu) = (d_1 - \mu)^2$. Then, $\rho_0^1 = (\mu_1 - \mu_0)/\sigma_1^2$, $\{\boldsymbol{\sigma}(\mathbf{a}_1)\} = \mathbf{I}^{11}$, and for $\mathbf{i} = i = 3$, $\boldsymbol{\Psi}_i^1 = 3\sigma_0^4$ and $\mathcal{E}_i^1 = 0$. Thus, $E_{\boldsymbol{\theta}_0}(\tilde{\mathbf{g}}_1, \tilde{\mathbf{h}}_1) = ((\mu_1 - \mu_0)\sigma_0^2/\sigma_1^2, -2\sigma_0^2 + \sigma_0^6 + 16\sigma_0^{10})$. As in (b), if one has good informative prior on $\mu$, then the hybrid estimator can have small-sample advantage over the MLE, and they are compatible in second-order behavior and computation. If we do not have sound information on $\sigma^2$, (c) often has better second-order property than (b).

For (d), let $\pi(\sigma^2)$ as in (b), $W(d_2, \sigma^2) = (d_2 - \sigma^2)^2$. Then, $\rho_0^2 = -\lambda_1$, $\{\boldsymbol{\sigma}(\mathbf{a}_2)\} = \mathbf{I}^{22}$, and for $\mathbf{i} = i = 3$, $\boldsymbol{\Psi}_i^2 = 12\sigma_0^8$ and $\mathcal{E}_i^2 = 0$. Thus $E_{\boldsymbol{\theta}_0}(\tilde{\mathbf{h}}_1, \tilde{\mathbf{g}}_1) = (0, -2\lambda_1\sigma_0^4 + \sigma_0^6 + 16\sigma_0^{10})$. In this case $\hat{\mu}_n$ is in closed form and $\check{\sigma}_n^2$ is not. If we infer $\tau = 1/\sigma^2$ and use a Gamma prior on it, then $(\hat{\mu}_n, \check{\tau}_n)$ has closed form, while the full Bayes estimator (b) is not. This has practical meaning when one seeks high-order accuracy in addition to computational advantage. As (b) is usually obtained by numerical approximation methods, such as Markov chain Monte Carlo, the inaccuracy of these methods is not easy to assess.

Next, we discuss small-sample advantage. Consider the $\mathbf{x}_i$'s are i.i.d. multivariate normal $N(\boldsymbol{\mu}, \boldsymbol{\Omega})$. Suppose we have good prior knowledge on $\boldsymbol{\Omega}$ in the form of a Wishart prior $\pi(\boldsymbol{\Omega})$, but relative ignorance about $\boldsymbol{\mu}$. So we use MLE for $\boldsymbol{\mu}$ and jointly, a Bayes estimate for $\boldsymbol{\Omega}$. Since in this case $\arg\sup_{\boldsymbol{\mu}} G_n(\mathbf{d}, \boldsymbol{\mu}) = \hat{\boldsymbol{\mu}} = \overline{\mathbf{x}}$ does not depend on $\mathbf{d}$, $\check{\boldsymbol{\Omega}}$ is just the Bayes solution given $\hat{\boldsymbol{\mu}} = \overline{\mathbf{x}}$, which can be evaluated numerically.

On the other hand, suppose we have good information about $\boldsymbol{\mu}$ summarized by the prior $N(\boldsymbol{\mu}_1, \boldsymbol{\Omega})$, but not enough knowledge on $\boldsymbol{\Omega}$. Note here the prior has the same unknown variance matrix as that in the data distribution, so that the Bayesian part has a closed-form solution. It can be checked that Theorem 2.1 still applies to this case. We want to use the good prior experience for $\boldsymbol{\mu}$, but a full Bayesian estimate for $(\boldsymbol{\mu}, \boldsymbol{\Omega})$ is not easy, so we estimate $\boldsymbol{\Omega}$ by the MLE. It is easy to see that for given $\boldsymbol{\Omega}$, the posterior on $\boldsymbol{\mu}$ is $N(n\overline{\mathbf{x}}/(n+1) + \boldsymbol{\mu}_1/(n+1), \boldsymbol{\Omega}/(n+1))$. With loss to be either absolute error, squared error or 0–1 error on $\boldsymbol{\mu}$, $\check{\boldsymbol{\mu}}_n$ is either the posterior median, mean or mode, which are all the same in this case and is independent of $\boldsymbol{\Omega}$. So we have $(\check{\boldsymbol{\mu}}_n, \hat{\boldsymbol{\Omega}}_n) = (\frac{n}{n+1}\overline{\mathbf{x}} + \frac{\boldsymbol{\mu}_1}{n+1}, \frac{1}{n+1}\sum_{i=1}^n (\mathbf{x}_i - \check{\boldsymbol{\mu}}_n)'(\mathbf{x}_i - \check{\boldsymbol{\mu}}_n) + \frac{1}{n+1}(\boldsymbol{\mu}_1 - \check{\boldsymbol{\mu}}_n)'(\boldsymbol{\mu}_1 - \check{\boldsymbol{\mu}}_n))$, which is given in closed form, while the full Bayes estimator $(\check{\boldsymbol{\mu}}_n, \check{\boldsymbol{\Omega}}_n)$ is not.

EXAMPLE 5. Using existing results, we give an application of the hybrid estimator in which neither the full MLE nor full Bayes estimator works. Let $x|\alpha$ be the model in [13] with distribution $P(A|\alpha)$ and density function

$$f(x|\alpha) = (1-\alpha)\frac{1}{\delta(\alpha)}f_0\left(\frac{x-\alpha}{\delta(\alpha)}\right) + \alpha f_1(x), \qquad \alpha \in [0,1],$$

where $f_0(x) = (1-|x|)\chi_{[-1,1]}(x)$, $f_1(x) = \chi_{[-1,1]}(x)/2$ and $\delta(\cdot)$ satisfies $\delta(0) = 1$, $0 < \delta(\alpha) \le 1-\alpha$ and $\delta(\alpha) \to 0$ as $\alpha \to 1$. Ferguson [13] shows that the MLE $\hat{\alpha}_n$ of $\alpha$ is not consistent; $\hat{\alpha}_n \to 1$ (a.s.) no matter what the true parameter $\alpha_0$ is if $\delta(\alpha) \to 0$ fast enough, in particular if $\delta(\alpha) = (1-\alpha)\exp(-(1-\alpha)^{-c}+1)$ with $c > 2$. On the other hand, it is easy to see that for this model Doob's general conditions (as stated in [31]) for the consistency of Bayes estimator are satisfied as follows: (1) The measurable spaces $\{\mathcal{X}, \mathcal{B}\}$ of $x$ and $\{\Lambda, \mathcal{U}\}$ of $\alpha$ are isomorphic to Borel field in a complete separable metric space; (2) For every $A \in \mathcal{B}$, $P(A|\cdot)$ is $\mathcal{U}$-measurable; (3) If $\alpha_1 \ne \alpha_2$ there exists a set $A \in \mathcal{B}$ such that $P(A|\alpha_1) \ne P(A|\alpha_2)$; (4) The prior $\pi(\cdot)$ has finite second moment. Then, the Bayes estimator $\check{\alpha}_n$ under quadratic loss is strongly consistent a.e. $(\pi)$.

On the other hand, let $y_1, \ldots, y_n$ i.i.d $y|\beta \sim U[0,1]\chi(\beta=1) + U[0, 2/\beta)\chi(1 < \beta < 2)$ be the model in Example 2 in [31] with the prior $\pi(\cdot)$ to be the Lebesgue measure on the Borel sets on $[1,2)$ and $U[0,a]$ be the uniform distribution on $[0,a]$. Denote $y_{(n)} = \max_{1 \le i \le n} y_i$, then the MLE $\hat{\beta}_n$ and Bayes estimator $\check{\beta}_n$ (under squared error loss) of $\beta$ are

$$\hat{\beta}_n = \chi(y_{(n)} \le 1) + \frac{2}{y_{(n)}}\chi(y_{(n)} > 1),$$

$$\check{\beta}_n = \frac{n+1}{n+2}\frac{2^{n+2}-1}{2^{n+1}-1}\chi(y_{(n)} \le 1) + \frac{\int_1^{2/y_{(n)}} \theta^{n+1}\,d\theta}{\int_1^{2/y_{(n)}} \theta^n\,d\theta}\chi(y_{(n)} > 1).$$

Schwartz [31] showed that the MLE is consistent while the Bayesian is not when $\beta = 1$ (although under some special prior, $\check{\beta}_n$ can be consistent).

Now let $x$ and $y$ be independent, $(x_1, y_1), \ldots, (x_n, y_n)$ be an i.i.d. sample of $(x, y)$ and the parameter be $\boldsymbol{\theta} = (\alpha, \beta)$. Then, neither the full Bayes estimator $(\check{\alpha}_n, \check{\beta}_n)$ nor the full MLE $(\hat{\alpha}_n, \hat{\beta}_n)$ of $\boldsymbol{\theta}$ is consistent while the hybrid estimator $(\check{\alpha}_n, \hat{\beta}_n)$ is.

EXAMPLE 6. As a last application, let's consider the problem mentioned in the Introduction. The data follows a mixture model $f(\mathbf{x}|\boldsymbol{\alpha},\boldsymbol{\beta}) = \sum_{j=1}^{k} \gamma_j \phi(\mathbf{x}|\boldsymbol{\alpha}_j, \boldsymbol{\Omega}_j)$, where $\phi(\cdot|\boldsymbol{\alpha}, \boldsymbol{\Omega})$ is the density of $N(\boldsymbol{\alpha}, \boldsymbol{\Omega})$. Assume that we have good knowledge on $\boldsymbol{\alpha} = (\boldsymbol{\alpha}_1, \ldots, \boldsymbol{\alpha}_k)$, as summarized by the prior density $\pi_j(\boldsymbol{\alpha}_j) \sim N(\boldsymbol{\alpha}_{j0}, \boldsymbol{\Omega}_{j0})$, $\pi(\boldsymbol{\alpha}) = \prod_{j=1}^{k} \pi(\boldsymbol{\alpha}_j)$, but not enough experience for the parameter $\boldsymbol{\beta} = (\boldsymbol{\beta}_1, \ldots, \boldsymbol{\beta}_k)'$, $\boldsymbol{\beta}_j = (\gamma_j, \boldsymbol{\Omega}_j)$ $(j = 1, \ldots, k)$. So we use a hybrid estimate with a Bayesian components on $\boldsymbol{\alpha}$ and the MLE on $\boldsymbol{\beta}$. To estimate the parameters in a mixture model, often it is more convenient if we use a complete data model. For this, let $I_{ij}$ be the membership indicator of $\mathbf{x}_i$, that is, $I_{ij} = 1$ if $\mathbf{x}_i$ is from the $j$th subdistribution, and $I_{ij} = 1$ and 0 otherwise. Let $I_i = (I_{i1}, \ldots, I_{ik})$, $\mathbf{y}_i = (\mathbf{x}_i, I_i)$ and $\mathbf{y}^n = (\mathbf{y}_1, \ldots, \mathbf{y}_n)$. Treating $I_1, \ldots, I_n$ as missing data, given the "complete data" $\mathbf{y}^n$ and $\boldsymbol{\beta}$, the posterior on $\boldsymbol{\alpha}$ is

$$\pi(\boldsymbol{\alpha}|\mathbf{y}^n, \boldsymbol{\beta}) \propto \pi(\boldsymbol{\alpha}) \prod_{i=1}^{n} \prod_{j=1}^{k} (\gamma_j \phi(\mathbf{x}_i|\boldsymbol{\alpha}_j, \boldsymbol{\Omega}_j))^{I_{ij}} := b(\mathbf{y}^n|\boldsymbol{\theta}),$$

and the corresponding logarithm is

$$l(\boldsymbol{\alpha}, \boldsymbol{\beta}|\mathbf{y}^n) = \sum_{i=1}^{n} \sum_{j=1}^{k} I_{ij}(\log \gamma_j + \log \phi(\mathbf{x}_i|\boldsymbol{\alpha}_j, \boldsymbol{\Omega}_j)) + \log \pi(\boldsymbol{\alpha}).$$

Using the 0–1 loss on $\boldsymbol{\alpha}$, its Bayesian solution is the posterior mode, so we are to maximize $l(\boldsymbol{\alpha}, \boldsymbol{\beta}|\mathbf{y}^n)$ over $(\boldsymbol{\alpha}, \boldsymbol{\beta})$. As typical, this leads to an EM algorithm. However, different from the common EM algorithm, here $l(\boldsymbol{\alpha}, \boldsymbol{\beta}|\mathbf{y}^n)$ is not a proper log-likelihood due to the extra term $\log \pi(\boldsymbol{\alpha})$. If we define $Q(\boldsymbol{\theta}'|\boldsymbol{\theta}) = E[\log b(\mathbf{y}^n|\boldsymbol{\theta}')|\mathbf{x}^n, \boldsymbol{\theta}]$, $H(\boldsymbol{\theta}'|\boldsymbol{\theta}) = E[\log g(\mathbf{y}^n|\boldsymbol{\theta}')|\mathbf{x}^n, \boldsymbol{\theta}]$, where $g(\mathbf{y}^n|\boldsymbol{\theta}) = b(\mathbf{y}^n|\boldsymbol{\theta})/a(\mathbf{x}^n|\boldsymbol{\theta})$ and $a(\mathbf{x}^n|\boldsymbol{\theta}) = \pi(\boldsymbol{\alpha}) \prod_{i=1}^{n} f(\mathbf{x}_i|\boldsymbol{\theta})$. Then $l(\boldsymbol{\theta}|\mathbf{y}^n) = Q(\boldsymbol{\theta}'|\boldsymbol{\theta}) - H(\boldsymbol{\theta}'|\boldsymbol{\theta})$. It is seen that $g(\mathbf{y}^n|\boldsymbol{\theta})$ is just the conditional density of $\mathbf{y}^n$ given $\mathbf{x}^n$, thus $Q(\cdot|\cdot)$ and $H(\cdot|\cdot)$ here play the same roles as they do in the standard EM algorithm [10], and $\boldsymbol{\theta}_n = (\check{\boldsymbol{\alpha}}_n, \hat{\boldsymbol{\beta}}_n)$ can be evaluated in closed form at each iteration, the details omitted here.

Below we give simulation results to compare the performances of the MLE $\hat{\boldsymbol{\theta}}_n$, Bayes estimator $\check{\boldsymbol{\theta}}_n$ and the hybrid estimator $\boldsymbol{\theta}_n$ for this example. We take $k = 3$, $\mathbf{x}_i$'s to be 1-dimensional and $\boldsymbol{\Omega}_j = \sigma_j^2$. The parameter vector in the model is now $\boldsymbol{\theta} = (\gamma_1, \gamma_2, \gamma_3, \alpha_1, \alpha_2, \alpha_3, \sigma_1^2, \sigma_2^2, \sigma_3^2)$. For the Bayes estimator, the prior is taken as $\pi(\boldsymbol{\theta}) = \pi(\boldsymbol{\alpha})\pi(\boldsymbol{\gamma})\pi(\boldsymbol{\sigma}^2)$. Since we do not have good knowledge on $(\boldsymbol{\gamma}, \boldsymbol{\sigma}^2)$, except that $0 < \gamma_j < 1$, $\sum_{j=1}^{3} \gamma_j = 1$, and $0 < \sigma_j^2 < 3$ $(j = 1, 2, 3)$. We use a noninformative prior on them, that is, $\pi(\boldsymbol{\gamma}) = \pi(\gamma_1)\pi(\gamma_2|\gamma_1)$ with $\pi(\gamma_1) \sim U(0, 2/3)$, $\pi(\gamma_2|\gamma_1) \sim U(0, 1 - \gamma_1)$, $\gamma_3 = 1 - \gamma_1 - \gamma_2$ and $\pi(\boldsymbol{\sigma}^2) \sim U[0, 3]^3$. To distinguish the Bayesian estimator from the hybrid estimate, we use the squared error loss for the former, which has no closed form, and we use Markov chain Monte Carlo sampling

18    A. YUAN

18    A. YUAN

Table 1
*Simulation results for three estimators*

| $n$ | $\theta_0$ | 0.190 | 0.540 | 0.270 | $-0.850$ | 0.220 | 1.350 | 0.450 | 0.200 | 0.860 |
|---|---|---|---|---|---|---|---|---|---|---|
| 100 | $\hat{\theta}_n$ | 0.281 | 0.463 | 0.255 | $-0.783$ | 0.238 | 1.803 | 0.401 | 0.139 | 0.480 |
|  |  | (0.163) | (0.130) | (0.184) | (0.064) | (0.134) | (0.059) | (0.289) | (0.483) | (0.164) |
|  | $\check{\theta}_n$ | 0.280 | 0.386 | 0.334 | $-0.649$ | 0.204 | 1.490 | 0.539 | 0.177 | 0.722 |
|  |  | (0.158) | (0.129) | (0.146) | (0.055) | (0.089) | (0.051) | (0.188) | (0.345) | (0.093) |
|  | $\theta_n$ | 0.277 | 0.427 | 0.296 | $-0.819$ | 0.210 | 1.637 | 0.348 | 0.104 | 0.575 |
|  |  | (0.168) | (0.135) | (0.169) | (0.071) | (0.149) | (0.059) | (0.370) | (0.649) | (0.131) |
| 300 | $\hat{\theta}_n$ | 0.201 | 0.531 | 0.267 | $-0.905$ | 0.239 | 1.371 | 0.454 | 0.169 | 1.049 |
|  |  | (0.107) | (0.069) | (0.090) | (0.027) | (0.076) | (0.022) | (0.132) | (0.204) | (0.062) |
|  | $\check{\theta}_n$ | 0.173 | 0.563 | 0.264 | $-1.005$ | 0.232 | 1.362 | 0.440 | 0.198 | 1.181 |
|  |  | (0.115) | (0.068) | (0.089) | (0.024) | (0.072) | (0.020) | (0.141) | (0.170) | (0.050) |
|  | $\theta_n$ | 0.201 | 0.531 | 0.268 | $-0.900$ | 0.237 | 1.373 | 0.447 | 0.167 | 1.021 |
|  |  | (0.108) | (0.069) | (0.090) | (0.028) | (0.077) | (0.022) | (0.140) | (0.209) | (0.065) |
| 1000 | $\hat{\theta}_n$ | 0.201 | 0.533 | 0.266 | $-0.797$ | 0.227 | 1.323 | 0.420 | 0.202 | 0.824 |
|  |  | (0.055) | (0.037) | (0.049) | (0.017) | (0.037) | (0.014) | (0.095) | (0.089) | (0.042) |
|  | $\check{\theta}_n$ | 0.185 | 0.556 | 0.259 | $-0.885$ | 0.215 | 1.314 | 0.405 | 0.220 | 0.883 |
|  |  | (0.059) | (0.036) | (0.050) | (0.016) | (0.038) | (0.013) | (0.104) | (0.080) | (0.037) |
|  | $\theta_n$ | 0.202 | 0.531 | 0.267 | $-0.801$ | 0.227 | 1.326 | 0.416 | 0.199 | 0.811 |
|  |  | (0.055) | (0.037) | (0.049) | (0.017) | (0.038) | (0.014) | (0.097) | (0.091) | (0.043) |

to compute it. The results are given in Table 1, in which $\theta_0$ is the true parameter value and the estimated standard deviations are in brackets. We see that when sample size is relatively small ($n = 100$), the hybrid estimator has better performance on $\alpha$, which may be due to the good knowledge on it. The Bayes estimator behaves better only on $\alpha_3$. As sample size increases, the performances of the three estimator are close, as anticipated.

*Some pros and cons of each method.* Here we discuss some advantages and disadvantages of the Bayesian and frequentist method in parametric inference. These known facts can help us in practice in the selection of the method to use. Our list is far from complete.

*Unbiasedness consideration.* It is known [6, 29] that there are essentially no unbiased Bayes procedures.

*Small-sample advantage.* When good prior knowledge about parameters is available, Bayesian estimate often has better small-sample advantage than the frequentist's, due to the information in the prior.

*High-order behavior.* Since the Bayes estimator, MLE and the hybrid estimator have the same first-order performance, if we want higher standards to select among them, usually their second-order terms will be evaluated, such as in Example 4.

*Prior selection.* In some cases we don't have sufficient knowledge for the prior on part of the parameters. Although one may use a noninformative prior [21] on these parameter components, so that a full Bayesian analysis can be performed, this often pays the price of small-sample bias and computational complexity.

*Feasibility.* Sometimes it is difficult to implement a full Bayesian or a full frequentist's analysis on all the parameters of interest, but relatively easy for parts of the parameters by one of the methods. For example, in a multi-parametric model [14], some of the parameters are the change points, the model is nondifferentiable at these points and to compute the MLE on this part of the parameters is infeasible.

*Multidimensionality and nuisance parameters.* In some models with high dimensional parameters or in the presence of nuisance parameters, either Bayes or frequentist's estimate may be difficult to compute. Various methods ([7, 15, 28], etc.) have been studied for this problem. A proper hybrid formulation may be among the options.

## APPENDIX

**Regularity conditions.** Throughout this paper we assume the densities are with respect to the Lebesque measure. In the following, conditions (A1)–(A3) are A2.1, A2.6 and A2.7 in [5]:

(A1) $\boldsymbol{\theta}$ belongs to an open subset of $R^d$.
(A2) Let $l(\mathbf{x}|\boldsymbol{\theta})$ be the log-likelihood. Assume $\partial l(\mathbf{x}|\boldsymbol{\theta})/\partial \boldsymbol{\theta}$ and $\partial^2 l(\mathbf{x}|\boldsymbol{\theta})/(\partial \boldsymbol{\theta}\, \partial \boldsymbol{\theta}')$ exist and are continuous in $\boldsymbol{\theta}$ for almost all $\mathbf{x}$.
(A3) $E_{\boldsymbol{\theta}}(\sup_{\boldsymbol{\eta}\in\boldsymbol{\Theta}}\|\partial^2 l(\mathbf{x}|\boldsymbol{\eta})/(\partial \boldsymbol{\theta}\, \partial \boldsymbol{\theta}')\|: \|\boldsymbol{\eta}-\boldsymbol{\theta}\|<\varepsilon(\boldsymbol{\theta}))<\infty$ for some $\varepsilon(\boldsymbol{\theta})>0$ and all $\boldsymbol{\theta}\in\boldsymbol{\Theta}$.

Let $P_{\boldsymbol{\theta}}$ be the data distribution given $\boldsymbol{\theta}\in\boldsymbol{\Theta}$, and $l_n(\mathbf{x}^n|\boldsymbol{\theta})=\frac{1}{n}\sum_{i=1}^n l(\mathbf{x}_i|\boldsymbol{\theta})$. Conditions (A4)–(A9) below are those of (1)–(6) in [34].

(A4) The metric space $(\boldsymbol{\Theta},d)$ is separable, where $d(\boldsymbol{\theta},\boldsymbol{\eta})=\|P_{\boldsymbol{\theta}}-P_{\boldsymbol{\eta}}\|$.
(A5) The functions $(l_n(\cdot|\boldsymbol{\theta}))_{\boldsymbol{\theta}\in\boldsymbol{\Theta}}, n\in N$, are separable and measurable.
(A6) $f(\cdot|\boldsymbol{\theta})$, $\boldsymbol{\theta}\in\boldsymbol{\Theta}$, are lower semicontinuous, that is, $\limsup_{n\to\infty} f(\cdot|\boldsymbol{\theta}_n)\leq f(\cdot|\boldsymbol{\theta})$ (a.e.) if $d(\boldsymbol{\theta}_n,\boldsymbol{\theta})\to 0$.
(A7) For every $\boldsymbol{\theta},\boldsymbol{\eta}\in\boldsymbol{\Theta}$, there is an open neighborhood $U_{\boldsymbol{\theta},\boldsymbol{\eta}}$ of $\boldsymbol{\eta}$ such that $E_{\boldsymbol{\theta}}(\inf_{\boldsymbol{\theta}'\in U_{\boldsymbol{\theta},\boldsymbol{\eta}}} l_n(\mathbf{x}^n|\boldsymbol{\theta}'))>-\infty$ for at least one $n$.
(A8) For every $\boldsymbol{\theta}\in\boldsymbol{\Theta}$ and $\varepsilon>0$, $\Pi(\boldsymbol{\eta}\in\boldsymbol{\Theta}: E_{\boldsymbol{\theta}}l(\mathbf{x}|\boldsymbol{\eta})<E_{\boldsymbol{\theta}}l(\mathbf{x}|\boldsymbol{\theta})+\varepsilon)>0$, where $\Pi(\cdot)$ is the distribution for $\pi(\cdot)$.

(A9) For every $\boldsymbol{\theta} \in \boldsymbol{\Theta}$ there is some $n_{\boldsymbol{\theta}}$ such that $P_{\boldsymbol{\theta}}^n(\mathbf{x}^n : \int \prod_{i=1}^n f(\mathbf{x}_i|\boldsymbol{\eta})\Pi(d\boldsymbol{\eta}) < \infty) = 1$ if $n \geq n_{\boldsymbol{\theta}}$.

Conditions (B1)–(B10) are multivariate versions of those of 1–10 in [17].

(B1) For $\boldsymbol{\theta} \neq \boldsymbol{\eta}$, $\int |f(\mathbf{x}|\boldsymbol{\theta}) - f(\mathbf{x}|\boldsymbol{\eta})|\, d\mathbf{x} > 0$.
(B2) For some $p_1 > 0$ and some compact set $\mathbf{K} \in \boldsymbol{\Theta}$, $\sup_{\boldsymbol{\theta} \in \mathbf{K}, \boldsymbol{\eta} \in \boldsymbol{\Theta}} \|\boldsymbol{\theta} - \boldsymbol{\eta}\|^{p_1} \times \int \sqrt{f(\mathbf{x}|\boldsymbol{\theta})f(\mathbf{x}|\boldsymbol{\eta})}\, d\mathbf{x} < \infty$.
(B3) $f(\mathbf{x}|\cdot)$ is continuous on $\boldsymbol{\Theta}^c$, the closure of $\boldsymbol{\Theta}$ on $R^d$ and has $k+2$ ($k \geq 1$) continuous derivatives on $\boldsymbol{\Theta}$.
(B4) (a) For some $b > 0$, and every compact $\mathbf{K} \in \boldsymbol{\Theta}$, $\sup_{\boldsymbol{\theta} \in \mathbf{K}} E_{\boldsymbol{\theta}}\|\mathbf{L}(\mathbf{x}|\boldsymbol{\theta})\|^{3 \vee (k+1+b)} < \infty$.
(b) For every compact $\mathbf{K} \in \boldsymbol{\Theta}$, $\max_{1 \leq |\mathbf{i}| \leq k} \sup_{\boldsymbol{\theta} \in \mathbf{K}} E_{\boldsymbol{\theta}}\|\mathbf{L_i}(\mathbf{x}|\boldsymbol{\theta})\|^{k+1} < \infty$.
(c) For every compact $\mathbf{K} \in \boldsymbol{\Theta}$, and for some $\varepsilon_1(\mathbf{K}) > 0$, $\max_{|\mathbf{i}|=k+1} \sup_{\boldsymbol{\theta} \in \mathbf{K}, \|\boldsymbol{\theta}-\boldsymbol{\eta}\| \leq \varepsilon_1(\mathbf{K})} E_{\boldsymbol{\theta}}\|\mathbf{L_i}(\mathbf{x}|\boldsymbol{\eta})\|^{(k+1)/2} < \infty$.
(B5) (a) For some $p_2 \geq 0$, $\sup_{\boldsymbol{\theta} \in \boldsymbol{\Theta}}(1 + \|\boldsymbol{\theta}\|^{p_2})^{-1}\|\mathbf{I}(\boldsymbol{\theta})\| < \infty$.
(b) $\mathbf{I}(\boldsymbol{\theta})$ is positive definite for $\boldsymbol{\theta} \in \boldsymbol{\Theta}$.
(B6) $\pi(\cdot)$ has $k$ continuous derivatives on $\boldsymbol{\Theta}$.
(B7) For some $p_3 > 0$, $\sup_{\boldsymbol{\theta} \in \boldsymbol{\Theta}}(1 + \|\boldsymbol{\theta}\|^{p_3})^{-1}\pi(\boldsymbol{\theta}) < \infty$.
(B8) $W(\cdot) \geq 0$, is convex, that is, for any $t \in [0,1]$ and $\mathbf{u}_1$ and $\mathbf{u}_2$, $W(t\mathbf{u}_1 + (1-t)\mathbf{u}_2) \leq tW(\mathbf{u}_1) + (1-t)W(\mathbf{u}_2)$.
(B9) For some $\mathbf{a} = (a_1, \ldots, a_d)'$ with $a_j > 1$ to be an even integer ($j = 1, \ldots, d$) and $\varepsilon_2 > 0$, $W(\boldsymbol{\theta}) = \sum_{j=1}^d \theta_j^{a_j}$ for $\|\boldsymbol{\theta}\| \leq \varepsilon_2$.
(B10) For some $p_4 > 0$, $\sup_{\boldsymbol{\theta} \in R^d}(1 + \|\boldsymbol{\theta}\|^{p_4})^{-1}W(\boldsymbol{\theta}) < \infty$.

Note for the 1-dimensional case in [17], (B9) is for some real $a > 1$, $\varepsilon_2 > 0$, $W(\theta) = |\theta|^a$ for $|\theta| \leq \varepsilon_2$. Here we require $\mathbf{a}$ to be a componentwise even integer, otherwise the computation will be unnecessarily involved.

PROOF OF THEOREM 2.1. There is a compact $\boldsymbol{\Omega}' \subset \boldsymbol{\Omega}$ such that $\boldsymbol{\beta}_0 \in \boldsymbol{\Omega}'$. Define a prior on $\boldsymbol{\beta}$ as $\pi(\boldsymbol{\beta}) = 1/L(\boldsymbol{\Omega}')$ for $\boldsymbol{\beta} \in \boldsymbol{\Omega}'$ and $= 0$ for $\boldsymbol{\beta} \in \boldsymbol{\Omega} \setminus \boldsymbol{\Omega}'$, where $L(\cdot)$ is the Lebesque measure on $R^{d_2}$. Then, define the prior $\pi(\boldsymbol{\theta}) = \pi(\boldsymbol{\alpha})\pi(\boldsymbol{\beta})$, the decision $\mathbf{d} = (\mathbf{d}_1, \mathbf{d}_2)$ for $(\boldsymbol{\alpha}, \boldsymbol{\beta})$. Assign the loss $W(\mathbf{d}, \boldsymbol{\theta}) = W(\mathbf{d}_1, \boldsymbol{\alpha}) \times V(\mathbf{d}_2, \boldsymbol{\beta})$, where, for small enough $\delta > 0$, $V(\mathbf{d}_2, \boldsymbol{\beta}) = V(\|\mathbf{d}_2 - \boldsymbol{\beta}\|) = 0$ if $\|\mathbf{d}_2 - \boldsymbol{\beta}\| < \delta$ and 1 otherwise. Here, without confusion, we used $W$ to denote both the loss on $\boldsymbol{\alpha}$ and that on $\boldsymbol{\theta}$. Denote $\pi(\boldsymbol{\theta}|\mathbf{x}^n)$ the posterior density for $\boldsymbol{\theta}$ under the above prior and $R_n = \int W(\mathbf{d}, \boldsymbol{\theta})\pi(\boldsymbol{\theta}|\mathbf{x}^n)\, d\boldsymbol{\theta}$ the posterior risk under the new prior and loss for $\boldsymbol{\theta}$. Then from (2.1), the hybrid estimate $\boldsymbol{\theta}_n$ is the Bayes estimate of $\boldsymbol{\theta}$ under the above new prior and loss. To see this, the Bayes estimate under the new setting is

$$\arg\min_{\mathbf{d}} \int W(\mathbf{d}, \boldsymbol{\theta})\pi(\boldsymbol{\theta}|\mathbf{x}^n)\, d\boldsymbol{\theta}$$

$$= \arg\min_{(\mathbf{d}_1,\mathbf{d}_2)} \int\int W(\mathbf{d}_1,\boldsymbol{\alpha})V(\mathbf{d}_2,\boldsymbol{\beta})\pi(\boldsymbol{\alpha},\boldsymbol{\beta}|\mathbf{x}^n)\,d\boldsymbol{\alpha}\,d\boldsymbol{\beta}$$

$$= \arg\min_{(\mathbf{d}_1,\mathbf{d}_2)} \int\int_{\|\boldsymbol{\beta}-\mathbf{d}_2\|\geq\delta} W(\mathbf{d}_1,\boldsymbol{\alpha})\pi(\boldsymbol{\alpha},\boldsymbol{\beta}|\mathbf{x}^n)\,d\boldsymbol{\alpha}\,d\boldsymbol{\beta}$$

$$= \arg\min,\max_{(\mathbf{d}_1,\boldsymbol{\beta})} \int W(\mathbf{d}_1,\boldsymbol{\alpha})f(\mathbf{x}^n|\boldsymbol{\alpha},\boldsymbol{\beta})\pi(\boldsymbol{\alpha})\,d\boldsymbol{\alpha} = (\check{\boldsymbol{\alpha}}_n,\hat{\boldsymbol{\beta}}_n),$$

as the integration over $\boldsymbol{\beta}$ is minimized when $\mathbf{d}_2$ is the $\boldsymbol{\beta}$-marginal posterior mode, in this case the $\boldsymbol{\beta}$-marginal MLE, while the first integration is minimized by the corresponding marginal Bayes estimator. The $\delta > 0$ above can be arbitrary and the result does not depend on it.

Under the given conditions, by Lemma 2.1 in [5], the MLE of $\boldsymbol{\theta}$ is consistent (a.s.); thus by Theorem 2.5 in [34], for any compact $\mathbf{M} \ni \boldsymbol{\theta}_0$,

$$P\left(\liminf_n \Pi(\mathbf{M}|\mathbf{x}^n) = 1\right) = 1,$$

where $\Pi(\cdot|\mathbf{x}^n)$ is the posterior distribution under the new prior $\pi(\boldsymbol{\theta})$.

Let $\boldsymbol{\theta}_0 = (\boldsymbol{\alpha}_0,\boldsymbol{\beta}_0)$ be the true parametric value. By the given conditions on $W(\cdot)$, for $\varepsilon > 0$, there is a $\delta > 0$ such that $W(\|\boldsymbol{\alpha}_0 - \boldsymbol{\alpha}\|) \leq \varepsilon$ as long as $\|\boldsymbol{\alpha}_0 - \boldsymbol{\alpha}\| \leq \delta$. Let $\mathbf{M} = (\boldsymbol{\theta}_0 \pm \delta\mathbf{1})$ be the $\delta$ neighborhood of $\boldsymbol{\theta}_0$, then $\sup_{\boldsymbol{\theta}\in\mathbf{M}} W(\boldsymbol{\theta}_0,\boldsymbol{\theta}) \leq \varepsilon$. Since $\pi(\cdot|\mathbf{x}^n) \to 0$ (a.s.) on $\mathbf{M}^c$, and for large $n$ it can be dominated by $\pi(\cdot|\mathbf{x}^k)$ on $\mathbf{M}^c$, so for large $n$ we have $\int_{\mathbf{M}^c} W(\boldsymbol{\theta}_0,\boldsymbol{\theta})\pi(\boldsymbol{\theta}|\mathbf{x}^n)\,d\boldsymbol{\theta} < \varepsilon$. Since $\boldsymbol{\theta}_n$ is Bayesian under the new prior and loss, (2.1) is rewritten as

$$\boldsymbol{\theta}_n = \arg\min_{\boldsymbol{\delta}} \int W(\boldsymbol{\delta},\boldsymbol{\theta})\pi(\boldsymbol{\theta}|\mathbf{x}^n)\,d\boldsymbol{\theta}.$$

So, for large $n$, we have the posterior risk

$$R_n = \int W(\boldsymbol{\theta}_n,\boldsymbol{\theta})\pi(\boldsymbol{\theta}|\mathbf{x}^n)\,d\boldsymbol{\theta} \leq \int W(\boldsymbol{\theta}_0,\boldsymbol{\theta})\pi(\boldsymbol{\theta}|\mathbf{x}^n)\,d\boldsymbol{\theta}$$

$$= \int_{\mathbf{M}} W(\boldsymbol{\theta}_0,\boldsymbol{\theta})\pi(\boldsymbol{\theta}|\mathbf{x}^n)\,d\boldsymbol{\theta} + \int_{\mathbf{M}^c} W(\boldsymbol{\theta}_0,\boldsymbol{\theta})\pi(\boldsymbol{\theta}|\mathbf{x}^n)\,d\boldsymbol{\theta}$$

$$\leq \int_{\mathbf{M}} W(\boldsymbol{\theta}_0,\boldsymbol{\theta})\pi(\boldsymbol{\theta}|\mathbf{x}^n)\,d\boldsymbol{\theta} + \varepsilon \leq \varepsilon\int \pi(\boldsymbol{\theta}|\mathbf{x}^n)\,d\boldsymbol{\theta} + \varepsilon = 2\varepsilon.$$

Since $\varepsilon > 0$ is arbitrary, we have $R_n \to 0$ (a.s.).

Suppose $\boldsymbol{\theta}_n$ is not consistent (a.s.) to $\boldsymbol{\theta}_0$, or $\limsup \|\boldsymbol{\theta}_n - \boldsymbol{\theta}_0\| \geq \varepsilon$ (a.s.) for some $\varepsilon > 0$. Then there is a sub-sequence $\{n_k\}$ such that either (a) $\lim_k \|\boldsymbol{\alpha}_{n_k} - \boldsymbol{\alpha}_0\| \geq \varepsilon/2$ and $\lim_k \|\boldsymbol{\beta}_{n_k} - \boldsymbol{\beta}_0\| \geq \varepsilon/2$; or (b) $\lim_k \|\boldsymbol{\alpha}_{n_k} - \boldsymbol{\alpha}_0\| \geq \varepsilon$ but $\lim_k \|\boldsymbol{\beta}_{n_k} - \boldsymbol{\beta}_0\| \to 0$; or (c) $\lim_k \|\boldsymbol{\alpha}_{n_k} - \boldsymbol{\alpha}_0\| \to 0$ but $\lim_k \|\boldsymbol{\beta}_{n_k} - \boldsymbol{\beta}_0\| \geq \varepsilon$. For case (a), let $\mathbf{M} = \{\boldsymbol{\theta}: \|\boldsymbol{\theta} - \boldsymbol{\theta}_0\| \leq \varepsilon/2\}$, then $\mathbf{M} \subset \mathbf{M}_1 := \{\boldsymbol{\theta}: \|\boldsymbol{\alpha} - \boldsymbol{\alpha}_{n_k}\| \geq \varepsilon/2; \|\boldsymbol{\beta} - \boldsymbol{\beta}_{n_k}\| \geq \varepsilon/2\}$. Note for $\boldsymbol{\theta} = (\boldsymbol{\alpha},\boldsymbol{\beta}) \in \mathbf{M}_1$, $W(\|\boldsymbol{\alpha} - \boldsymbol{\alpha}_{n_k}\|) \geq W(\varepsilon/2)$,

and $\Pi(\mathbf{M}|\mathbf{x}^{n_k}) \to 1$ (a.s.) by the previous result. Also, by our choice of the prior on $\boldsymbol{\beta}$, we have for some $0 < c < 1$,

$$\int_{\mathbf{M}_1} V(\boldsymbol{\beta}_{n_k}, \boldsymbol{\beta}) \pi(\boldsymbol{\theta}|\mathbf{x}^{n_k}) d\boldsymbol{\theta} \geq \int_{\mathbf{M}} V(\boldsymbol{\beta}_{n_k}, \boldsymbol{\beta}) \pi(\boldsymbol{\theta}|\mathbf{x}^{n_k}) d\boldsymbol{\theta} \geq c\Pi(\mathbf{M}|\mathbf{x}^{n_k}).$$

So we have, for all $n_k$,

$$R_{n_k} \geq \int_{\mathbf{M}_1} W(\boldsymbol{\theta}_{n_k}, \boldsymbol{\theta}) \pi(\boldsymbol{\theta}|\mathbf{x}^{n_k}) d\boldsymbol{\theta}$$

$$\geq W(\varepsilon/2) \int_{\mathbf{M}_1} V(\boldsymbol{\beta}_{n_k}, \boldsymbol{\beta}) \pi(\boldsymbol{\theta}|\mathbf{x}^{n_k}) d\boldsymbol{\theta}$$

$$\geq W(\varepsilon/2) \int_{\mathbf{M}} V(\boldsymbol{\beta}_{n_k}, \boldsymbol{\beta}) \pi(\boldsymbol{\theta}|\mathbf{x}^{n_k}) d\boldsymbol{\theta}$$

$$\geq cW(\varepsilon/2) \Pi(\mathbf{M}|\mathbf{x}^{n_k}) \to cW(\varepsilon/2) > 0,$$

which is a contradiction to the fact that $R_n \to 0$ (a.s.), and so (a) cannot be true.

For case (b), let $\Pi_{1,n_k}(\cdot|\mathbf{x}^{n_k})$ be the $\boldsymbol{\alpha}$-marginal posterior distribution with $\boldsymbol{\beta}$ evaluated at $\boldsymbol{\beta}_{n_k}$. Similarly as before, let $\mathbf{M} = (\boldsymbol{\alpha}_0 \pm \delta \mathbf{1})$. Since $\boldsymbol{\beta}_{n_k} \to \boldsymbol{\beta}_0$, we have $P(\lim_k \inf \Pi_{1,n_k}(\mathbf{M}|\mathbf{x}^{n_k}) = 1) = 1$, and as before, the posterior risk

$$R_{n_k} := \int W(\boldsymbol{\alpha}_{n_k}, \boldsymbol{\alpha}) \pi(\boldsymbol{\alpha}, \boldsymbol{\beta}_{n_k}|\mathbf{x}^{n_k}) d\boldsymbol{\alpha}$$

$$\leq \int_{\mathbf{M}} W(\boldsymbol{\alpha}_0, \boldsymbol{\alpha}) \pi(\boldsymbol{\alpha}, \boldsymbol{\beta}_{n_k}|\mathbf{x}^{n_k}) d\boldsymbol{\alpha}$$

$$+ \int_{\mathbf{M}^c} W(\boldsymbol{\alpha}_0, \boldsymbol{\alpha}) \pi(\boldsymbol{\alpha}, \boldsymbol{\beta}_{n_k}|\mathbf{x}^{n_k}) d\boldsymbol{\alpha}$$

$$\leq \varepsilon \int \pi(\boldsymbol{\alpha}, \boldsymbol{\beta}_{n_k}|\mathbf{x}^{n_k}) d\boldsymbol{\alpha} + \varepsilon \leq 2\varepsilon,$$

i.e. $R_{n_k} \to 0$. On the other hand, for the fixed $\varepsilon > 0$, let $\mathbf{M} = \{\boldsymbol{\alpha} : \|\boldsymbol{\alpha} - \boldsymbol{\alpha}_0\| \leq \varepsilon/2\}$, then $\mathbf{M} \subset \mathbf{M}_1 := \{\boldsymbol{\alpha} : \|\boldsymbol{\alpha} - \boldsymbol{\alpha}_{n_k}\| \geq \varepsilon/2\}$. We have

$$R_{n_k} \geq \int_{\mathbf{M}_1} W(\boldsymbol{\alpha}_{n_k}, \boldsymbol{\alpha}) \pi(\boldsymbol{\alpha}, \boldsymbol{\beta}_{n_k}|\mathbf{x}^{n_k}) d\boldsymbol{\alpha} \geq W(\varepsilon/2) \Pi_{1,n_k}(\mathbf{M}|\mathbf{x}^{n_k}) \to W(\varepsilon/2) > 0,$$

a contradiction.

By similar argument on the $\boldsymbol{\beta}$-margin, it is seen that case (c) cannot be true. $\square$

PROOF OF THEOREM 2.2. We check the main steps, for the vector version, of the proof of Theorem 1 in [17]. Lemmas 1–6 of Gusev [17] there can be checked similarly. Let $\hat{\boldsymbol{\theta}}'_n = \sqrt{n}(\hat{\boldsymbol{\theta}}_n - \boldsymbol{\theta}_0)$. By definition of $\hat{\boldsymbol{\theta}}_n$, we have

(A.1) $$\mathbf{S_0}(\hat{\boldsymbol{\theta}}_n) + n^{-1/2} \boldsymbol{\rho_0}(\hat{\boldsymbol{\theta}}_n) = \mathbf{0}.$$

Also,

$$\mathbf{S_0}(\hat{\boldsymbol{\theta}}_n) = \mathbf{S_0}(\boldsymbol{\theta}_0 + n^{-1/2}\hat{\boldsymbol{\theta}}'_n) \overset{k}{\sim} \sum_{r=0}^{k-1} n^{-r/2} \sum_{|\mathbf{i}|=r} \frac{\langle (\hat{\boldsymbol{\theta}}'_n)^{\mathbf{i}} \rangle}{\mathbf{i}!} \mathbf{S_i}$$

$$\overset{k}{\sim} \sum_{r=0}^{k-1} n^{-r/2} \left( \sum_{|\mathbf{i}|=r} \frac{\langle (\hat{\boldsymbol{\theta}}'_n)^{\mathbf{i}} \rangle}{\mathbf{i}!} \boldsymbol{\Delta_i} + \sum_{|\mathbf{i}|=r+1} \frac{\langle (\hat{\boldsymbol{\theta}}'_n)^{\mathbf{i}} \rangle}{\mathbf{i}!} \mathbf{E_i} \right)$$

and

$$n^{-1/2} \boldsymbol{\rho_0}(\hat{\boldsymbol{\theta}}_n) \overset{k}{\sim} \sum_{r=1}^{k-1} n^{-r/2} \sum_{|\mathbf{i}|=r-1} \frac{\langle (\hat{\boldsymbol{\theta}}'_n)^{\mathbf{i}} \rangle}{\mathbf{i}!} \boldsymbol{\rho_i}.$$

In the above, we used the relationship $\mathbf{S_i} = \boldsymbol{\Delta_i} + n^{1/2} \mathbf{E_i}$, and $\mathbf{E_0} = \mathbf{0}$.

Note $\mathbf{I} = -(\mathbf{E_i} : |\mathbf{i}| = 1)'$, the Fisher information matrix, so $\sum_{|\mathbf{i}|=1} \langle (\hat{\boldsymbol{\theta}}'_n)^{\mathbf{i}} \rangle \mathbf{E_i} = -\mathbf{I}\hat{\boldsymbol{\theta}}'_n$. Now, we rewrite (A.1) as

$$\boldsymbol{\Delta_0} - \mathbf{I}\hat{\boldsymbol{\theta}}'_n + \sum_{r=1}^{k-1} n^{-r/2} \left( \sum_{|\mathbf{i}|=r-1} \frac{\langle (\hat{\boldsymbol{\theta}}'_n)^{\mathbf{i}} \rangle}{\mathbf{i}!} \boldsymbol{\rho_i} + \sum_{|\mathbf{i}|=r} \frac{\langle (\hat{\boldsymbol{\theta}}'_n)^{\mathbf{i}} \rangle}{\mathbf{i}!} \boldsymbol{\Delta_i} + \sum_{|\mathbf{i}|=r+1} \frac{\langle (\hat{\boldsymbol{\theta}}'_n)^{\mathbf{i}} \rangle}{\mathbf{i}!} \mathbf{E_i} \right)$$
(A.2)
$$\overset{k}{\sim} \mathbf{0}.$$

Consider, with $t = 0$, the term $\sum_{s+t=r} \sum_{|\mathbf{i}|=t+1} \mathbf{E_i} \sum_{|\mathbf{l}|=s} \sum_{(0,\mathbf{l},\mathbf{i})} \prod_{v}^{s} \frac{\langle \mathbf{H}_v^{\mathbf{i}_v} \rangle}{\mathbf{i}_v!}$. Note $\mathbf{i}$ can only take one of the vectors $\mathbf{e}_j$'s. First we take $\mathbf{i} = \mathbf{e}_1$, since $|\mathbf{l}| = r$ it is easy to check that the only nonempty integer vector sets satisfy the definition of $\sum_{(0,\mathbf{l},\mathbf{i})}$ is $\{\mathbf{l} = r\mathbf{e}_1\}$, and the only $\mathbf{i}_v$ in this set is $\{\mathbf{i}_r = \mathbf{e}_1\}$. Similarly, for $\mathbf{i} = \mathbf{e}_2$, the only nonempty integer vector sets satisfy the definition of $\sum_{(0,\mathbf{l},\mathbf{i})}$ is $\{\mathbf{l} = r\mathbf{e}_2\}$, and the only $\mathbf{i}_v$ in this set is $\{\mathbf{i}_r = \mathbf{e}_2\}, \ldots$, so we have, with $t = 0$,

$$\sum_{s+t=r} \sum_{|\mathbf{i}|=t+1} \mathbf{E_i} \sum_{|\mathbf{l}|=s} \sum_{(0,s,\mathbf{l},\mathbf{i})} \prod_{v=0}^{s} \frac{\langle \mathbf{H}_v^{\mathbf{i}_v} \rangle}{\mathbf{i}_v!} = \sum_{|\mathbf{i}|=1} \mathbf{E_i} \sum_{|\mathbf{l}|=r} \sum_{(0,s,\mathbf{l},\mathbf{i})} \prod_{v=0}^{r} \frac{\langle \mathbf{H}_v^{\mathbf{i}_v} \rangle}{\mathbf{i}_v!}$$

$$= \sum_{|\mathbf{i}|=1} \mathbf{E_i} \frac{\langle \mathbf{H}_r^{\mathbf{i}} \rangle}{1!} = -\mathbf{I}\mathbf{H}_r.$$

Thus, the expression for $\mathbf{H}_r$ in Theorem 2.2 is rewritten as

$$\sum_{s+t=r} \left( \sum_{|\mathbf{i}|=t-1} \boldsymbol{\rho_i} \sum_{|\mathbf{l}|=s} \sum_{(0,s,\mathbf{l},\mathbf{i})} \prod_{v=0}^{s} \frac{\langle \mathbf{H}_v^{\mathbf{i}_v} \rangle}{\mathbf{i}_v!} + \sum_{|\mathbf{i}|=t} \boldsymbol{\Delta_i} \sum_{|\mathbf{l}|=s} \sum_{(0,s,\mathbf{l},\mathbf{i})} \prod_{v=0}^{s} \frac{\langle \mathbf{H}_v^{\mathbf{i}_v} \rangle}{\mathbf{i}_v!} \right.$$
(A.3)
$$\left. + \sum_{|\mathbf{i}|=t+1} \mathbf{E_i} \sum_{|\mathbf{l}|=s} \sum_{(0,s,\mathbf{l},\mathbf{i})} \prod_{v=0}^{s} \frac{\langle \mathbf{H}_v^{\mathbf{i}_v} \rangle}{\mathbf{i}_v!} \right) \overset{k}{\sim} \mathbf{0}.$$

Set $\hat{\boldsymbol{\theta}}_n'' = \sum_{i=0}^{k-1} \mathbf{H}_i n^{-i/2} := (\hat{\theta}_{n1}'', \ldots, \hat{\theta}_{nd}'')'$, so $\langle (\hat{\boldsymbol{\theta}}_n'')^{\mathbf{i}} \rangle = \prod_{j=1}^{d} (\hat{\theta}_{nj}'')^{i_j}$.

For integers $r \geq a$ and integers $m, l \geq 0$, $I_1(a, r, m, l)$ denote the set of nonnegative integers $(i_1, \ldots, i_r)$,

$$I_1(a, r, m, l) = \left\{ (i_1, \ldots, i_r) : \sum_{v=a}^{r} v i_v = m, \ \sum_{v=a}^{r} i_v = l \right\}.$$

Write $\mathbf{H}_i = (H_{i,1}, \ldots, H_{i,d})'$ and $\mathbf{i} = (i_1, \ldots, i_d)'$. Note $\hat{\theta}_{n,j}'' = \sum_{r=0}^{k-1} n^{-r/2} H_{r,j}$, and

$$(\hat{\theta}_{n,j}'')^{i_j} \stackrel{k}{\sim} i_j! \sum_{r=0}^{k-1} n^{-r/2} \sum_{I_1(0,r,r,i_j)} \prod_{v=0}^{r} H_{v,j}^{t_v}/t_v!.$$

It can be checked that

$$\prod_{j=1}^{d} \sum_{I_1(0,r,r,i_j)} \prod_{v=0}^{r} H_{v,j}^{t_v}/t_v! = \sum_{|\mathbf{l}|=r} \sum_{(0,r,\mathbf{l},\mathbf{i})} \prod_{v=0}^{r} \langle \mathbf{H}_v^{\mathbf{i}_v} \rangle / \mathbf{i}_v!,$$

thus

$$\langle (\hat{\boldsymbol{\theta}}_n'')^{\mathbf{i}} \rangle \stackrel{k}{\sim} \mathbf{i}! \sum_{r=0}^{k-1} n^{-r/2} \prod_{j=1}^{d} \sum_{I_1(0,r,r,i_j)} \prod_{v=0}^{r} H_{v,j}^{t_v}/t_v!$$

$$= \mathbf{i}! \sum_{r=0}^{k-1} n^{-r/2} \sum_{|\mathbf{l}|=r} \sum_{(0,r,\mathbf{l},\mathbf{i})} \prod_{v=0}^{r} \langle \mathbf{H}_v^{\mathbf{i}_v} \rangle / \mathbf{i}_v!.$$

Note (A.2) still holds with $\hat{\boldsymbol{\theta}}_n'$ replaced by $\hat{\boldsymbol{\theta}}_n''$, and using (A.3) we get

$$\boldsymbol{\Delta_0} - \mathbf{I}\hat{\boldsymbol{\theta}}_n'' + \sum_{r=1}^{k-1} n^{-r/2} \left( \sum_{|\mathbf{i}|=r-1} \frac{\langle (\hat{\boldsymbol{\theta}}_n'')^{\mathbf{i}} \rangle}{\mathbf{i}!} \boldsymbol{\rho_i} + \sum_{|\mathbf{i}|=r} \frac{\langle (\hat{\boldsymbol{\theta}}_n'')^{\mathbf{i}} \rangle}{\mathbf{i}!} \boldsymbol{\Delta_i} + \sum_{|\mathbf{i}|=r+1} \frac{\langle (\hat{\boldsymbol{\theta}}_n'')^{\mathbf{i}} \rangle}{\mathbf{i}!} \mathbf{E_i} \right)$$

$$= \sum_{r=0}^{k-1} n^{-r/2} \left( \sum_{|\mathbf{i}|=r-1} \frac{\langle (\hat{\boldsymbol{\theta}}_n'')^{\mathbf{i}} \rangle}{\mathbf{i}!} \boldsymbol{\rho_i} + \sum_{|\mathbf{i}|=r} \frac{\langle (\hat{\boldsymbol{\theta}}_n'')^{\mathbf{i}} \rangle}{\mathbf{i}!} \boldsymbol{\Delta_i} + \sum_{|\mathbf{i}|=r+1} \frac{\langle (\hat{\boldsymbol{\theta}}_n'')^{\mathbf{i}} \rangle}{\mathbf{i}!} \mathbf{E_i} \right)$$

$$\stackrel{k}{\sim} \sum_{r=0}^{k-1} \sum_{s=0}^{k-1} n^{-(r+s)/2} \left( \sum_{|\mathbf{i}|=r-1} \boldsymbol{\rho_i} \sum_{|\mathbf{l}|=s} \sum_{(0,s,\mathbf{l},\mathbf{i})} \prod_{v=0}^{s} \frac{\langle \mathbf{H}_v^{\mathbf{i}_v} \rangle}{\mathbf{i}_v!} \right.$$

$$+ \sum_{|\mathbf{i}|=r} \boldsymbol{\Delta_i} \sum_{|\mathbf{l}|=s} \sum_{(0,s,\mathbf{l},\mathbf{i})} \prod_{v=0}^{s} \frac{\langle \mathbf{H}_v^{\mathbf{i}_v} \rangle}{\mathbf{i}_v!}$$

(A.4) $$\left. + \sum_{|\mathbf{i}|=r+1} \mathbf{E_i} \sum_{|\mathbf{l}|=s} \sum_{(0,s,\mathbf{l},\mathbf{i})} \prod_{v=0}^{s} \frac{\langle \mathbf{H}_v^{\mathbf{i}_v} \rangle}{\mathbf{i}_v!} \right)$$

$$\overset{k}{\sim} \sum_{r=0}^{k-1} n^{-r/2} \sum_{s+t=r} \left( \sum_{|\mathbf{i}|=t-1} \boldsymbol{\rho}_{\mathbf{i}} \sum_{|\mathbf{l}|=s} \sum_{(0,s,\mathbf{l},\mathbf{i})} \prod_{v=0}^{s} \frac{\langle \mathbf{H}_v^{\mathbf{i}_v} \rangle}{\mathbf{i}_v!} \right.$$

$$+ \sum_{|\mathbf{i}|=t} \boldsymbol{\Delta}_{\mathbf{i}} \sum_{|\mathbf{l}|=s} \sum_{(0,s,\mathbf{l},\mathbf{i})} \prod_{v=0}^{s} \frac{\langle \mathbf{H}_v^{\mathbf{i}_v} \rangle}{\mathbf{i}_v!}$$

$$\left. + \sum_{|\mathbf{i}|=t+1} \mathbf{E}_{\mathbf{i}} \sum_{|\mathbf{l}|=s} \sum_{(0,s,\mathbf{l},\mathbf{i})} \prod_{v=0}^{s} \frac{\langle \mathbf{H}_v^{\mathbf{i}_v} \rangle}{\mathbf{i}_v!} \right) \overset{k}{\sim} \mathbf{0}.$$

Now (A.2) minus the left-hand side of (A.4), similarly as in [17], we get

$$\mathbf{0} \overset{k}{\sim} \mathbf{I}(\hat{\boldsymbol{\theta}}'_n - \hat{\boldsymbol{\theta}}''_n) + O(n^{-1/2}) \mathbf{I}(\hat{\boldsymbol{\theta}}'_n - \hat{\boldsymbol{\theta}}''_n) = (1 + O(n^{-1/2})) \mathbf{I}(\hat{\boldsymbol{\theta}}'_n - \hat{\boldsymbol{\theta}}''_n).$$

Thus, $\hat{\boldsymbol{\theta}}'_n \overset{k}{\sim} \hat{\boldsymbol{\theta}}''_n$, and Theorem 2.2 is proved. $\square$

Recall $\hat{\boldsymbol{\theta}}'_n = \sqrt{n}(\hat{\boldsymbol{\theta}}_n - \boldsymbol{\theta}_0)$, and define

$$Z_n(\boldsymbol{\theta}) = \left( \prod_{i=1}^{n} \frac{f(x_i | \boldsymbol{\theta}_0 + \boldsymbol{\theta} n^{-1/2})}{f(x_i | \boldsymbol{\theta}_0)} \right) \frac{\pi(\boldsymbol{\theta}_0 + \boldsymbol{\theta} n^{-1/2})}{\pi(\boldsymbol{\theta}_0)}.$$

We first extend a result in [17] to the multivariate case.

LEMMA 1. *Under the conditions of Theorem 2.2, we have*

$$\frac{Z_n(\boldsymbol{\theta} + \hat{\boldsymbol{\theta}}'_n)}{Z_n(\hat{\boldsymbol{\theta}}'_n)} = \exp\left(-\frac{1}{2} \boldsymbol{\theta}' \mathbf{I} \boldsymbol{\theta}\right) \left(1 + \sum_{r=1}^{k-1} n^{-r/2} \sum_{|\mathbf{i}|=2}^{3r} \langle \boldsymbol{\theta}^{\mathbf{i}} \rangle N_{\mathbf{i},r} \right) + O_p(n^{-k/2}),$$

*where*

$$N_{\mathbf{i},r} = \sum_{I_2(r,\mathbf{i})} \prod_{v=1}^{r} \sum_{I_1(2,v+2,\mathbf{k}_v,i_v)} \prod_{|\mathbf{j}|=2}^{v+2} \frac{F_{\mathbf{j},v}^{u_{\mathbf{j}}}}{u_{\mathbf{j}}!(\mathbf{j}!)^{u_{\mathbf{j}}}};$$

*in the above, the summations are for $i_v \in I_0(1,r,r)$, $\mathbf{k}_v \in I_2(r,\mathbf{i})$, and for each fixed $v$, $u_{\mathbf{j}} \in I_1(2,v+2,\mathbf{k}_v,i_v)$, with the notation $I_0(1,r,r)$, $I_2(r,\mathbf{i})$ and $I_1(m,r,\mathbf{k}_v,i_v)$ given at the end of the proof, and*

$$F_{\mathbf{i},r} = \sum_{t+s=r, t \geq 1} \sum_{\mathbf{j} \geq \mathbf{i}, |\mathbf{j}|=t} \varrho_{\mathbf{j}} \sum_{|\mathbf{l}|=s} \sum_{(0,s,\mathbf{l},\mathbf{j}-\mathbf{i})} \prod_{v=0}^{s} \langle \mathbf{H}_v^{\mathbf{i}_v} \rangle / \mathbf{i}_v!$$

$$+ \sum_{t+s=r+1, t \geq 2} \sum_{\mathbf{j} \geq \mathbf{i}, |\mathbf{j}|=t} \delta_{\mathbf{j}} \sum_{|\mathbf{l}|=s} \sum_{(0,s,\mathbf{l},\mathbf{j}-\mathbf{i})} \prod_{v=0}^{s} \langle \mathbf{H}_v^{\mathbf{i}_v} \rangle / \mathbf{i}_v!$$

$$+ \sum_{t+s=r+2, t \geq 3} \sum_{\mathbf{j} \geq \mathbf{i}, |\mathbf{j}|=t} \mathcal{E}_{\mathbf{j}} \sum_{|\mathbf{l}|=s} \sum_{(0,s,\mathbf{l},\mathbf{j}-\mathbf{i})} \prod_{v=0}^{s} \langle \mathbf{H}_v^{\mathbf{i}_v} \rangle / \mathbf{i}_v!.$$

PROOF. As in the proof of Theorem 2 in [17], we have

$$\log Z_n(\boldsymbol{\theta} + \hat{\boldsymbol{\theta}}'_n) \stackrel{k}{\sim} \sum_{r=0}^{k} n^{-r/2} \sum_{|\mathbf{j}|=r+1} \frac{\langle(\boldsymbol{\theta} + \hat{\boldsymbol{\theta}}'_n)^{\mathbf{j}}\rangle}{\mathbf{j}!} \mathcal{S}_{\mathbf{j}}$$

$$+ \sum_{r=1}^{k-1} n^{-r/2} \sum_{|\mathbf{j}|=r} \frac{\langle(\boldsymbol{\theta} + \hat{\boldsymbol{\theta}}'_n)^{\mathbf{j}}\rangle}{\mathbf{j}!} \varrho_{\mathbf{j}}$$

$$= \sum_{r=0}^{k} n^{-r/2} \sum_{|\mathbf{j}|=r+1} \mathcal{S}_{\mathbf{j}} \sum_{\mathbf{i} \leq \mathbf{j}} \frac{\langle\boldsymbol{\theta}^{\mathbf{i}}\rangle\langle(\hat{\boldsymbol{\theta}}'_n)^{\mathbf{j}-\mathbf{i}}\rangle}{\mathbf{i}!(\mathbf{j}-\mathbf{i})!}$$

$$+ \sum_{r=1}^{k-1} n^{-r/2} \sum_{|\mathbf{j}|=r} \varrho_{\mathbf{j}} \sum_{\mathbf{i} \leq \mathbf{j}} \frac{\langle\boldsymbol{\theta}^{\mathbf{i}}\rangle\langle(\hat{\boldsymbol{\theta}}'_n)^{\mathbf{j}-\mathbf{i}}\rangle}{\mathbf{i}!(\mathbf{j}-\mathbf{i})!}$$

and

$$\log Z_n(\hat{\boldsymbol{\theta}}'_n) \stackrel{k}{\sim} \sum_{r=0}^{k} n^{-r/2} \sum_{|\mathbf{j}|=r+1} \frac{\langle(\hat{\boldsymbol{\theta}}'_n)^{\mathbf{j}}\rangle}{\mathbf{j}!} \mathcal{S}_{\mathbf{j}} + \sum_{r=1}^{k-1} n^{-r/2} \sum_{|\mathbf{j}|=r} \frac{\langle(\hat{\boldsymbol{\theta}}'_n)^{\mathbf{j}}\rangle}{\mathbf{j}!} \varrho_{\mathbf{j}}.$$

These give

(A.5)
$$\log \frac{Z_n(\boldsymbol{\theta} + \hat{\boldsymbol{\theta}}'_n)}{Z_n(\hat{\boldsymbol{\theta}}'_n)} \stackrel{k}{\sim} \sum_{r=0}^{k} n^{-r/2} \sum_{|\mathbf{j}|=r+1} \mathcal{S}_{\mathbf{j}} \sum_{\mathbf{i} \leq \mathbf{j}, |\mathbf{i}|=1} \frac{\langle\boldsymbol{\theta}^{\mathbf{i}}\rangle\langle(\hat{\boldsymbol{\theta}}'_n)^{\mathbf{j}-\mathbf{i}}\rangle}{(\mathbf{j}-\mathbf{i})!}$$

$$+ \sum_{r=1}^{k-1} n^{-r/2} \sum_{|\mathbf{j}|=r} \varrho_{\mathbf{j}} \sum_{\mathbf{i} \leq \mathbf{j}, |\mathbf{i}|=1} \frac{\langle\boldsymbol{\theta}^{\mathbf{i}}\rangle\langle(\hat{\boldsymbol{\theta}}'_n)^{\mathbf{j}-\mathbf{i}}\rangle}{(\mathbf{j}-\mathbf{i})!}$$

$$+ \sum_{r=1}^{k} n^{-r/2} \sum_{|\mathbf{j}|=r+1} \mathcal{S}_{\mathbf{j}} \sum_{\mathbf{i} \leq \mathbf{j}, |\mathbf{i}| \geq 2} \frac{\langle\boldsymbol{\theta}^{\mathbf{i}}\rangle\langle(\hat{\boldsymbol{\theta}}'_n)^{\mathbf{j}-\mathbf{i}}\rangle}{\mathbf{i}!(\mathbf{j}-\mathbf{i})!}$$

$$+ \sum_{r=2}^{k-1} n^{-r/2} \sum_{|\mathbf{j}|=r} \varrho_{\mathbf{j}} \sum_{\mathbf{i} \leq \mathbf{j}, |\mathbf{i}| \geq 2} \frac{\langle\boldsymbol{\theta}^{\mathbf{i}}\rangle\langle(\hat{\boldsymbol{\theta}}'_n)^{\mathbf{j}-\mathbf{i}}\rangle}{\mathbf{i}!(\mathbf{j}-\mathbf{i})!}.$$

Recall for $|\mathbf{i}| = 1$ we have

$$0 = \mathcal{S}_{\mathbf{i}}(\hat{\boldsymbol{\theta}}_n) + n^{-1/2}\varrho_{\mathbf{i}}(\hat{\boldsymbol{\theta}}_n)$$
$$= \mathcal{S}_{\mathbf{i}}(\boldsymbol{\theta}_0 + n^{-1/2}\hat{\boldsymbol{\theta}}'_n) + n^{-1/2}\varrho_{\mathbf{i}}(\boldsymbol{\theta}_0 + n^{-1/2}\hat{\boldsymbol{\theta}}'_n),$$

$$\mathcal{S}_{\mathbf{i}}(\boldsymbol{\theta}_0 + n^{-1/2}\hat{\boldsymbol{\theta}}'_n) \stackrel{k}{\sim} \sum_{r=0}^{k} n^{-r/2} \sum_{|\mathbf{j}|=r} \mathcal{S}_{\mathbf{i}+\mathbf{j}} \frac{\langle(\hat{\boldsymbol{\theta}}'_n)^{\mathbf{j}}\rangle}{\mathbf{j}!}$$

and

$$\varrho_{\mathbf{i}}(\boldsymbol{\theta}_0 + n^{-1/2}\hat{\boldsymbol{\theta}}'_n) \stackrel{k}{\sim} \sum_{r=1}^{k-1} n^{-r/2} \sum_{|\mathbf{j}|=r-1} \varrho_{\mathbf{i}+\mathbf{j}} \frac{\langle (\hat{\boldsymbol{\theta}}'_n)^{\mathbf{j}} \rangle}{\mathbf{j}!},$$

so the first two summations in the right-hand side of (A.5) together is

$$\sum_{|\mathbf{i}|=1} \langle \boldsymbol{\theta}^{\mathbf{i}} \rangle \left( \sum_{r=0}^{k} n^{-r/2} \sum_{|\mathbf{j}|=r} \mathcal{S}_{\mathbf{i}+\mathbf{j}} \frac{\langle (\hat{\boldsymbol{\theta}}'_n)^{\mathbf{j}} \rangle}{\mathbf{j}!} + \sum_{r=1}^{k-1} n^{-r/2} \sum_{|\mathbf{j}|=r-1} \varrho_{\mathbf{i}+\mathbf{j}} \frac{\langle (\hat{\boldsymbol{\theta}}'_n)^{\mathbf{j}} \rangle}{\mathbf{j}!} \right) = 0.$$

Using $\mathcal{S}_{\mathbf{j}} = \delta_{\mathbf{j}} + \sqrt{n} \mathcal{E}_{\mathbf{j}}$, (A.5) is now

(A.6)
$$\begin{aligned}
\log \frac{Z_n(\boldsymbol{\theta}+\hat{\boldsymbol{\theta}}'_n)}{Z_n(\hat{\boldsymbol{\theta}}'_n)} &\stackrel{k}{\sim} \sum_{r=2}^{k-1} n^{-r/2} \sum_{|\mathbf{j}|=r} \varrho_{\mathbf{j}} \sum_{\mathbf{i} \leq \mathbf{j}, |\mathbf{i}| \geq 2} \frac{\langle \boldsymbol{\theta}^{\mathbf{i}} \rangle \langle (\hat{\boldsymbol{\theta}}'_n)^{\mathbf{j}-\mathbf{i}} \rangle}{\mathbf{i}!(\mathbf{j}-\mathbf{i})!} \\
&\quad + \sum_{r=1}^{k-1} n^{-r/2} \sum_{|\mathbf{j}|=r+1} \delta_{\mathbf{j}} \sum_{\mathbf{i} \leq \mathbf{j}, |\mathbf{i}| \geq 2} \frac{\langle \boldsymbol{\theta}^{\mathbf{i}} \rangle \langle (\hat{\boldsymbol{\theta}}'_n)^{\mathbf{j}-\mathbf{i}} \rangle}{\mathbf{i}!(\mathbf{j}-\mathbf{i})!} \\
&\quad + \sum_{r=0}^{k-1} n^{-r/2} \sum_{|\mathbf{j}|=r+2} \mathcal{E}_{\mathbf{j}} \sum_{\mathbf{i} \leq \mathbf{j}, |\mathbf{i}| \geq 2} \frac{\langle \boldsymbol{\theta}^{\mathbf{i}} \rangle \langle (\hat{\boldsymbol{\theta}}'_n)^{\mathbf{j}-\mathbf{i}} \rangle}{\mathbf{i}!(\mathbf{j}-\mathbf{i})!} \\
&= \sum_{r=1}^{k-1} n^{-r/2} \sum_{|\mathbf{j}|=r} \varrho_{\mathbf{j}} \sum_{\mathbf{i} \leq \mathbf{j}, |\mathbf{i}| \geq 2} \frac{\langle \boldsymbol{\theta}^{\mathbf{i}} \rangle \langle (\hat{\boldsymbol{\theta}}'_n)^{\mathbf{j}-\mathbf{i}} \rangle}{\mathbf{i}!(\mathbf{j}-\mathbf{i})!} \\
&\quad + \sum_{r=1}^{k-1} n^{-r/2} \sum_{|\mathbf{j}|=r+1} \delta_{\mathbf{j}} \sum_{\mathbf{i} \leq \mathbf{j}, |\mathbf{i}| \geq 2} \frac{\langle \boldsymbol{\theta}^{\mathbf{i}} \rangle \langle (\hat{\boldsymbol{\theta}}'_n)^{\mathbf{j}-\mathbf{i}} \rangle}{\mathbf{i}!(\mathbf{j}-\mathbf{i})!} \\
&\quad + \sum_{r=1}^{k-1} n^{-r/2} \sum_{|\mathbf{j}|=r+2} \mathcal{E}_{\mathbf{j}} \sum_{\mathbf{i} \leq \mathbf{j}, |\mathbf{i}| \geq 2} \frac{\langle \boldsymbol{\theta}^{\mathbf{i}} \rangle \langle (\hat{\boldsymbol{\theta}}'_n)^{\mathbf{j}-\mathbf{i}} \rangle}{\mathbf{i}!(\mathbf{j}-\mathbf{i})!} + \sum_{|\mathbf{j}|=2} \mathcal{E}_{\mathbf{j}} \frac{\langle \boldsymbol{\theta}^{\mathbf{j}} \rangle}{\mathbf{j}!}.
\end{aligned}$$

In the above we used the fact that for $r=1$, $\sum_{|\mathbf{j}|=r} \varrho_{\mathbf{j}} \sum_{\mathbf{i} \leq \mathbf{j}, |\mathbf{i}| \geq 2} \cdots = 0$, as it is a summation over empty set, thus $\sum_{r=2}^{k-1} n^{-r/2} \sum_{|\mathbf{j}|=r} \varrho_{\mathbf{j}} \sum_{\mathbf{i} \leq \mathbf{j}, |\mathbf{i}| \geq 2} \cdots$ can be rewritten as $\sum_{r=1}^{k-1} n^{-r/2} \sum_{|\mathbf{j}|=r} \varrho_{\mathbf{j}} \sum_{\mathbf{i} \leq \mathbf{j}, |\mathbf{i}| \geq 2} \cdots$. Note

$$\sum_{|\mathbf{j}|=2} \mathcal{E}_{\mathbf{j}} \frac{\langle \boldsymbol{\theta}^{\mathbf{j}} \rangle}{\mathbf{j}!} = -\frac{1}{2} \boldsymbol{\theta}' \mathbf{I} \boldsymbol{\theta}.$$

Also, as in the proof of Theorem 2.2,

$$\langle (\hat{\boldsymbol{\theta}}'_n)^{\mathbf{j}-\mathbf{i}} \rangle \stackrel{k}{\sim} \langle (\hat{\boldsymbol{\theta}}''_n)^{\mathbf{j}-\mathbf{i}} \rangle \stackrel{k}{\sim} (\mathbf{j}-\mathbf{i})! \sum_{r=0}^{k-1} n^{-r/2} \sum_{|\mathbf{l}|=r} \sum_{(0,r,\mathbf{l},\mathbf{j}-\mathbf{i})} \prod_{v=0}^{r} \langle \mathbf{H}_v^{\mathbf{i}_v} \rangle / \mathbf{i}_v!.$$

Plugging in the above into (A.6) and rearranging terms, we get

$$\log \frac{Z_n(\boldsymbol{\theta} + \hat{\boldsymbol{\theta}}'_n)}{Z_n(\hat{\boldsymbol{\theta}}'_n)}$$

$$\overset{k}{\sim} -\boldsymbol{\theta}'\mathbf{I}\boldsymbol{\theta}/2$$

$$+ \sum_{r=1}^{k-1} n^{-r/2} \Bigg( \sum_{t+s=r, t\geq 1} \sum_{|\mathbf{i}|=2}^{r} \frac{\langle \boldsymbol{\theta^i} \rangle}{\mathbf{i}!} \sum_{\mathbf{j} \geq \mathbf{i}, |\mathbf{j}|=t} \varrho_{\mathbf{j}} \sum_{|\mathbf{l}|=s} \sum_{(0,s,\mathbf{l},\mathbf{j}-\mathbf{i})} \prod_{v=0}^{s} \frac{\langle \mathbf{H}_v^{\mathbf{i}_v} \rangle}{\mathbf{i}_v!}$$

$$+ \sum_{t+s=r+1, t\geq 2} \sum_{|\mathbf{i}|=2}^{r+1} \frac{\langle \boldsymbol{\theta^i} \rangle}{\mathbf{i}!} \sum_{\mathbf{j} \geq \mathbf{i}, |\mathbf{j}|=t} \delta_{\mathbf{j}} \sum_{|\mathbf{l}|=s} \sum_{(0,s,\mathbf{l},\mathbf{j}-\mathbf{i})} \prod_{v=0}^{s} \frac{\langle \mathbf{H}_v^{\mathbf{i}_v} \rangle}{\mathbf{i}_v!}$$

$$+ \sum_{t+s=r+2, t\geq 3} \sum_{|\mathbf{i}|=2}^{r+2} \frac{\langle \boldsymbol{\theta^i} \rangle}{\mathbf{i}!} \sum_{\mathbf{j} \geq \mathbf{i}, |\mathbf{j}|=t} \mathcal{E}_{\mathbf{j}} \sum_{|\mathbf{l}|=s} \sum_{(0,s,\mathbf{l},\mathbf{j}-\mathbf{i})} \prod_{v=0}^{s} \frac{\langle \mathbf{H}_v^{\mathbf{i}_v} \rangle}{\mathbf{i}_v!} \Bigg)$$

$$= -\boldsymbol{\theta}'\mathbf{I}\boldsymbol{\theta}/2$$

$$+ \sum_{r=1}^{k-1} n^{-r/2} \sum_{|\mathbf{i}|=2}^{r+2} \frac{\langle \boldsymbol{\theta^i} \rangle}{\mathbf{i}!} \Bigg( \sum_{t+s=r, t\geq 1} \sum_{\mathbf{j} \geq \mathbf{i}, |\mathbf{j}|=t} \varrho_{\mathbf{j}} \sum_{|\mathbf{l}|=s} \sum_{(0,s,\mathbf{l},\mathbf{j}-\mathbf{i})} \prod_{v=0}^{s} \frac{\langle \mathbf{H}_v^{\mathbf{i}_v} \rangle}{\mathbf{i}_v!}$$

$$+ \sum_{t+s=r+1, t\geq 2} \sum_{\mathbf{j} \geq \mathbf{i}, |\mathbf{j}|=t} \delta_{\mathbf{j}} \sum_{|\mathbf{l}|=s} \sum_{(0,s,\mathbf{l},\mathbf{j}-\mathbf{i})} \prod_{v=0}^{s} \frac{\langle \mathbf{H}_v^{\mathbf{i}_v} \rangle}{\mathbf{i}_v!}$$

$$+ \sum_{t+s=r+2, t\geq 3} \sum_{\mathbf{j} \geq \mathbf{i}, |\mathbf{j}|=t} \mathcal{E}_{\mathbf{j}} \sum_{|\mathbf{l}|=s} \sum_{(0,s,\mathbf{l},\mathbf{j}-\mathbf{i})} \prod_{v=0}^{s} \frac{\langle \mathbf{H}_v^{\mathbf{i}_v} \rangle}{\mathbf{i}_v!} \Bigg)$$

$$= -\frac{1}{2}\boldsymbol{\theta}'\mathbf{I}\boldsymbol{\theta} + \sum_{r=1}^{k-1} n^{-r/2} \sum_{|\mathbf{i}|=2}^{r+2} \frac{\langle \boldsymbol{\theta^i} \rangle}{\mathbf{i}!} F_{\mathbf{i},r}.$$

In the above we used the fact that, for $|\mathbf{i}| \geq r+1$, the summation $\sum_{t+s=r} \sum_{\mathbf{j} \geq \mathbf{i}, |\mathbf{j}|=t} \sum_{|\mathbf{l}|=s}$ is over an empty set, thus the first term inside the bracket above originally is $\sum_{r=1}^{k-1} n^{-r/2} \sum_{t+s=r} \sum_{|\mathbf{i}|=2}^{t} \langle \boldsymbol{\theta^i} \rangle / \mathbf{i}! \cdots$ and can be rewritten as $\sum_{r=1}^{k-1} n^{-r/2} \sum_{|\mathbf{i}|=2}^{r+2} \langle \boldsymbol{\theta^i} \rangle / \mathbf{i}! \cdots$. For the same reason, the second term inside the bracket above originally is $\sum_{r=1}^{k-1} n^{-r/2} \sum_{t+s=r+1} \sum_{|\mathbf{i}|=2}^{t} \langle \boldsymbol{\theta^i} \rangle / \mathbf{i}! \cdots$ and can be rewritten as $\sum_{r=1}^{k-1} n^{-r/2} \sum_{|\mathbf{i}|=2}^{r+2} \langle \boldsymbol{\theta^i} \rangle / \mathbf{i}! \cdots$. Now as in [17] we have

$$\frac{Z_n(\boldsymbol{\theta} + \hat{\boldsymbol{\theta}}'_n)}{Z_n(\hat{\boldsymbol{\theta}}'_n)} \overset{k}{\sim} \exp\left(-\frac{1}{2}\boldsymbol{\theta}'\mathbf{I}\boldsymbol{\theta}\right) \exp\left(\sum_{r=1}^{k-1} n^{-r/2} \sum_{|\mathbf{i}|=2}^{r+2} \frac{\langle \boldsymbol{\theta^i} \rangle}{\mathbf{i}!} F_{\mathbf{i},r}\right),$$

$$\exp\left(\sum_{r=1}^{k-1} n^{-r/2} \sum_{|\mathbf{i}|=2}^{r+2} \frac{\langle\boldsymbol{\theta}^{\mathbf{i}}\rangle}{\mathbf{i}!} F_{\mathbf{i},r}\right) \stackrel{k}{\sim} 1 + \sum_{r=1}^{k-1} n^{-r/2} \sum_{I_0(1,r,r)} \prod_{v=1}^{r} \left(\sum_{|\mathbf{i}|=2}^{v+2} \frac{\langle\boldsymbol{\theta}^{\mathbf{i}}\rangle}{\mathbf{i}!} F_{\mathbf{i},v}\right)^{i_v} \Big/ i_v!,$$

the second summation on the right-hand side above is for $(i_1, \ldots, i_r) \in I_0(1, r, r)$. Also,

$$\sum_{I_0(1,r,r)} \prod_{v=1}^{r} \left(\sum_{|\mathbf{i}|=2}^{v+2} \frac{\langle\boldsymbol{\theta}^{\mathbf{i}}\rangle}{\mathbf{i}!} F_{\mathbf{i},v}\right)^{i_v} \Big/ i_v!$$

$$= \sum_{I_0(1,r,r)} \prod_{v=1}^{r} \sum_{|\mathbf{k}_v|=2i_v}^{(v+2)i_v} \langle\boldsymbol{\theta}^{\mathbf{k}_v}\rangle \sum_{I_1(2,v+2,\mathbf{k}_v,i_v)} \prod_{|\mathbf{j}|=2}^{v+2} \frac{F_{\mathbf{j},v}^{u_\mathbf{j}}}{u_\mathbf{j}!(\mathbf{j}!)^{u_\mathbf{j}}} = \sum_{|\mathbf{i}|=2}^{3r} \langle\boldsymbol{\theta}^{\mathbf{i}}\rangle L_{\mathbf{i},r};$$

in the left-hand side above the summations are for $(i_1, \ldots, i_r) \in I_0(1, r, r)$, and for each given $v$ and $\mathbf{k}_v$, $u_\mathbf{j} \in I_1(2, v+2, \mathbf{k}_v, i_v)$. Now we have

$$\exp\left(\sum_{r=1}^{k-1} n^{-r/2} \sum_{|\mathbf{i}|=2}^{r+2} \frac{\langle\boldsymbol{\theta}^{\mathbf{i}}\rangle}{\mathbf{i}!} F_{\mathbf{i},r}\right) \stackrel{k}{\sim} 1 + \sum_{r=1}^{k-1} n^{-r/2} \sum_{|\mathbf{i}|=2}^{3r} \langle\boldsymbol{\theta}^{\mathbf{i}}\rangle N_{\mathbf{i},r}.$$

In the definition of $L_{\mathbf{i},r}$, $I_0(1,r,r) = \bigcup_{l \geq 0} I_1(1,r,r,l)$, where $I_1(1,r,r,l)$ is defined in the proof of Theorem 2.2. For given $(i_1, \ldots, i_r) \in I_0(1, r, r)$, and integer $d$-vector $\mathbf{i}$, $I_2(r, \mathbf{i})$ is the collection of integer $d$-vectors $\mathbf{k}_1, \ldots, \mathbf{k}_r$,

$$I_2(r, \mathbf{i}) = \left\{(\mathbf{k}_1, \ldots, \mathbf{k}_r) : \sum_{v=1}^{r} |\mathbf{k}_v| = \mathbf{i},\ 2i_v \leq |\mathbf{k}_v| \leq (v+2)i_v,\ v = 1, \ldots, r\right\}.$$

Given integer $d$-vector $\mathbf{k}$, integers $i$ and $m \leq r$, $I_1(m, r, \mathbf{k}, i)$ is the collection of integers $u_\mathbf{j}$ indexed by a integer $d$-vector $\mathbf{j}$,

$$I_1(m, r, \mathbf{k}, i) = \left\{u_\mathbf{j} : \sum_{|\mathbf{j}|=m}^{r} \mathbf{j} u_\mathbf{j} = \mathbf{k},\ \sum_{|\mathbf{j}|=m}^{r} u_\mathbf{j} = i\right\}. \qquad \square$$

PROOF OF THEOREM 2.3. By Theorem 2.2, we only need to prove

(A.7) $$\sqrt{n}(\check{\boldsymbol{\theta}}_n - \hat{\boldsymbol{\theta}}_n) = \sum_{r=1}^{k-1} n^{-r/2} \mathbf{Q}_r + O_P(n^{-k/2}).$$

Denote $\mathbf{W}^{(1)}(\boldsymbol{\theta}) = (\partial W(\boldsymbol{\theta})/\partial \theta_1, \ldots, \partial W(\boldsymbol{\theta})/\partial \theta_d)'$, $\mathbf{d}_n = \sqrt{n}(\check{\boldsymbol{\theta}}_n - \hat{\boldsymbol{\theta}}_n)$ and $\hat{\boldsymbol{\theta}}'_n = \sqrt{n}(\hat{\boldsymbol{\theta}}_n - \boldsymbol{\theta}_0)$. We only need to point out the main modifications to the proof of Theorem 4 in [17]. In place (4.5) of [17] there we have

$$\int \mathbf{W}^{(1)}\left(\frac{\mathbf{d}_n - \boldsymbol{\theta}}{\sqrt{n}}\right) \frac{Z_n(\boldsymbol{\theta} + \hat{\boldsymbol{\theta}}'_n)}{Z_n(\hat{\boldsymbol{\theta}}'_n)} d\boldsymbol{\theta} = \mathbf{0},$$

corresponding to (4.7) of [17] there we have

$$\|\mathbf{W}^{(1)}(\mathbf{u})\| \leq \left(\sum_{i=1}^{d}(W(\mathbf{u}+\mathbf{e}_i)+W(\mathbf{u})+W(\mathbf{u}-\mathbf{e}_i))^2\right)^{1/2} \quad \text{a.s.}$$

By Condition 9, for $\|(\mathbf{d}_n - \boldsymbol{\theta})/\sqrt{n}\| \leq \varepsilon_2$, we can replace $\mathbf{W}^{(1)}((\mathbf{d}_n - \boldsymbol{\theta})/\sqrt{n})$ by the vector $\mathbf{a}((\mathbf{d}_n - \boldsymbol{\theta})/\sqrt{n})^{\mathbf{a}-\mathbf{1}}$, and (4.10) of [17] there is replaced by, for $\delta > 0$,

$$\int_{\|\boldsymbol{\theta}\| \leq n^{\delta/2}} (\mathbf{d}_n - \boldsymbol{\theta})^{\mathbf{a}-\mathbf{1}} \frac{Z_n(\boldsymbol{\theta}+\hat{\boldsymbol{\theta}}'_n)}{Z_n(\hat{\boldsymbol{\theta}}'_n)} d\boldsymbol{\theta} \overset{k}{\sim} \mathbf{0}.$$

As in [17], $\mathbf{d}_n \overset{1}{\sim} \mathbf{0}$. Define $N_{\mathbf{0},0} = 1$, thus $1 = n^{-0/2} \langle \boldsymbol{\theta}^{\mathbf{0}} \rangle N_{\mathbf{0},0}$, let $|\mathbf{I}|$ be the determinant of $\mathbf{I}$, by Lemma 1 the above is

$$\mathbf{0} \overset{k}{\sim} \sum_{r=0}^{k-1} n^{-r/2} \sum_{|\mathbf{i}|=2(1 \wedge r)}^{3r} (2\pi)^{d/2} |\mathbf{I}|^{-1/2} N_{\mathbf{i},r}$$

$$\times \int_{\|\boldsymbol{\theta}\| \leq n^{\delta/2}} (\mathbf{d}_n - \boldsymbol{\theta})^{\mathbf{a}-\mathbf{1}} \langle \boldsymbol{\theta}^{\mathbf{i}} \rangle \frac{|\mathbf{I}|^{1/2}}{(2\pi)^{d/2}} \exp\left(-\frac{1}{2}\boldsymbol{\theta}'\mathbf{I}\boldsymbol{\theta}\right) d\boldsymbol{\theta}.$$

Also, as in [17], for each $\mathbf{i}$ and $k > 0$,

$$\int_{\|\boldsymbol{\theta}\| > n^{\delta/2}} \frac{|\mathbf{I}|^{1/2}}{(2\pi)^{d/2}} \exp\left(-\frac{1}{2}\boldsymbol{\theta}'\mathbf{I}\boldsymbol{\theta}\right) d\boldsymbol{\theta} \sim o_P(n^{-k/2}).$$

Define the $R^d$ to $R^d$ function $\boldsymbol{\Psi}_{\mathbf{i}}(\cdot) = (\psi_{\mathbf{i},1}(\cdot), \ldots, \psi_{\mathbf{i},d}(\cdot))'$ as

$$\boldsymbol{\Psi}_{\mathbf{i}}(\mathbf{u}) = \int \boldsymbol{\theta}^{\mathbf{a}-\mathbf{1}} \langle (\boldsymbol{\theta}+\mathbf{u})^{\mathbf{i}} \rangle \frac{|\mathbf{I}|^{1/2}}{(2\pi)^{d/2}} \exp\left(-\frac{1}{2}(\boldsymbol{\theta}+\mathbf{u})'\mathbf{I}(\boldsymbol{\theta}+\mathbf{u})\right) d\boldsymbol{\theta},$$

and $\boldsymbol{\Psi}_{\mathbf{i}}^{(\mathbf{j})}(\cdot) = (\psi_{\mathbf{i},1}^{(\mathbf{j})}(\cdot), \ldots, \psi_{\mathbf{i},d}^{(\mathbf{j})}(\cdot))'$ with $\psi_{\mathbf{i},k}^{(\mathbf{j})}(\mathbf{u}) = \partial^{|\mathbf{j}|}\psi_{\mathbf{i},k}(\mathbf{u})/\partial \mathbf{u}^{\mathbf{j}}$ ($k=1,\ldots,d$).

On the right-hand side of the "$\mathbf{0} \overset{k}{\sim}$" relationship above, leave out the constant $(2\pi)^{d/2}|\mathbf{I}|^{-1/2}$ and multiply the nonsingular matrix $(\{\boldsymbol{\sigma}(\mathbf{a})\}\mathbf{I})^{-1}$, the "$\mathbf{0} \overset{k}{\sim}$" relationship remains. The choice of this matrix will be clear when we prove the expression for $\mathbf{Q}_r$ $(1 \leq r \leq k-1)$ later. We will show below that $\boldsymbol{\Psi}_{\mathbf{i}}^{(\mathbf{j})}(\mathbf{0}) = \mathbf{0}$, for $|\mathbf{i}+\mathbf{j}|$ even, so the previous relationship is rewritten as

$$\mathbf{0} \overset{k}{\sim} \sum_{r=0}^{k-1} n^{-r/2} \sum_{|\mathbf{i}|=2(1 \wedge r)}^{3r} N_{\mathbf{i},r} (\{\boldsymbol{\sigma}(\mathbf{a})\}\mathbf{I})^{-1} \boldsymbol{\Psi}_{\mathbf{i}}(\mathbf{d}_n)$$

$$(\text{A.8}) \overset{k}{\sim} \sum_{r=0}^{k-1} n^{-r/2} \sum_{|\mathbf{i}|=2(1 \wedge r)}^{3r} N_{\mathbf{i},r} \sum_{|\mathbf{j}|=0}^{k-1-r} \frac{\langle \mathbf{d}_n^{\mathbf{j}} \rangle}{\mathbf{j}!} (\{\boldsymbol{\sigma}(\mathbf{a})\}\mathbf{I})^{-1} \boldsymbol{\Psi}_{\mathbf{i}}^{(\mathbf{j})}(\mathbf{0})$$

$$= \sum_{r=0}^{k-1} n^{-r/2} \sum_{|\mathbf{j}|=0}^{k-1-r} \langle \mathbf{d}_n^{\mathbf{j}} \rangle \sum_{|\mathbf{i}|=2(1\wedge r)}^{3r} \frac{1}{\mathbf{j}!} N_{\mathbf{i},r}(\{\boldsymbol{\sigma}(\mathbf{a})\}\mathbf{I})^{-1}\boldsymbol{\Psi}_{\mathbf{i}}^{(\mathbf{j})}(\mathbf{0})$$

$$= \sum_{r=0}^{k-1} n^{-r/2} \sum_{|\mathbf{j}|=0}^{k-1-r} \langle \mathbf{d}_n^{\mathbf{j}} \rangle \sum_{|\mathbf{s}|=2(1\wedge r)+|\mathbf{j}|}^{3r+|\mathbf{j}|} \frac{1}{\mathbf{j}!} N_{\mathbf{s-j},r}(\{\boldsymbol{\sigma}(\mathbf{a})\}\mathbf{I})^{-1}\boldsymbol{\Psi}_{\mathbf{s-j}}^{(\mathbf{j})}(\mathbf{0})$$

$$= \sum_{r=0}^{k-1} n^{-r/2} \sum_{|\mathbf{j}|=0}^{k-1-r} \langle \mathbf{d}_n^{\mathbf{j}} \rangle \sum_{|\mathbf{s}|\in\langle 2(1\wedge r)+|\mathbf{j}|, 3r+|\mathbf{j}|\rangle} \frac{1}{\mathbf{j}!} N_{\mathbf{s-j},r}(\{\boldsymbol{\sigma}(\mathbf{a})\}\mathbf{I})^{-1}\boldsymbol{\Psi}_{\mathbf{s-j}}^{(\mathbf{j})}(\mathbf{0})$$

$$= \sum_{r=0}^{k-1} n^{-r/2} \sum_{|\mathbf{j}|=0}^{k-1-r} \langle \mathbf{d}_n^{\mathbf{j}} \rangle \mathbf{M}_{\mathbf{j},r},$$

where $\{\boldsymbol{\sigma}(\mathbf{a})\}$ is the $d \times d$ diagonal matrix with $r$th element be the $a_r$th moment for $\theta_r$ [with $\boldsymbol{\theta} \sim N(\mathbf{0}, \mathbf{I}^{-1})$]. Denote $\mathbf{I} = (i_{st})$, $\mathbf{I}^{-1} = (i^{st})$, using the joint moment formula for multivariate normal distribution, we have

$$\psi_{\mathbf{i},k}(\mathbf{0}) = \frac{|\mathbf{I}|^{1/2}}{(2\pi)^{d/2}} \int \theta_k^{a_k-1} \langle \boldsymbol{\theta}^{\mathbf{i}} \rangle \exp\left(-\frac{1}{2}\boldsymbol{\theta}'\mathbf{I}\boldsymbol{\theta}\right) d\boldsymbol{\theta}$$

$$= \begin{cases} 0, & \text{if } |\mathbf{i}| \text{ even,} \\ \sigma((\mathbf{a}-\mathbf{1})\mathbf{e}_k+\mathbf{i}), & \text{otherwise,} \end{cases}$$

where $\sigma(\mathbf{a}) = E\langle \boldsymbol{\theta}^{\mathbf{a}} \rangle$ is the joint $\mathbf{a}$th moment of $\boldsymbol{\theta} \sim N(\mathbf{0}, \mathbf{I}^{-1})$. Note

$$\frac{\partial(\boldsymbol{\theta}+\mathbf{u})'\mathbf{I}(\boldsymbol{\theta}+\mathbf{u})}{2\partial u_k}\bigg|_{\mathbf{u}=\mathbf{0}} = \mathbf{e}_k'\mathbf{I}\boldsymbol{\theta} = \sum_{r=1}^{d} i_{kr}\theta_r,$$

$$\frac{\partial^2(\boldsymbol{\theta}+\mathbf{u})'\mathbf{I}(\boldsymbol{\theta}+\mathbf{u})}{2\partial u_k\,\partial u_r}\bigg|_{\mathbf{u}=\mathbf{0}} = \mathbf{e}_k'\mathbf{I}\mathbf{e}_r = i_{kr}$$

and $\partial^{|\mathbf{i}|}[(\boldsymbol{\theta}+\mathbf{u})'\mathbf{I}(\boldsymbol{\theta}+\mathbf{u})]/\partial \mathbf{u}^{\mathbf{i}}|_{\mathbf{u}=\mathbf{0}} = 0$, for $|\mathbf{i}| > 2$, we have

$$\psi_{\mathbf{i},k}^{(\mathbf{j})}(\mathbf{0}) = \frac{|\mathbf{I}|^{1/2}}{(2\pi)^{d/2}} \int P_{\mathbf{i}+(\mathbf{a}-\mathbf{1})\mathbf{e}_k}^{(\mathbf{j})}(\boldsymbol{\theta}) \exp\left(-\frac{1}{2}\boldsymbol{\theta}'\mathbf{I}\boldsymbol{\theta}\right) d\boldsymbol{\theta} = P_{\mathbf{i}+(\mathbf{a}-\mathbf{1})\mathbf{e}_k}^{(\mathbf{j})}(\sigma),$$

where $P_{\mathbf{i}+(\mathbf{a}-\mathbf{1})\mathbf{e}_k}^{(\mathbf{j})}(\boldsymbol{\theta})$ is the multivariate polynomial in $\boldsymbol{\theta}$ given by the relationship [note $\theta_k^{a_k-1} = \langle \boldsymbol{\theta}^{(\mathbf{a}-\mathbf{1})\mathbf{e}_k} \rangle$]

$$\frac{\partial^{|\mathbf{j}|}[\theta_k^{a_k-1}\langle(\boldsymbol{\theta}+\mathbf{u})^{\mathbf{i}}\rangle \exp(-1/2(\boldsymbol{\theta}+\mathbf{u})'\mathbf{I}(\boldsymbol{\theta}+\mathbf{u}))]}{\partial \mathbf{u}^{\mathbf{j}}}\bigg|_{\mathbf{u}=\mathbf{0}}$$

$$= P_{\mathbf{i}+(\mathbf{a}-\mathbf{1})\mathbf{e}_k}^{(\mathbf{j})}(\boldsymbol{\theta}) \exp\left(-\frac{1}{2}\boldsymbol{\theta}'\mathbf{I}\boldsymbol{\theta}\right),$$

and $P_{\mathbf{i}+(\mathbf{a}-\mathbf{1})\mathbf{e}_k}^{(\mathbf{j})}(\sigma) = E(P_{\mathbf{i}+(\mathbf{a}-\mathbf{1})\mathbf{e}_k}^{(\mathbf{j})}(\boldsymbol{\theta}))$ with $\boldsymbol{\theta} \sim N(\mathbf{0}, \mathbf{I}^{-1})$ is the corresponding vector polynomial in the moments $\sigma_{st}$'s. Especially, $\psi_{\mathbf{0},k}^{(\mathbf{j})}(\mathbf{0}) = 0$ ($k =$ $1, \ldots, d)$, for $|\mathbf{j}|$ even, or $\boldsymbol{\Psi}^{(\mathbf{j})}_{\mathbf{0},k}(\mathbf{0}) = \mathbf{0}$, for $|\mathbf{j}|$ even, and it can be shown that $\boldsymbol{\Psi}^{(\mathbf{j})}_{\mathbf{i},k}(\mathbf{0}) = \mathbf{0}$, if $|\mathbf{i}+\mathbf{j}|$ is even. We get

$$\mathbf{M}_{\mathbf{j},r} = (\{\boldsymbol{\sigma}(\mathbf{a})\}\mathbf{I})^{-1} \sum_{|\mathbf{i}|=2(1\wedge r)}^{3r} N_{\mathbf{i},r} \frac{1}{\mathbf{j}!} \boldsymbol{\Psi}^{(\mathbf{j})}_{\mathbf{i}}(\mathbf{0})$$

and $\mathbf{M}_{\mathbf{j},r} = \mathbf{0}$ for $|\mathbf{i}+\mathbf{j}|$ even. Note the recursive relationship for the $\mathbf{M}_{\mathbf{i},r}$'s in Theorem 2.3 can be rewritten as

$$(A.9) \quad \sum_{m=0}^{r} \sum_{|\mathbf{i}|=0}^{m} \mathbf{i}! \mathbf{M}_{\mathbf{i},r-m} \sum_{|\mathbf{l}|=m} \sum_{(1,m,\mathbf{l},\mathbf{i})} \prod_{v=1}^{m} \langle \mathbf{Q}^{\mathbf{i}_v}_v \rangle / \mathbf{i}_v! = \mathbf{0}, \quad 1 \le r \le k-1.$$

To see this, let $\mathbf{d}'_n = \sum_{j=1}^{k-1} n^{-j/2} \mathbf{Q}_j$; as in the proof of Theorem 2.2, we have

$$\langle (\mathbf{d}'_n)^{\mathbf{i}} \rangle \overset{k}{\sim} \mathbf{i}! \sum_{r=|\mathbf{i}|}^{k-1} n^{-r/2} \sum_{|\mathbf{l}|=r} \sum_{(1,r,\mathbf{l},\mathbf{i})} \prod_{v=1}^{r} \langle \mathbf{Q}^{\mathbf{i}_v}_v \rangle / \mathbf{i}_v!.$$

Denote $\mathbf{P}^{(\mathbf{j})}_{\mathbf{i}+\mathbf{a}-\mathbf{1}}(\sigma) = (P^{(\mathbf{j})}_{\mathbf{i}+(\mathbf{a}-\mathbf{1})\mathbf{e}_1}(\sigma), \ldots, P^{(\mathbf{j})}_{\mathbf{i}+(\mathbf{a}-\mathbf{1})\mathbf{e}_d}(\sigma))'$. Take $\mathbf{i}=\mathbf{0}$ in the above, we have $\sum_{|\mathbf{l}|=0} \sum_{(1,0,\mathbf{l},\mathbf{0})} \prod_{v=1} \langle \mathbf{Q}^{\mathbf{i}_v}_v \rangle / \mathbf{i}_v! = 1$ symbolically. Note $\mathbf{M}_{\mathbf{0},0} = \mathbf{0}$, and recall we defined $N_{\mathbf{0},0} = 1$, so $\mathbf{M}_{\mathbf{e}_r,0} = N_{\mathbf{0},0}(\{\boldsymbol{\sigma}(\mathbf{a})\mathbf{I}\})^{-1} \mathbf{P}^{(\mathbf{e}_r)}_{\mathbf{a}-\mathbf{1}}(\sigma) = -(\{\boldsymbol{\sigma}(\mathbf{a})\}\mathbf{I})^{-1} \times (\{\boldsymbol{\sigma}(\mathbf{a})\}\mathbf{I})_r$, where $(\{\boldsymbol{\sigma}(\mathbf{a})\}\mathbf{I})_r$ is the $r$th column of $\{\boldsymbol{\sigma}(\mathbf{a})\}\mathbf{I}$; this is the reason we multiply it in the previous $\overset{k}{\sim}$ relationship, and we see that $\{\boldsymbol{\sigma}(\mathbf{a})\}\mathbf{I} = -(\mathbf{P}^{(\mathbf{e}_1)}_{\mathbf{a}-\mathbf{1}}(\sigma), \ldots, \mathbf{P}^{(\mathbf{e}_d)}_{\mathbf{a}-\mathbf{1}}(\sigma))$. Note $\sum_{r=1}^{k}(\{\boldsymbol{\sigma}(\mathbf{a})\}\mathbf{I})_r \langle \mathbf{Q}^{\mathbf{e}_r}_1 \rangle = (\{\boldsymbol{\sigma}(\mathbf{a})\}\mathbf{I})\mathbf{Q}_1$. So, for $r=1$, (A.9) is

$$\mathbf{0} = \mathbf{M}_{\mathbf{0},1} + \sum_{|\mathbf{i}|=0}^{1} \mathbf{M}_{\mathbf{i},0} \sum_{|\mathbf{l}|=1} \sum_{(1,1,\mathbf{l},\mathbf{i})} \langle \mathbf{Q}^{\mathbf{i}}_1 \rangle = \mathbf{M}_{\mathbf{0},1} + \sum_{|\mathbf{i}|=0}^{1} \mathbf{M}_{\mathbf{i},0} \langle \mathbf{Q}^{\mathbf{i}}_1 \rangle$$

$$= \mathbf{M}_{\mathbf{0},1} + \mathbf{M}_{\mathbf{0},0} + \sum_{r=1}^{d} \mathbf{M}_{\mathbf{e}_r,0} \langle \mathbf{Q}^{\mathbf{e}_r}_1 \rangle = \mathbf{M}_{\mathbf{0},1} - \mathbf{Q}_1;$$

this gives the equivalence of (A.9) and the expression for $\mathbf{Q}_1$ given in Theorem 2.3.

For $r>1$, since $\mathbf{M}_{\mathbf{j},0} = \mathbf{0}$ for $|\mathbf{j}|$ even, (A.9) is rewritten as

$$\mathbf{0} = \sum_{m=0}^{r-1} \sum_{|\mathbf{i}|=0}^{m} \mathbf{i}! \mathbf{M}_{\mathbf{i},r-m} \sum_{|\mathbf{l}|=m} \sum_{(1,m,\mathbf{l},\mathbf{i})} \prod_{v=1}^{m} \langle \mathbf{Q}^{\mathbf{i}_v}_v \rangle / \mathbf{i}_v!$$

$$+ \sum_{|\mathbf{i}|\in\langle 3,r\rangle} \mathbf{i}! \mathbf{M}_{\mathbf{i},0} \sum_{|\mathbf{l}|=r} \sum_{(1,r,\mathbf{l},\mathbf{i})} \prod_{v=1}^{r} \langle \mathbf{Q}^{\mathbf{i}_v}_v \rangle / \mathbf{i}_v! + \sum_{|\mathbf{i}|=1} \mathbf{M}_{\mathbf{i},0} \sum_{|\mathbf{l}|=r} \sum_{(1,r,\mathbf{l},\mathbf{i})} \prod_{v=1}^{r} \langle \mathbf{Q}^{\mathbf{i}_v}_v \rangle / \mathbf{i}_v!.$$

As before, $\sum_{|\mathbf{l}|=0}\sum_{(1,0\mathbf{l},\mathbf{0})}\prod_{v=1}\langle\mathbf{Q}_v^{\mathbf{i}_v}\rangle/\mathbf{i}_v! = 1$, so the first term above is $\mathbf{M}_{\mathbf{0},r} + \sum_{m=1}^{r-1}\sum_{|\mathbf{i}|=1}^{m}\mathbf{i}!\mathbf{M}_{\mathbf{i},r-m}\sum_{|\mathbf{l}|=m}\sum_{(1,m,\mathbf{l},\mathbf{i})}\prod_{v=1}^{m}\langle\mathbf{Q}_v^{\mathbf{i}_v}\rangle/\mathbf{i}_v!$. For the second term, since $r > 1$, $3 \leq |\mathbf{i}| \leq |\mathbf{l}|$, the set $(1,r,\mathbf{l},\mathbf{i})$ is empty if $\mathbf{i}_r \neq \mathbf{0}$, thus the factor $\prod_{v=1}^{r}\langle\mathbf{Q}_v^{\mathbf{i}_v}\rangle/\mathbf{i}_v!$ in the second term is in fact $\prod_{v=1}^{r-1}\langle\mathbf{Q}_v^{\mathbf{i}_v}\rangle/\mathbf{i}_v!$. For the third term, $\mathbf{i} \in \{\mathbf{e}_1,\ldots,\mathbf{e}_d\}$ and $|\mathbf{l}| = r$, if $\mathbf{i} = \mathbf{e}_k$, the set $(1,r,\mathbf{l},\mathbf{e}_k)$ is nonempty only if $\mathbf{l} = r\mathbf{e}_k$, and in this case $(1,r,\mathbf{l},\mathbf{i}) = \{(\mathbf{i}_1,\ldots,\mathbf{i}_{r-1},\mathbf{i}_r) = (\mathbf{0},\ldots,\mathbf{0},\mathbf{e}_k)\}$. Similarly as before, the third term above is $\sum_{k=1}^{d}\mathbf{M}_{\mathbf{e}_k,0}\langle\mathbf{Q}_r^{\mathbf{e}_k}\rangle = -(\{\boldsymbol{\sigma}(\mathbf{a})\}\mathbf{I})^{-1}(\{\boldsymbol{\sigma}(\mathbf{a})\}\mathbf{I})_r\langle\mathbf{Q}_r^{\mathbf{e}_k}\rangle = -\mathbf{Q}_r$. Now we get

$$\mathbf{0} = \mathbf{M}_{\mathbf{0},r} + \sum_{m=1}^{r-1}\sum_{|\mathbf{i}|=1}^{m}\mathbf{i}!\mathbf{M}_{\mathbf{i},r-m}\sum_{|\mathbf{l}|=m}\sum_{(1,m,\mathbf{l},\mathbf{i})}\prod_{v=1}^{m}\langle\mathbf{Q}_v^{\mathbf{i}_v}\rangle/\mathbf{i}_v!$$

$$+ \sum_{|\mathbf{i}|\in\langle 3,r\rangle}\mathbf{i}!\mathbf{M}_{\mathbf{i},0}\sum_{|\mathbf{l}|=r}\sum_{(1,r-1,\mathbf{l},\mathbf{i})}\prod_{v=1}^{r-1}\langle\mathbf{Q}_v^{\mathbf{i}_v}\rangle/\mathbf{i}_v! - \mathbf{Q}_r,$$

and the equivalence of (A.9) and the recursive relationship for the $\mathbf{M}_{\mathbf{i},r}$'s in Theorem 2.3 is proved.

We now show that $\mathbf{d}_n'$ also satisfies (A.8). In fact, by (A.9) we have

$$\sum_{r=0}^{k-1}n^{-r/2}\sum_{|\mathbf{i}|=0}^{k-1-r}\langle(\mathbf{d}_n')^{\mathbf{i}}\rangle\mathbf{M}_{\mathbf{i},r}$$

$$\overset{k}{\sim}\sum_{r=0}^{k-1}n^{-r/2}\sum_{|\mathbf{i}|=0}^{k-1-r}\mathbf{M}_{\mathbf{i},r}\left(\mathbf{i}!\sum_{s=|\mathbf{i}|}^{k-1}n^{-s/2}\sum_{|\mathbf{l}|=s}\sum_{(1,s,\mathbf{l},\mathbf{i})}\prod_{v=1}^{s}\langle\mathbf{Q}_v^{\mathbf{i}_v}\rangle/\mathbf{i}_v!\right)$$

$$\overset{k}{\sim}\sum_{r=0}^{k-1}n^{-r/2}\sum_{s=0}^{k-1}n^{-s/2}\sum_{|\mathbf{i}|=s}^{k-1-r}\mathbf{i}!\mathbf{M}_{\mathbf{i},r}\sum_{|\mathbf{l}|=s}\sum_{(1,s,\mathbf{l},\mathbf{i})}\prod_{v=1}^{s}\langle\mathbf{Q}_v^{\mathbf{i}_v}\rangle/\mathbf{i}_v!$$

$$\overset{k}{\sim}\sum_{r=0}^{k-1}n^{-t/2}\sum_{r+s=t}\sum_{|\mathbf{i}|=0}^{(k-1-r)\vee s}\mathbf{i}!\mathbf{M}_{\mathbf{i},r}\sum_{|\mathbf{l}|=s}\sum_{(1,s,\mathbf{l},\mathbf{i})}\prod_{v=1}^{s}\langle\mathbf{Q}_v^{\mathbf{i}_v}\rangle/\mathbf{i}_v!.$$

Note $t \leq k-1$ and $r+s=t$ implies $k-1-r \geq s$, so $(k-1-r)\vee s = s$, and the above is

$$\sum_{r=0}^{k-1}n^{-t/2}\sum_{s=0}^{t}\sum_{|\mathbf{i}|=0}^{s}\mathbf{i}!\mathbf{M}_{\mathbf{i},t-s}\sum_{|\mathbf{l}|=s}\sum_{(1,s,\mathbf{l},\mathbf{i})}\prod_{v=1}^{s}\langle\mathbf{Q}_v^{\mathbf{i}_v}\rangle/\mathbf{i}_v! = \mathbf{0}.$$

The remaining proofs are similar to those in [17] and are omitted.

Lastly, for general loss function $W(\cdot)$ satisfies Conditions (B8) and (B10) with various derivatives and $\|\mathbf{W}^{(1)}(n^{-1/2}(\mathbf{d}_n - \boldsymbol{\theta}))\| \leq C\|\boldsymbol{\theta} + \hat{\boldsymbol{\theta}}_n'\|^{\gamma}$ for some $0 < C < \infty$, and $\gamma > 0$ [17]. Then similarly as before, we have

$$\boldsymbol{\Psi}_{\mathbf{i}}(\mathbf{u}) = \int \mathbf{W}^{(1)}(\boldsymbol{\theta})\langle(\boldsymbol{\theta}+\mathbf{u})^{\mathbf{i}}\rangle\frac{|\mathbf{I}|^{1/2}}{(2\pi)^{d/2}}\exp\left(-\frac{1}{2}(\boldsymbol{\theta}+\mathbf{u})'\mathbf{I}(\boldsymbol{\theta}+\mathbf{u})\right)d\boldsymbol{\theta},$$

$$\mathbf{\Psi_i^{(j)}(0)} = \int \frac{\partial^{|\mathbf{j}|}}{\partial \mathbf{u^j}} \left[ \mathbf{W}^{(1)}(\boldsymbol{\theta}) \langle (\boldsymbol{\theta} + \mathbf{u})^{\mathbf{i}} \rangle \frac{|\mathbf{I}|^{1/2}}{(2\pi)^{d/2}} \exp\left(-\frac{1}{2}(\boldsymbol{\theta} + \mathbf{u})'\mathbf{I}(\boldsymbol{\theta} + \mathbf{u})\right) \right]\bigg|_{\mathbf{u}=\mathbf{0}} d\boldsymbol{\theta}.$$

Denote $\mathbf{W}^{(1)} = (W_1, \ldots, W_d)'$, assume $W_k(1) \neq 0$ $(k = 1, \ldots, d)$, and define
$$\boldsymbol{\Sigma} = (\tilde{\sigma}_{kj}),$$

$$\tilde{\sigma}_{kj} = -\frac{|\mathbf{I}|^{1/2}}{(2\pi)^{d/2}} \int \frac{\partial}{\partial u_j} \left[ \frac{W_k(\boldsymbol{\theta})}{W_k(1)} \exp\left(-\frac{1}{2}(\boldsymbol{\theta} + \mathbf{u})'\mathbf{I}(\boldsymbol{\theta} + \mathbf{u})\right) \right]\bigg|_{\mathbf{u}=\mathbf{0}} d\boldsymbol{\theta}.$$

[The previous situation is a special case with $\tilde{\sigma}_{kj} = E_{\boldsymbol{\theta}}(\theta_k^{a_k-1} \sum_{r=1}^{d} i_{rj}\theta_r) = i_{rk} E_{\boldsymbol{\theta}}(\theta_k^{a_k})$ and $\boldsymbol{\theta} \sim N(\mathbf{0}, \mathbf{I}^{-1})$, or $\boldsymbol{\Sigma} = \{\boldsymbol{\sigma}(\mathbf{a})\}\mathbf{I}$.] Then, we have $\mathbf{G}_0 = \mathbf{H}_0$, $\mathbf{G}_r = \mathbf{H}_r + \mathbf{Q}_r$ $(1 \leq r \leq k-1)$, with

$$\mathbf{Q}_r = \sum_{m=0}^{r-1} \sum_{|\mathbf{i}|=0}^{m} \mathbf{i}! \mathbf{M}_{\mathbf{i},r-m} \sum_{|\mathbf{l}|=m} \sum_{(1,m,\mathbf{l},\mathbf{i})} \prod_{v=1}^{m} \langle \mathbf{Q}_v^{\mathbf{i}_v} \rangle / \mathbf{i}_v!$$

$$+ \sum_{|\mathbf{i}|=2}^{r} \mathbf{i}! \mathbf{M}_{\mathbf{i},0} \sum_{|\mathbf{l}|=r} \sum_{(1,r-1,\mathbf{l},\mathbf{i})} \prod_{v=1}^{r-1} \langle \mathbf{Q}_v^{\mathbf{i}_v} \rangle / \mathbf{i}_v!,$$

$$\mathbf{M}_{\mathbf{j},r} = (\boldsymbol{\Sigma})^{-1} \sum_{|\mathbf{i}|=2(1 \wedge r)+|\mathbf{j}|}^{3r+|\mathbf{j}|} \frac{1}{\mathbf{j}!} N_{\mathbf{i}-\mathbf{j},r} \mathbf{\Psi}_{\mathbf{i}-\mathbf{j}}^{(\mathbf{i})}(\mathbf{0}).$$

In the above, we assume $\boldsymbol{\Sigma}$ to be nonsingular. Here we no longer have $\mathbf{\Psi}_{\mathbf{i}}^{(\mathbf{j})}(\mathbf{0}) = 0$ for $|\mathbf{i} + \mathbf{j}|$ even, thus the above formula for $\mathbf{M}_{\mathbf{j},r}$ has more terms than stated in Theorem 2.3, and each term has more complicated form. $\square$

PROOF OF THEOREM 2.4. For $\boldsymbol{\theta}_n = (\check{\boldsymbol{\alpha}}_n, \hat{\boldsymbol{\beta}}_n)$ given in (2.1), although it can be formulated as a joint Bayesian estimator with an additional 0–1 error loss and a constant prior on $\boldsymbol{\beta}$, Theorem 2.3 cannot be applied to get its expansion, as Condition (B10) there excludes such loss.

We first outline the idea of the proof. Denote $\mathbf{H}_r^\circ = (\mathbf{h}_{r1}', \mathbf{h}_{r2}')'$, and $\mathbf{G}_r = (\mathbf{g}_{r1}', \mathbf{g}_{r2}')'$ is the component notations of the $r$th order term in the expansion of $\sqrt{n}(\hat{\boldsymbol{\alpha}}_n' - \boldsymbol{\alpha}_0', \hat{\boldsymbol{\beta}}_n' - \boldsymbol{\beta}_0')'$ and $\sqrt{n}(\check{\boldsymbol{\alpha}}_n' - \boldsymbol{\alpha}_0', \check{\boldsymbol{\beta}}_n' - \boldsymbol{\beta}_0')'$. In Theorem 2.2, $\mathbf{H}_r^\circ$ is a function of $\mathbf{H}_0^\circ, \ldots, \mathbf{H}_{r-1}^\circ$. In terms of the components, $\mathbf{h}_{r1}$ is a function of $\mathbf{h}_{01}, \ldots, \mathbf{h}_{r-1,1}$ and $\mathbf{h}_{02}, \ldots, \mathbf{h}_{r-1,2}$, with some evaluations at $\boldsymbol{\theta}_0$, similarly for $\mathbf{h}_{r2}$. We denote these functions as, for $r = 1, \ldots, k-1$,

$$\mathbf{h}_{r1} = \mathcal{H}_{r1}(\mathbf{h}_{01}, \mathbf{h}_{02}, \ldots, \mathbf{h}_{r-1,1}, \mathbf{h}_{r-1,2} | \boldsymbol{\theta}_0),$$
$$\mathbf{h}_{r2} = \mathcal{H}_{r2}(\mathbf{h}_{01}, \mathbf{h}_{02}, \ldots, \mathbf{h}_{r-1,1}, \mathbf{h}_{r-1,2} | \boldsymbol{\theta}_0).$$

By Theorem 2.3, the components $\mathbf{g}_{r1}$ and $\mathbf{g}_{r2}$ are also functions of $\mathbf{g}_{01}, \ldots, \mathbf{g}_{r-1,1}$ and $\mathbf{g}_{02}, \ldots, \mathbf{g}_{r-1,2}$. We denote them as, for $r = 1, \ldots, k-1$,

$$\mathbf{g}_{r1} = \mathcal{G}_{r1}(\mathbf{g}_{01}, \mathbf{g}_{02}, \ldots, \mathbf{g}_{r-1,1}, \mathbf{g}_{r-1,2} | \boldsymbol{\theta}_0),$$
$$\mathbf{g}_{r2} = \mathcal{G}_{r2}(\mathbf{g}_{01}, \mathbf{g}_{02}, \ldots, \mathbf{g}_{r-1,1}, \mathbf{g}_{r-1,2} | \boldsymbol{\theta}_0).$$

The expansion in Theorem 2.3 can be obtained by another way. Fix $\check{\boldsymbol{\beta}}_n$ in the posterior and expand $\sqrt{n}(\check{\boldsymbol{\alpha}}_n - \boldsymbol{\alpha}_0)$ to get

$$\sqrt{n}(\check{\boldsymbol{\alpha}}_n - \boldsymbol{\alpha}_0)|_{\check{\boldsymbol{\beta}}_n} = \sum_{r=0}^{k-1} n^{-r/2} \bar{\mathbf{g}}_{r1}(\check{\boldsymbol{\beta}}_n) + O_p(n^{-k/2}),$$

where $\bar{\mathbf{g}}_{r1}(\cdot)$ is the $d_1$-dimensional version of $\mathbf{G}_r(\cdot)$ and with $\boldsymbol{\theta}_0$ replaced by $(\boldsymbol{\alpha}_0, \check{\boldsymbol{\beta}}_n)$. Now expand $\sqrt{n}(\check{\boldsymbol{\beta}}_n - \boldsymbol{\beta}_0)$ in the $\bar{\mathbf{g}}_{r1}(\cdot)$'s above; this gives the same expansion as that in Theorem 2.3, so we must have

$$\sum_{r=0}^{k-1} n^{-r/2} \bar{\mathbf{g}}_{r1}(\check{\boldsymbol{\beta}}_n) = \sum_{r=0}^{k-1} n^{-r/2} \mathcal{G}_{r1}(\mathbf{g}_{01}, \mathbf{g}_{02}, \ldots, \mathbf{g}_{r1}, \mathbf{g}_{r2}) + O_p(n^{-k/2}).$$

On the other hand, for the hybrid estimator $(\check{\boldsymbol{\alpha}}_n, \hat{\boldsymbol{\beta}}_n)$ we can expand the two components simultaneously or componentwise, and the two expansions are the same. We fix $\hat{\boldsymbol{\beta}}_n$, first expand $\sqrt{n}(\check{\boldsymbol{\alpha}}_n - \boldsymbol{\alpha}_0)$ and we have

$$\sqrt{n}(\check{\boldsymbol{\alpha}}_n - \boldsymbol{\alpha}_0)|_{\hat{\boldsymbol{\beta}}_n} = \sum_{r=0}^{k-1} n^{-r/2} \bar{\mathbf{g}}_{r1}(\hat{\boldsymbol{\beta}}_n) + O_p(n^{-k/2}).$$

Then, expand $\sqrt{n}(\hat{\boldsymbol{\beta}}_n - \boldsymbol{\beta}_0)$ in the $\bar{\mathbf{g}}_{r1}(\cdot)$'s. Comparing the procedures for $\sqrt{n}(\check{\boldsymbol{\alpha}}_n - \boldsymbol{\alpha}_0)|_{\check{\boldsymbol{\beta}}_n}$ and $\sqrt{n}(\check{\boldsymbol{\alpha}}_n - \boldsymbol{\alpha}_0)|_{\hat{\boldsymbol{\beta}}_n}$, and note the statuses of $\check{\boldsymbol{\beta}}_n$ and $\hat{\boldsymbol{\beta}}_n$ in the $\bar{\mathbf{g}}_{r1}(\cdot)$'s are the same, except that the former expands in terms of $\mathbf{g}_{01}, \mathbf{g}_{02}, \ldots, \mathbf{g}_{r1}, \mathbf{g}_{r2}$ and the latter in terms of $\mathbf{g}_{01}, \mathbf{h}_{02}, \ldots, \mathbf{g}_{r1}, \mathbf{h}_{r2}$. So we get

$$\sum_{r=0}^{k-1} n^{-r/2} \bar{\mathbf{g}}_{r1}(\hat{\boldsymbol{\beta}}_n) = \sum_{r=0}^{k-1} n^{-r/2} \mathcal{G}_{r1}(\mathbf{g}_{01}, \mathbf{h}_{02}, \ldots, \mathbf{g}_{r1}, \mathbf{h}_{r2}) + O_p(n^{-k/2})$$

$$= \sum_{r=0}^{k-1} n^{-r/2} [\mathcal{H}_{r1}(\mathbf{g}_{01}, \mathbf{h}_{02}, \ldots, \mathbf{g}_{r1}, \mathbf{h}_{r2}) + \mathbf{t}_r] + O_p(n^{-k/2}).$$

The last step above is by the recursive relationship in Theorem 2.3.

In the same way, for the MLE $(\hat{\boldsymbol{\alpha}}_n, \hat{\boldsymbol{\beta}}_n)$, we have

$$\sqrt{n}(\hat{\boldsymbol{\beta}}_n - \boldsymbol{\beta}_0)|_{\hat{\boldsymbol{\alpha}}_n} = \sum_{r=0}^{k-1} n^{-r/2} \bar{\mathbf{h}}_{r2}(\hat{\boldsymbol{\alpha}}_n) + O_p(n^{-k/2}),$$

where $\bar{\mathbf{h}}_{r2}(\cdot)$ is the $d_2$-dimensional version of $\mathbf{H}_r$ and with $\boldsymbol{\theta}_0$ replaced by $(\hat{\boldsymbol{\alpha}}_n, \boldsymbol{\beta}_0)$. As before, we have

$$\sum_{r=0}^{k-1} n^{-r/2} \bar{\mathbf{h}}_{r2}(\hat{\boldsymbol{\alpha}}_n) = \sum_{r=0}^{k-1} n^{-r/2} \mathcal{H}_{r2}(\mathbf{h}_{01}, \mathbf{h}_{02}, \ldots, \mathbf{h}_{r1}, \mathbf{h}_{r2}) + O_p(n^{-k/2}).$$

Similarly, for the hybrid estimator $(\check{\boldsymbol{\alpha}}_n, \hat{\boldsymbol{\beta}}_n)$, we have

$$\sqrt{n}(\hat{\boldsymbol{\beta}}_n - \boldsymbol{\beta}_0)|_{\check{\boldsymbol{\alpha}}_n} = \sum_{r=0}^{k-1} n^{-r/2} \bar{\mathbf{h}}_{r2}(\check{\boldsymbol{\alpha}}_n) + O_p(n^{-k/2})$$

$$= \sum_{r=0}^{k-1} n^{-r/2} \mathcal{H}_{r2}(\mathbf{g}_{01}, \mathbf{h}_{02}, \ldots, \mathbf{g}_{r1}, \mathbf{h}_{r2}) + O_p(n^{-k/2}).$$

These give

$$\begin{pmatrix} \mathbf{g}_r \\ \mathbf{h}_r \end{pmatrix} = \begin{pmatrix} \mathcal{H}_{r1}(\mathbf{g}_0, \mathbf{h}_0, \ldots, \mathbf{g}_{r-1}, \mathbf{h}_{r-1}|\boldsymbol{\theta}_0) \\ \mathcal{H}_{r2}(\mathbf{g}_0, \mathbf{h}_0, \ldots, \mathbf{g}_{r-1}, \mathbf{h}_{r-1}|\boldsymbol{\theta}_0) \end{pmatrix} + \begin{pmatrix} \mathbf{t}_r \\ \mathbf{0} \end{pmatrix},$$

which is the recursive formula in the theorem.

Below, we go into some details of the above sketch. We need the following notation. Let $\mathbf{I} = (\mathbf{I}_{ij})_{1 \leq i,j \leq 2}$ be the partition of the Fisher information matrix, where $\mathbf{I}_{11}$ corresponding to the block for $\boldsymbol{\alpha}$ and $\mathbf{I}_{22}$ for $\boldsymbol{\beta}$. For nonnegative integer $d_1$-vector $\mathbf{i}$ and nonnegative integer $d_2$-vector $\mathbf{j}$, define

$$\boldsymbol{\rho}^1(\boldsymbol{\alpha}) := (\rho^1_1(\boldsymbol{\alpha}), \ldots, \rho^1_{d_1}(\boldsymbol{\alpha}))' = \left( \frac{\partial}{\partial \alpha_1} \log \pi(\boldsymbol{\alpha}), \ldots, \frac{\partial}{\partial \alpha_{d_1}} \log \pi(\boldsymbol{\alpha}) \right)',$$

$$\boldsymbol{\rho}^1_{\mathbf{i}}(\boldsymbol{\alpha}) = \left( \frac{\partial^{|\mathbf{i}|}}{\partial \boldsymbol{\alpha}^{\mathbf{i}}} \rho^1_1(\boldsymbol{\alpha}), \ldots, \frac{\partial^{|\mathbf{i}|}}{\partial \boldsymbol{\alpha}^{\mathbf{i}}} \rho^1_{d_1}(\boldsymbol{\alpha}) \right)',$$

$$\mathbf{L}^1(x|\boldsymbol{\alpha}, \boldsymbol{\beta}) = \left( \frac{\partial}{\partial \alpha_1} \log f(x|\boldsymbol{\alpha}, \boldsymbol{\beta}), \ldots, \frac{\partial}{\partial \alpha_{d_1}} \log f(x|\boldsymbol{\alpha}, \boldsymbol{\beta}) \right)'$$

$$:= (L^1_1(x|\boldsymbol{\alpha}, \boldsymbol{\beta}), \ldots, L^1_{d_1}(x|\boldsymbol{\alpha}, \boldsymbol{\beta}))',$$

$$\mathbf{L}^1_{(\mathbf{i};\mathbf{j})}(x|\boldsymbol{\alpha}, \boldsymbol{\beta}) = \left( \frac{\partial^{|\mathbf{i}+\mathbf{j}|}}{\partial \boldsymbol{\alpha}^{\mathbf{i}} \partial \boldsymbol{\beta}^{\mathbf{j}}} L^1_1(x|\boldsymbol{\alpha}, \boldsymbol{\beta}), \ldots, \frac{\partial^{|\mathbf{i}+\mathbf{j}|}}{\partial \boldsymbol{\alpha}^{\mathbf{i}} \partial \boldsymbol{\beta}^{\mathbf{j}}} L^1_{d_1}(x|\boldsymbol{\alpha}, \boldsymbol{\beta}) \right)',$$

$$\mathbf{L}^2(x|\boldsymbol{\alpha}, \boldsymbol{\beta}) = \left( \frac{\partial}{\partial \beta_1} \log f(x|\boldsymbol{\alpha}, \boldsymbol{\beta}), \ldots, \frac{\partial}{\partial \beta_{d_2}} \log f(x|\boldsymbol{\alpha}, \boldsymbol{\beta}) \right)'$$

$$:= (L^2_1(x|\boldsymbol{\alpha}, \boldsymbol{\beta}), \ldots, L^2_{d_2}(x|\boldsymbol{\alpha}, \boldsymbol{\beta}))',$$

$$\mathbf{L}^2_{(\mathbf{i};\mathbf{j})}(x|\boldsymbol{\alpha}, \boldsymbol{\beta}) = \left( \frac{\partial^{|\mathbf{i}+\mathbf{j}|}}{\partial \boldsymbol{\alpha}^{\mathbf{i}} \partial \boldsymbol{\beta}^{\mathbf{j}}} L^2_1(x|\boldsymbol{\alpha}, \boldsymbol{\beta}), \ldots, \frac{\partial^{|\mathbf{i}+\mathbf{j}|}}{\partial \langle \boldsymbol{\alpha}^{\mathbf{i}} \rangle \partial \langle \boldsymbol{\beta}^{\mathbf{j}} \rangle} L^2_{d_2}(x|\boldsymbol{\alpha}, \boldsymbol{\beta}) \right)',$$

$$\mathbf{S}^1_{(\mathbf{i};\mathbf{j})}(\boldsymbol{\alpha}, \boldsymbol{\beta}) = \frac{1}{\sqrt{n}} \sum_{j=1}^{n} \mathbf{L}^1_{(\mathbf{i};\mathbf{j})}(x_j|\boldsymbol{\alpha}, \boldsymbol{\beta}), \qquad \mathbf{S}^1_{(\mathbf{i};\mathbf{j})}(\boldsymbol{\beta}) = \frac{1}{\sqrt{n}} \sum_{j=1}^{n} \mathbf{L}^1_{(\mathbf{i};\mathbf{j})}(x_j|\boldsymbol{\alpha}_0, \boldsymbol{\beta}),$$

$$\boldsymbol{\Delta}^1_{(\mathbf{i};\mathbf{j})}(\boldsymbol{\alpha}, \boldsymbol{\beta}) = \frac{1}{\sqrt{n}} \sum_{j=1}^{n} (\mathbf{L}^1_{(\mathbf{i};\mathbf{j})}(X_j|\boldsymbol{\alpha}, \boldsymbol{\beta}) - E_{(\boldsymbol{\alpha}, \boldsymbol{\beta})} \mathbf{L}^1_{(\mathbf{i};\mathbf{j})}(X|\boldsymbol{\alpha}, \boldsymbol{\beta})),$$

$$\mathbf{\Delta}^1_{(\mathbf{i};\mathbf{j})}(\boldsymbol{\beta}) = \frac{1}{\sqrt{n}} \sum_{j=1}^n (\mathbf{L}^1_{(\mathbf{i};\mathbf{j})}(X_j|\boldsymbol{\alpha}_0,\boldsymbol{\beta}) - E_{(\boldsymbol{\alpha}_0,\boldsymbol{\beta})} \mathbf{L}^1_{(\mathbf{i};\mathbf{j})}(X|\boldsymbol{\alpha}_0,\boldsymbol{\beta})),$$

$$\mathbf{E}^1_{(\mathbf{i};\mathbf{j})}(\boldsymbol{\beta}) = E_{(\boldsymbol{\alpha}_0,\boldsymbol{\beta})} \mathbf{L}^1_{(\mathbf{i};\mathbf{j})}(X|\boldsymbol{\alpha}_0,\boldsymbol{\beta}), \qquad \mathbf{E}^1_{(\mathbf{i};\mathbf{j})} = E_{(\boldsymbol{\alpha}_0,\boldsymbol{\beta}_0)} \mathbf{L}^1_{(\mathbf{i};\mathbf{j})}(X|\boldsymbol{\alpha}_0,\boldsymbol{\beta}_0).$$

$$\mathbf{S}^1_{(\mathbf{i};\mathbf{j})} = \mathbf{S}^1_{(\mathbf{i};\mathbf{j})}(\boldsymbol{\alpha}_0,\boldsymbol{\beta}_0), \qquad \mathbf{\Delta}^1_{(\mathbf{i};\mathbf{j})} = \mathbf{\Delta}^1_{(\mathbf{i};\mathbf{j})}(\boldsymbol{\alpha}_0,\boldsymbol{\beta}_0),$$

$$\mathbf{S}^2_{(\mathbf{i};\mathbf{j})}(\boldsymbol{\alpha},\boldsymbol{\beta}) = \frac{1}{\sqrt{n}} \sum_{j=1}^n \mathbf{L}^2_{(\mathbf{i};\mathbf{j})}(X_j|\boldsymbol{\alpha},\boldsymbol{\beta}),$$

$$\mathbf{S}^2_{(\mathbf{i};\mathbf{j})}(\boldsymbol{\alpha}) = \frac{1}{\sqrt{n}} \sum_{j=1}^n \mathbf{L}^2_{(\mathbf{i};\mathbf{j})}(X_j|\boldsymbol{\alpha},\boldsymbol{\beta}_0),$$

and define $\mathbf{\Delta}^2_{(\mathbf{i};\mathbf{j})}(\boldsymbol{\alpha},\boldsymbol{\beta})$, $\mathbf{\Delta}^2_{(\mathbf{i};\mathbf{j})}(\boldsymbol{\alpha})$, $\mathbf{E}^2_{(\mathbf{i};\mathbf{j})}(\boldsymbol{\alpha})$ and $\mathbf{S}^2_{(\mathbf{i};\mathbf{j})}$ accordingly.

We first give the expression for $\mathbf{h}_0$; the expression for $\mathbf{g}_0$ will be outlined later. Fix $\check{\boldsymbol{\alpha}}_n$, note here $\boldsymbol{\rho}^1_0$ does not depend on $\boldsymbol{\beta}$, $\mathbf{E}^2_{(0;0)} = \mathbf{0}$ and $\mathbf{E}^2_{(0;0)}(\check{\boldsymbol{\alpha}}_n) \neq \mathbf{0}$. Let $\hat{\boldsymbol{\beta}}'_n = \sqrt{n}(\hat{\boldsymbol{\beta}}_n - \boldsymbol{\beta}_0)$, and $\check{\boldsymbol{\alpha}}'_n = \sqrt{n}(\check{\boldsymbol{\alpha}}_n - \boldsymbol{\alpha}_0)$ (not to be confused with the transpose). As in the proof of Theorem 2.2, we have

$$\mathbf{0} = \mathbf{S}^2_{(0;0)}(\check{\boldsymbol{\alpha}}_n, \hat{\boldsymbol{\beta}}_n) = \mathbf{S}^2_{(0,0)}(\check{\boldsymbol{\alpha}}_n, \boldsymbol{\beta}_0 + n^{-1/2}\hat{\boldsymbol{\beta}}'_n)$$

$$\overset{k}{\sim} \sum_{r=0}^{k-1} n^{-r/2} \sum_{|\mathbf{j}|=r} \frac{\langle (\boldsymbol{\beta}'_n)^{\mathbf{j}} \rangle}{\mathbf{j}!} \mathbf{S}^2_{(0;\mathbf{j})}(\check{\boldsymbol{\alpha}}_n)$$

$$\overset{k}{\sim} n^{1/2} \mathbf{E}^2_{(0;0)}(\check{\boldsymbol{\alpha}}_n)$$

$$+ \sum_{r=0}^{k-1} n^{-r/2} \left( \sum_{|\mathbf{j}|=r} \frac{\langle (\hat{\boldsymbol{\beta}}'_n)^{\mathbf{j}} \rangle}{\mathbf{j}!} \mathbf{\Delta}^2_{(0;\mathbf{j})}(\check{\boldsymbol{\alpha}}_n) + \sum_{|\mathbf{j}|=r+1} \frac{\langle (\hat{\boldsymbol{\beta}}'_n)^{\mathbf{j}} \rangle}{\mathbf{j}!} \mathbf{E}^2_{(0;\mathbf{j})}(\check{\boldsymbol{\alpha}}_n) \right)$$

$$\overset{k}{\sim} \sum_{r=0}^{k-1} n^{-r/2} \left( \sum_{|\mathbf{i}|=r+1} \frac{\langle (\check{\boldsymbol{\alpha}}'_n)^{\mathbf{i}} \rangle}{\mathbf{i}!} \mathbf{E}^2_{(\mathbf{i};0)} \right.$$

$$+ \sum_{s+t=r} \sum_{|\mathbf{j}|=s} \frac{\langle (\hat{\boldsymbol{\beta}}'_n)^{\mathbf{j}} \rangle}{\mathbf{j}!} \sum_{|\mathbf{i}|=t} \frac{\langle (\check{\boldsymbol{\alpha}}'_n)^{\mathbf{i}} \rangle}{\mathbf{i}!} \mathbf{\Delta}^2_{(\mathbf{i};\mathbf{j})}$$

$$\left. + \sum_{s+t=r} \sum_{|\mathbf{j}|=s+1} \frac{\langle (\hat{\boldsymbol{\beta}}'_n)^{\mathbf{j}} \rangle}{\mathbf{j}!} \sum_{|\mathbf{i}|=t} \frac{\langle (\check{\boldsymbol{\alpha}}'_n)^{\mathbf{i}} \rangle}{\mathbf{i}!} \mathbf{E}^2_{(\mathbf{i};\mathbf{j})} \right)$$

$$= \mathbf{\Delta}^2_{(0;0)} - \mathbf{I}_{21} \check{\boldsymbol{\alpha}}'_n - \mathbf{I}_{22} \hat{\boldsymbol{\beta}}'_n$$

(A.10) $\qquad + \sum_{r=1}^{k-1} n^{-r/2} \left( \sum_{|\mathbf{i}|=r+1} \frac{\langle (\check{\boldsymbol{\alpha}}'_n)^{\mathbf{i}} \rangle}{\mathbf{i}!} \mathbf{E}^2_{(\mathbf{i};0)} \right.$

$$+ \sum_{s+t=r} \sum_{|\mathbf{j}|=s} \frac{\langle(\hat{\boldsymbol{\beta}}'_n)^{\mathbf{j}}\rangle}{\mathbf{j}!} \sum_{|\mathbf{i}|=t} \frac{\langle(\check{\boldsymbol{\alpha}}'_n)^{\mathbf{i}}\rangle}{\mathbf{i}!} \boldsymbol{\Delta}^2_{(\mathbf{i};\mathbf{j})}$$

$$+ \sum_{s+t=r} \sum_{|\mathbf{j}|=s+1} \frac{\langle(\hat{\boldsymbol{\beta}}'_n)^{\mathbf{j}}\rangle}{\mathbf{j}!} \sum_{|\mathbf{i}|=t} \frac{\langle(\check{\boldsymbol{\alpha}}'_n)^{\mathbf{i}}\rangle}{\mathbf{i}!} \mathbf{E}^2_{(\mathbf{i};\mathbf{j})} \Big),$$

or $\mathbf{I}_{21}\check{\boldsymbol{\alpha}}'_n + \mathbf{I}_{22}\hat{\boldsymbol{\beta}}'_n = \boldsymbol{\Delta}^2_{(\mathbf{0},\mathbf{0})} + O_p(n^{-1/2})$. Similarly, fix $\hat{\boldsymbol{\beta}}_n$ and expand $\check{\boldsymbol{\alpha}}_n$, we will have $\mathbf{I}_{11}\check{\boldsymbol{\alpha}}_n + \mathbf{I}_{12}\hat{\boldsymbol{\beta}}_n = \boldsymbol{\Delta}^1_{(\mathbf{0},\mathbf{0})} + O_p(n^{-1/2})$. So we get $\sqrt{n}\mathbf{I}\binom{\check{\boldsymbol{\alpha}}_n - \boldsymbol{\alpha}_0}{\hat{\boldsymbol{\beta}}_n - \boldsymbol{\beta}_0} = \binom{\boldsymbol{\Delta}^1_{(\mathbf{0},\mathbf{0})}}{\boldsymbol{\Delta}^2_{(\mathbf{0},\mathbf{0})}} + O_p(n^{-1/2})$, and note $(\boldsymbol{\Delta}^1_{(\mathbf{0};\mathbf{0})}, \boldsymbol{\Delta}^2_{(\mathbf{0};\mathbf{0})})' = \boldsymbol{\Delta}_{\mathbf{0}}$; this gives the expression for $\mathbf{g}_0$ and $\mathbf{h}_0$.

To prove the expansions for the $\mathbf{g}_r$'s and $\mathbf{h}_r$'s, we use induction. We only need to prove those for the $\mathbf{g}'_r$'s and $\mathbf{h}'_r$'s, where $\binom{\mathbf{g}'_r}{\mathbf{h}'_r} = \mathbf{I}\binom{\mathbf{g}_r}{\mathbf{h}_r}$. Now we prove

$$\mathbf{h}'_r = \sum_{s+t=r-1} \sum_{|\mathbf{i}|=s+1} \mathbf{E}^2_{(\mathbf{i};\mathbf{0})} \sum_{|\mathbf{l}|=t} \sum_{(0,t,\mathbf{l},\mathbf{i})} \prod_{v=0}^{t} \langle \mathbf{g}^{\mathbf{i}_v}_v \rangle / \mathbf{i}_v!$$

$$+ \sum_{a+b+c=r-1} \sum_{s+t=c} \sum_{|\mathbf{j}|=s} \sum_{|\mathbf{i}|=t} \boldsymbol{\Delta}^2_{(\mathbf{i};\mathbf{j})} \sum_{|\mathbf{l}|=a} \sum_{(0,a,\mathbf{l},\mathbf{i})} \prod_{v=0}^{a} \langle \mathbf{g}^{\mathbf{i}_v}_v \rangle / \mathbf{i}_v!$$

$$\times \sum_{|\mathbf{l}|=b} \sum_{(0,b,\mathbf{l},\mathbf{j})} \prod_{v=0}^{b} \langle \mathbf{h}^{\mathbf{i}_v}_v \rangle / \mathbf{i}_v!$$

$$+ \sum_{a+b+c=r-1} \sum_{s+t=c} \sum_{|\mathbf{j}|=s+1} \sum_{|\mathbf{i}|=t} \mathbf{E}^2_{(\mathbf{i};\mathbf{j})} \sum_{|\mathbf{l}|=a} \sum_{(0,a,\mathbf{l},\mathbf{i})} \prod_{v=0}^{a} \langle \mathbf{g}^{\mathbf{i}_v}_v \rangle / \mathbf{i}_v!$$

$$\times \sum_{|\mathbf{l}|=b} \sum_{(0,b,\mathbf{l},\mathbf{j})} \prod_{v=0}^{b} \langle \mathbf{h}^{\mathbf{i}_v}_v \rangle / \mathbf{i}_v!,$$

$\mathbf{g}'_r = \mathbf{h}'_r + \mathbf{q}'_r$ ($1 \le r \le k-1$), and the $\mathbf{q}'_r$'s are outlined later.

Since in the following $t = r$, we get

$$\sum_{s+t=r} \sum_{|\mathbf{i}|=s+1} \mathbf{E}^2_{(\mathbf{i};\mathbf{0})} \sum_{|\mathbf{l}|=t} \sum_{(0,t,\mathbf{l},\mathbf{i})} \prod_{v=0}^{t} \langle \mathbf{g}^{\mathbf{i}_v}_v \rangle / \mathbf{i}_v! = \sum_{|\mathbf{i}|=1} \mathbf{E}^2_{(\mathbf{i};\mathbf{0})} \langle \mathbf{g}^{\mathbf{i}}_r \rangle = -\mathbf{I}_{21}\mathbf{g}_r;$$

similarly, in the following, we take $b = r$ (note $a = r$ will result in summations over empty sets) and we get

$$\sum_{a+b+c=r} \sum_{s+t=c} \sum_{|\mathbf{j}|=s+1} \sum_{|\mathbf{i}|=t} \mathbf{E}^2_{(\mathbf{i};\mathbf{j})} \sum_{|\mathbf{l}|=a} \sum_{(0,a,\mathbf{l},\mathbf{i})} \prod_{v=0}^{a} \langle \mathbf{g}^{\mathbf{i}_v}_v \rangle / \mathbf{i}_v! \sum_{|\mathbf{l}|=b} \sum_{(0,b,\mathbf{l},\mathbf{j})} \prod_{v=0}^{b} \langle \mathbf{h}^{\mathbf{i}_v}_v \rangle / \mathbf{i}_v!$$

$$= -\mathbf{I}_{22}\mathbf{h}_r,$$

and the previous expression for $\tilde{\mathbf{h}}_r = \mathbf{I}_{21}\mathbf{g}_r + \mathbf{I}_{22}\mathbf{h}_r$ $(1 \leq r \leq k-1)$ is rewritten as

$$\sum_{s+t=r}\sum_{|\mathbf{i}|=s+1}\mathbf{E}^2_{(\mathbf{i};\mathbf{0})}\sum_{|\mathbf{l}|=t}\sum_{(0,t,\mathbf{l},\mathbf{i})}\prod_{v=0}^{t}\langle \mathbf{g}_v^{\mathbf{i}_v}\rangle/\mathbf{i}_v!$$

$$+ \sum_{a+b+c=r}\sum_{s+t=c}\sum_{|\mathbf{j}|=s}\sum_{|\mathbf{i}|=t}\boldsymbol{\Delta}^2_{(\mathbf{i};\mathbf{j})}\sum_{|\mathbf{l}|=a}\sum_{(0,a,\mathbf{l},\mathbf{i})}\prod_{v=0}^{a}\langle \mathbf{g}_v^{\mathbf{i}_v}\rangle/\mathbf{i}_v!$$

$$(\text{A.11}) \qquad\qquad \times \sum_{|\mathbf{l}|=b}\sum_{(0,b,\mathbf{l},\mathbf{j})}\prod_{v=0}^{b}\langle \mathbf{h}_v^{\mathbf{i}_v}\rangle/\mathbf{i}_v!$$

$$+ \sum_{a+b+c=r}\sum_{s+t=c}\sum_{|\mathbf{j}|=s+1}\sum_{|\mathbf{i}|=t}\mathbf{E}^2_{(\mathbf{i};\mathbf{j})}\sum_{|\mathbf{l}|=a}\sum_{(0,a,\mathbf{l},\mathbf{i})}\prod_{v=0}^{a}\langle \mathbf{g}_v^{\mathbf{i}_v}\rangle/\mathbf{i}_v!$$

$$\times \sum_{|\mathbf{l}|=b}\sum_{(0,b,\mathbf{l},\mathbf{j})}\prod_{v=0}^{b}\langle \mathbf{h}_v^{\mathbf{i}_v}\rangle/\mathbf{i}_v! = \mathbf{0}.$$

Let $\check{\boldsymbol{\alpha}}_n'' = \sum_{r=0}^{k-1} n^{-r/2}\mathbf{g}_r$ and $\hat{\boldsymbol{\beta}}_n'' = \sum_{r=0}^{k-1} n^{-r/2}\mathbf{h}_r$, then

$$\langle (\check{\boldsymbol{\alpha}}_n'')^{\mathbf{i}}\rangle \overset{k}{\sim} \mathbf{i}!\sum_{r=0}^{k-1} n^{-r/2}\sum_{|\mathbf{l}|=r}\sum_{(0,r,\mathbf{l},\mathbf{i})}\prod_{v=0}^{r}\langle \mathbf{g}_v^{\mathbf{i}_v}\rangle/\mathbf{i}_v!,$$

$$\langle (\hat{\boldsymbol{\beta}}_n'')^{\mathbf{j}}\rangle \overset{k}{\sim} \mathbf{j}!\sum_{r=0}^{k-1} n^{-r/2}\sum_{|\mathbf{l}|=r}\sum_{(0,r,\mathbf{l},\mathbf{j})}\prod_{v=0}^{r}\langle \mathbf{h}_v^{\mathbf{i}_v}\rangle/\mathbf{i}_v!.$$

By (A.10) and (A.11), we have

$$\boldsymbol{\Delta}^2_{(\mathbf{0};\mathbf{0})} - \mathbf{I}_{21}\check{\boldsymbol{\alpha}}_n'' - \mathbf{I}_{22}\hat{\boldsymbol{\beta}}_n''$$

$$+ \sum_{r=1}^{k-1} n^{-r/2}\Bigg(\sum_{|\mathbf{i}|=r+1}\frac{\langle(\check{\boldsymbol{\alpha}}_n'')^{\mathbf{i}}\rangle}{\mathbf{i}!}\mathbf{E}^2_{(\mathbf{i};\mathbf{0})}$$

$$+ \sum_{s+t=r}\sum_{|\mathbf{j}|=s}\frac{\langle(\hat{\boldsymbol{\beta}}_n')^{\mathbf{j}}\rangle}{\mathbf{j}!}\sum_{|\mathbf{i}|=t}\frac{\langle(\check{\boldsymbol{\alpha}}_n'')^{\mathbf{i}}\rangle}{\mathbf{i}!}\boldsymbol{\Delta}^2_{(\mathbf{i};\mathbf{j})}$$

$$+ \sum_{s+t=r}\sum_{|\mathbf{j}|=s+1}\frac{\langle(\hat{\boldsymbol{\beta}}_n'')^{\mathbf{j}}\rangle}{\mathbf{j}!}\sum_{|\mathbf{i}|=t}\frac{\langle(\check{\boldsymbol{\alpha}}_n'')^{\mathbf{i}}\rangle}{\mathbf{i}!}\mathbf{E}^2_{(\mathbf{i};\mathbf{j})}\Bigg)$$

$$(\text{A.12}) = \sum_{r=0}^{k-1} n^{-r/2}\Bigg(\sum_{|\mathbf{i}|=r+1}\frac{\langle(\check{\boldsymbol{\alpha}}_n'')^{\mathbf{i}}\rangle}{\mathbf{i}!}\mathbf{E}^2_{(\mathbf{i};\mathbf{0})} + \sum_{s+t=r}\sum_{|\mathbf{j}|=s}\frac{\langle(\hat{\boldsymbol{\beta}}_n'')^{\mathbf{j}}\rangle}{\mathbf{j}!}\sum_{|\mathbf{i}|=t}\frac{\langle(\check{\boldsymbol{\alpha}}_n'')^{\mathbf{i}}\rangle}{\mathbf{i}!}\boldsymbol{\Delta}^2_{(\mathbf{i};\mathbf{j})}$$

$$+ \sum_{s+t=r} \sum_{|\mathbf{j}|=s+1} \frac{\langle (\hat{\boldsymbol{\beta}}_n'')^{\mathbf{j}} \rangle}{\mathbf{j}!} \sum_{|\mathbf{i}|=t} \frac{\langle (\check{\boldsymbol{\alpha}}_n'')^{\mathbf{i}} \rangle}{\mathbf{i}!} \mathbf{E}^2_{(\mathbf{i};\mathbf{j})} \Bigg)$$

$$\overset{k}{\sim} \sum_{r=0}^{k-1} n^{-r/2} \Bigg( \sum_{s+t=r} \sum_{|\mathbf{i}|=s+1} \mathbf{E}^2_{(\mathbf{i};\mathbf{0})} \sum_{|\mathbf{l}|=t} \sum_{(0,t,\mathbf{l},\mathbf{i})} \prod_{v=0}^{t} \langle \mathbf{g}_v^{\mathbf{i}_v} \rangle / \mathbf{i}_v!$$

$$+ \sum_{a+b+c=r} \sum_{s+t=c} \sum_{|\mathbf{j}|=s} \sum_{|\mathbf{i}|=t} \boldsymbol{\Delta}^2_{(\mathbf{i};\mathbf{j})} \sum_{|\mathbf{l}|=a} \sum_{(0,a,\mathbf{l},\mathbf{i})} \prod_{v=0}^{a} \langle \mathbf{g}_v^{\mathbf{i}_v} \rangle / \mathbf{i}_v!$$

$$\times \sum_{|\mathbf{l}|=b} \sum_{(0,b,\mathbf{l},\mathbf{j})} \prod_{v=0}^{b} \langle \mathbf{h}_v^{\mathbf{i}_v} \rangle / \mathbf{i}_v!$$

$$+ \sum_{a+b+c=r} \sum_{s+t=c} \sum_{|\mathbf{j}|=s+1} \sum_{|\mathbf{i}|=t} \mathbf{E}^2_{(\mathbf{i};\mathbf{j})} \sum_{|\mathbf{l}|=a} \sum_{(0,a,\mathbf{l},\mathbf{i})} \prod_{v=0}^{a} \langle \mathbf{g}_v^{\mathbf{i}_v} \rangle / \mathbf{i}_v!$$

$$\times \sum_{|\mathbf{l}|=b} \sum_{(0,b,\mathbf{l},\mathbf{j})} \prod_{v=0}^{b} \langle \mathbf{h}_v^{\mathbf{i}_v} \rangle / \mathbf{i}_v! \Bigg)$$

$$= \mathbf{0}.$$

Thus, (A.10) minus (A.12) will give the expression for $\mathbf{h}'_r$.

Now, we outline the expressions for the $\mathbf{g}'_r$'s. We need to modify the result in Lemma 1. For fixed $\hat{\boldsymbol{\beta}}_n$, define

$$Z_n(\boldsymbol{\alpha}, \hat{\boldsymbol{\beta}}_n) = \left( \prod_{i=1}^{n} \frac{f(x_i | \boldsymbol{\alpha}_0 + \boldsymbol{\alpha} n^{-1/2}, \hat{\boldsymbol{\beta}}_n)}{f(x_i | \boldsymbol{\alpha}_0, \hat{\boldsymbol{\beta}}_n)} \right) \frac{\pi(\boldsymbol{\alpha}_0 + \boldsymbol{\alpha} n^{-1/2})}{\pi(\boldsymbol{\alpha}_0)}.$$

Define $\mathcal{E}^1_{(\mathbf{i};\mathbf{j})}(\cdot)$, $\delta^1_{(\mathbf{i};\mathbf{j})}(\cdot)$ and $\varrho^1_{\mathbf{j}}$ accordingly. We have

$$\mathcal{E}^1_{(\mathbf{i};\mathbf{0})}(\hat{\boldsymbol{\beta}}_n) \overset{k}{\sim} \sum_{r=0}^{k-1} n^{-r/2} \sum_{|\mathbf{j}|=r} \mathcal{E}^1_{(\mathbf{i};\mathbf{j})} \langle (\hat{\boldsymbol{\beta}}_n)^{\mathbf{j}} \rangle / \mathbf{j}!$$

$$= \sum_{r=0}^{k-1} n^{-r/2} \sum_{|\mathbf{j}|=r} \sum_{s+t=r} \sum_{|\mathbf{j}|=s} \mathcal{E}^1_{(\mathbf{i};\mathbf{j})} \sum_{|\mathbf{l}|=t} \sum_{(0,\mathbf{l},\mathbf{j})} \prod_{v=0}^{t} \langle (\mathbf{h}_v)^{\mathbf{i}_v} \rangle / \mathbf{i}_v!$$

and

$$\delta^1_{(\mathbf{i};\mathbf{0})}(\hat{\boldsymbol{\beta}}_n) \overset{k}{\sim} \sum_{r=0}^{k-1} n^{-r/2} \sum_{|\mathbf{j}|=r} \delta^1_{(\mathbf{i};\mathbf{j})} \langle (\hat{\boldsymbol{\beta}}_n)^{\mathbf{j}} \rangle / \mathbf{j}!$$

$$= \sum_{r=0}^{k-1} n^{-r/2} \sum_{|\mathbf{j}|=r} \sum_{s+t=r} \sum_{|\mathbf{j}|=s} \delta^1_{(\mathbf{i};\mathbf{j})} \sum_{|\mathbf{l}|=t} \sum_{(0,t,\mathbf{l},\mathbf{j})} \prod_{v=0}^{t} \langle (\mathbf{h}_v)^{\mathbf{i}_v} \rangle / \mathbf{i}_v!.$$

From these we can get expansion for $F_{\mathbf{i},r}(\hat{\boldsymbol{\beta}}_n)$ by the relationship

$$F_{\mathbf{i},r}(\hat{\boldsymbol{\beta}}_n) = \sum_{t+s=r}\sum_{\mathbf{j}\geq\mathbf{i},|\mathbf{j}|=t}\varrho_{\mathbf{j}}^1\sum_{|\mathbf{l}|=s}\sum_{(0,s,\mathbf{l},\mathbf{j}-\mathbf{i})}\prod_{v=0}^{s}\langle\mathbf{h}_v^{\mathbf{i}_v}\rangle/\mathbf{i}_v!$$
$$+\sum_{t+s=r+1}\sum_{\mathbf{j}\geq\mathbf{i},|\mathbf{j}|=t}\delta_{(\mathbf{j};\mathbf{0})}^1(\hat{\boldsymbol{\beta}}_n)\sum_{|\mathbf{l}|=s}\sum_{(0,s,\mathbf{l},\mathbf{j}-\mathbf{i})}\prod_{v=0}^{s}\langle\mathbf{h}_v^{\mathbf{i}_v}\rangle/\mathbf{i}_v!$$
$$+\sum_{t+s=r+2}\sum_{\mathbf{j}\geq\mathbf{i},|\mathbf{j}|=t}\mathcal{E}_{(\mathbf{j};\mathbf{0})}^1(\hat{\boldsymbol{\beta}}_n)\sum_{|\mathbf{l}|=s}\sum_{(0,s,\mathbf{l},\mathbf{j}-\mathbf{i})}\prod_{v=0}^{s}\langle\mathbf{h}_v^{\mathbf{i}_v}\rangle/\mathbf{i}_v!.$$

With these $F_{\mathbf{i},r}(\hat{\boldsymbol{\beta}}_n)$'s we can get the expansion for

$$N_{\mathbf{i},r}(\hat{\boldsymbol{\beta}}_n) = \sum_{I_2(r,|\mathbf{i}|)}\prod_{v=1}^{r}\sum_{|\mathbf{j}|=2}^{r+2}\sum_{I_1(2,v+2,|\mathbf{k}_v|,i_v)}\prod_{u=2}^{v+2}\frac{F_{\mathbf{j},v}^{i_u}(\hat{\boldsymbol{\beta}}_n)}{i_u!(\mathbf{j}!)^{i_u}}.$$

Also, note in this case

$$-\frac{1}{2}\boldsymbol{\alpha}'\mathbf{I}_{11}(\boldsymbol{\alpha}_0,\hat{\boldsymbol{\beta}}_n)\boldsymbol{\alpha} = \sum_{|\mathbf{i}|=2}\mathcal{E}_{(\mathbf{i};\mathbf{0})}^1(\hat{\boldsymbol{\beta}}_n)\frac{\langle\boldsymbol{\alpha}^{\mathbf{i}}\rangle}{\mathbf{i}!}$$
$$\overset{k}{\sim}\sum_{r=0}^{k-1}n^{-r/2}\sum_{|\mathbf{i}|=2}\sum_{s+t=r}\sum_{|\mathbf{j}|=s}\mathcal{E}_{(\mathbf{i};\mathbf{j})}^1\sum_{|\mathbf{l}|=t}\sum_{(0,t,\mathbf{l},\mathbf{j})}\prod_{v=0}^{t}\langle\mathbf{h}_v^{\mathbf{i}_v}\rangle/\mathbf{i}_v!.$$

Then we can get the expansion for $Z_n(\boldsymbol{\alpha}+\hat{\boldsymbol{\alpha}}_n',\hat{\boldsymbol{\beta}}_n)/Z_n(\hat{\boldsymbol{\alpha}}_n',\hat{\boldsymbol{\beta}}_n)$. Going through the remaining part in the proof of Theorem 2.3, we will get two corresponding relationships of (A.10) and (A.12) for $\mathbf{g}_r$; Taking these together, as in the proof of Theorem 2.2, we get $\mathbf{0}\overset{k}{\sim}(1+O(n^{-1/2}))\mathbf{I}(\check{\boldsymbol{\alpha}}_n'-\check{\boldsymbol{\alpha}}_n'',\hat{\boldsymbol{\beta}}_n'-\hat{\boldsymbol{\beta}}_n'')'$, and thus $(\check{\boldsymbol{\alpha}}_n',\hat{\boldsymbol{\beta}}_n')'\overset{k}{\sim}(\check{\boldsymbol{\alpha}}_n'',\hat{\boldsymbol{\beta}}_n'')'$. Other proof details can be similarly obtained and are omitted.

As an alternative simple, but not rigorous, justification in the proof of Theorem 2.3, replace $\boldsymbol{\rho}_{\mathbf{i}}$, $\mathbf{a}$ and $\mathbf{I}$ with $(\boldsymbol{\rho}_{\mathbf{i}}^{1'},\mathbf{0}')'$, $\mathbf{a}_1$ and $\mathbf{I}_{11}$, set $\pi(\boldsymbol{\theta})=\pi(\boldsymbol{\alpha})\pi(\boldsymbol{\beta})$, with $\pi(\boldsymbol{\beta})$ being constant, and $W(\mathbf{d},\boldsymbol{\theta})=W(\mathbf{d}_1,\boldsymbol{\alpha})V(\mathbf{d}_2,\boldsymbol{\beta})$, with $V(\cdot,\cdot)$ be the 0–1 loss. Note $\mathbf{V}^{(1)}=\mathbf{0}$ a.e. $(\boldsymbol{\beta})$. Then, find the expression for $\mathbf{M}_{\mathbf{j},r}=(\mathbf{M}_{\mathbf{j},r}^{'1},\mathbf{0}')'$ as the way to the end of the proof in that theorem. □

PROOF OF THE FACT. (i) is a special case of (ii) with $\boldsymbol{\rho}_{\mathbf{i}}=\mathbf{0}$ $(|\mathbf{i}|=0,1)$.

(ii) The key is to find out the set $(0,s,\mathbf{l},\mathbf{i})$ for given $\mathbf{l}$ and $\mathbf{i}$. It is empty if $\mathbf{l}\neq\mathbf{0}$ and $\mathbf{i}=\mathbf{0}$. For $\mathbf{H}_1$, $s+t=1$ we must have $(s,t)=(0,1)$ or $(1,0)$. If $(s,t)=(0,1)$, for the first term $\mathbf{i}=\mathbf{l}=\mathbf{0}$ and the set $(0,0,\mathbf{l},\mathbf{i})=\{\mathbf{i}_0=\mathbf{0}\}$, so the first term is $\mathbf{I}^{-1}\boldsymbol{\rho}_{\mathbf{0}}$. For the second term, $|\mathbf{i}|=1$, so $\mathbf{i}=\mathbf{e}_j$ for some $j$, $\mathbf{l}=\mathbf{0}$ and the set $(0,0,\mathbf{l},\mathbf{i})=\{\mathbf{i}_0=\mathbf{e}_j\}$, and so the second term is $\mathbf{I}^{-1}\sum_j \boldsymbol{\Delta}_{\mathbf{e}_j}\langle\mathbf{H}_0^{\mathbf{e}_j}\rangle = \mathbf{I}^{-1}\sum_{|\mathbf{i}|=1}\boldsymbol{\Delta}_\mathbf{i}\langle\mathbf{H}_0^\mathbf{i}\rangle$. For the third term, $|\mathbf{i}|=2$, $\mathbf{l}=\mathbf{0}$ and the set $(0,0,\mathbf{l},\mathbf{i}) = \{\mathbf{i}_0 = \mathbf{i}\}$; this term is $\mathbf{I}^{-1}\sum_{|\mathbf{i}|=2}\mathbf{E}_\mathbf{i}\langle\mathbf{H}_0^\mathbf{i}\rangle/\mathbf{i}!$. If $(s,t)=(1,0)$, the summation $\sum_{|\mathbf{i}|=-1}$ in the first term is empty, and the set $(0,1,\mathbf{l},\mathbf{i})$ in the second term is empty. Also, since $t=0$, the summation $\sum_{|\mathbf{i}|=t+1,t>0}$ in the third term is empty. These give the expression for $\mathbf{H}_1$.

For $\mathbf{H}_2$, the case $(s,t) = (2,0)$ corresponds to empty summations. So we only consider $(s,t) = (0,2)$ or $(1,1)$. When $(s,t) = (0,2)$, for the first term, $|\mathbf{i}|=1$, $\mathbf{l}=\mathbf{0}$ and the set $(0,0,\mathbf{l},\mathbf{i}) = \{\mathbf{i}_0=\mathbf{i}\}$; this term is $\mathbf{I}^{-1}\sum_{|\mathbf{i}|=1}\boldsymbol{\rho}_\mathbf{i}\langle\mathbf{H}_0^\mathbf{i}\rangle$. For the second term we have $\mathbf{i} = 2\mathbf{e}_j$ for some $j$ or $\mathbf{i} = \mathbf{e}_j+\mathbf{e}_l$ for some $j\neq l$. Since $\mathbf{l}=\mathbf{0}$, the set $(0,0,\mathbf{l},\mathbf{i}) = \{\mathbf{i}_0=\mathbf{i}\}$ and this term is $\mathbf{I}^{-1}\sum_{|\mathbf{i}|=2}\boldsymbol{\Delta}_\mathbf{i}\langle\mathbf{H}_0^\mathbf{i}\rangle/\mathbf{i}!$. For the third term, $|\mathbf{i}|=3$, $\mathbf{l}=\mathbf{0}$ and $(0,0,\mathbf{l},\mathbf{i}) = \{\mathbf{i}_0=\mathbf{i}\}$, resulting in $\mathbf{I}^{-1}\sum_{|\mathbf{i}|=3}\mathbf{E}_\mathbf{i}\langle\mathbf{H}_0^\mathbf{i}\rangle/\mathbf{i}!$. When $(s,t)=(1,1)$, for the first term, $\mathbf{i}=\mathbf{0}$, $|\mathbf{l}|=1$ and the set $(0,1,\mathbf{l},\mathbf{i})$ is empty. For the second term, we have $|\mathbf{i}|=|\mathbf{l}|=1$. If $\mathbf{l}\neq\mathbf{i}$, $(0,1,\mathbf{l},\mathbf{i})$ is empty. If $\mathbf{l}=\mathbf{i}$, $(0,1,\mathbf{l},\mathbf{i}) = \{(\mathbf{i}_0,\mathbf{i}_1) = (\mathbf{0},\mathbf{i})\}$, this term is $\mathbf{I}^{-1}\sum_{|\mathbf{i}|=1}\boldsymbol{\Delta}_\mathbf{i}\langle\mathbf{H}_1^\mathbf{i}\rangle$. For the third term, $|\mathbf{i}|=2$ and so it has the form $\mathbf{i}=\mathbf{e}_j+\mathbf{e}_l$ for some $j,l$ and $|\mathbf{l}|=1$. It is easily checked that if $\mathbf{l}=\mathbf{e}_j$ or $\mathbf{e}_l$, $(0,1,\mathbf{l},\mathbf{i}) = \{(\mathbf{i}_0,\mathbf{i}_1) = (\mathbf{e}_j,\mathbf{e}_l)$ or $(\mathbf{e}_l,\mathbf{e}_j)\}$, otherwise $(0,1,\mathbf{l},\mathbf{i})$ is empty. So this term is $\mathbf{I}^{-1}\sum_{j,l=1}^d \mathbf{E}_{\mathbf{e}_j+\mathbf{e}_l}(\langle\mathbf{H}_0^{\mathbf{e}_j}\mathbf{H}_1^{\mathbf{e}_l}\rangle + \langle\mathbf{H}_0^{\mathbf{e}_l}\mathbf{H}_1^{\mathbf{e}_j}\rangle)/2$.

Note for $d=1$, it is easy to see that $\mathbf{H}_0$, $\mathbf{H}_1$ and $\mathbf{G}_1$ (for $a=2$) above coincides with those corresponding on page 496 in [17]. $h_2$ there has two extra terms $\rho_1 E_3/I^3 + \rho_1 \Delta_2/I^2$. These two terms come from $\rho_1 \sum_{I_1(0,1,1,0)} \prod_{v=0}^1 h_v^{i_v}/i_v!$ in his formula. Obviously, $I_1(0,1,1,0)$ is an empty set by definition. So these extra terms should not be there.

(iii) By Theorem 2.3, $\mathbf{G}_1 = \mathbf{H}_1 + \mathbf{Q}_1$, $\mathbf{Q}_1 = \mathbf{M}_{\mathbf{0},1}$, and it is easily checked that $\mathbf{M}_{\mathbf{0},1} = (\{\boldsymbol{\sigma}(\mathbf{a})\})^{-1}\sum_{|\mathbf{i}|=3} N_{\mathbf{i},1}\boldsymbol{\Psi}_\mathbf{i}$, and note $\boldsymbol{\Psi}_\mathbf{i} = \boldsymbol{\Psi}_\mathbf{i}^{(\mathbf{0})}(\mathbf{0})$. To evaluate $N_{\mathbf{i},1}$, for $|\mathbf{i}|=3$, note

$$N_{\mathbf{i},1} = \sum_{I_2(1,\mathbf{i})}\sum_{I_1(2,3,\mathbf{k}_1,i_1)}\prod_{|\mathbf{j}|=2}^{3}\frac{F_{\mathbf{j},1}^{u_\mathbf{j}}}{u_\mathbf{j}!(\mathbf{j}!)^{u_\mathbf{j}}}.$$

To get $I_2(1,\mathbf{i})$, we first find the corresponding $I_0(1,1,1) = \bigcup_{r\geq 0} I_1(1,1,1,r)$. It is easy to see that $I_1(1,1,1,1) = \{1\}$ and $I_1(1,1,1,r)$ is empty for $r\neq 1$, so $I_0(1,1,1,) = \{i_1=1\}$, and $I_2(1,\mathbf{i}) = \{\mathbf{k}_1 : \mathbf{k}_1 = \mathbf{i}, 2i_1 \leq |\mathbf{k}_1| \leq 3i_1\} = \{\mathbf{k}_1 : \mathbf{k}_1 = \mathbf{i}\}$. It is easy to see that $I_1(2,3,\mathbf{k}_1,i_1) = I_1(2,3,\mathbf{i},1) = \{u_\mathbf{j} : \sum_{|\mathbf{j}|=2}^3 \mathbf{j}u_\mathbf{j} = \mathbf{i}, \sum_{|\mathbf{j}|=2}^3 u_\mathbf{j} = 1\}$. There are $d^2$ of $u_\mathbf{j}$'s with $|\mathbf{j}|=2$, and $d^3$ of $u_\mathbf{j}$'s with $|\mathbf{j}|=3$, but only one of them can be 1; the rest are zeros, so $I_1(2,3,\mathbf{i},1) = \{u_\mathbf{j} : u_\mathbf{i} = 1, u_\mathbf{j}=0, \text{ for } \mathbf{j}\neq\mathbf{i}, 2\leq |\mathbf{j}|\leq 3\}$. These give, for $|\mathbf{i}|=3$, $N_{\mathbf{i},1} = \frac{F_{\mathbf{i},1}}{\mathbf{i}!}$. Refer to the definition of $F_{\mathbf{i},r}$ in Lemma 1; since $|\mathbf{i}|=3$, it is easy to see that the first two terms in $F_{\mathbf{i},1}$ are zeros as they are summations over empty sets. Also, $|\mathbf{i}|=3$, the constraints $t\geq 3$, $\mathbf{j}\geq\mathbf{i}$ and $|\mathbf{j}|=t$ gives $s=0$, $\mathbf{j}=\mathbf{i}$ and $(0, 1, \mathbf{j} - \mathbf{i}) = \{\mathbf{i}_0 = \mathbf{0}\}$, and so $F_{\mathbf{i},1} = \mathcal{E}_\mathbf{i}$. Thus

$$\mathbf{G}_1 = \mathbf{H}_1 + \mathbf{M}_{\mathbf{0},1} = \mathbf{H}_1 + (\{\boldsymbol{\sigma}(\mathbf{a})\}\mathbf{I})^{-1} \sum_{|\mathbf{i}|=3} \boldsymbol{\Psi}_\mathbf{i} \frac{\mathcal{E}_\mathbf{i}}{\mathbf{i}!}.$$

$\mathbf{G}_2 = \mathbf{H}_2 + \mathbf{Q}_2$, $\mathbf{Q}_2 = \mathbf{M}_{\mathbf{0},2} + \sum_{|\mathbf{i}|=1} \mathbf{M}_{\mathbf{i},1} \sum_{|\mathbf{l}|=1} \sum_{(1,1,\mathbf{l},\mathbf{i})} \langle \mathbf{Q}_1^\mathbf{i} \rangle / \mathbf{i}! = \mathbf{M}_{\mathbf{0},2} + \sum_{|\mathbf{i}|=1} \mathbf{M}_{\mathbf{i},1} \langle \mathbf{Q}_1^\mathbf{i} \rangle / \mathbf{i}!$. Note $\mathbf{M}_{\mathbf{0},2} = (\{\boldsymbol{\sigma}(\mathbf{a})\}\mathbf{I})^{-1} \sum_{\mathbf{i} \in \langle 2,6 \rangle} N_{\mathbf{i},2} \boldsymbol{\Psi}_\mathbf{i} = (\{\boldsymbol{\sigma}(\mathbf{a})\}\mathbf{I})^{-1} \times (\sum_{|\mathbf{i}|=3} N_{\mathbf{i},2} \boldsymbol{\Psi}_\mathbf{i} + \sum_{|\mathbf{i}|=5} N_{\mathbf{i},2} \boldsymbol{\Psi}_\mathbf{i})$ and, for $|\mathbf{i}| = 1$, $\mathbf{M}_{\mathbf{i},1} = (\{\boldsymbol{\sigma}(\mathbf{a})\}\mathbf{I})^{-1} \times \sum_{|\mathbf{j}|=3} N_{\mathbf{j}-\mathbf{i},1} \boldsymbol{\Psi}_{\mathbf{j}-\mathbf{i}}^{(\mathbf{i})}$. For $|\mathbf{i}| = 3$,

$$N_{\mathbf{i},2} = \sum_{I_2(2,\mathbf{i})} \prod_{v=1}^{2} \sum_{I_1(2,v+2,\mathbf{k}_v,i_v)} \prod_{|\mathbf{j}|=2}^{v+2} \frac{F_{\mathbf{j},v}^{u_\mathbf{j}}}{u_\mathbf{j}!(\mathbf{j}!)^{u_\mathbf{j}}}.$$

For $I_2(2,\mathbf{i})$, the corresponding $I_0(1,2,2) = \bigcup_{r \geq 0} I_1(1,2,2,r) = I_1(1,2,2,1) = \{(i_1, i_2) = (0, 1)\}$, so $I_2(2, \mathbf{i}) = \{(\mathbf{k}_1, \mathbf{k}_2) : \mathbf{k}_1 + \mathbf{k}_2 = \mathbf{i}, |\mathbf{k}_1| = 0, 2 \leq |\mathbf{k}_2| \leq 4\} = \{(\mathbf{k}_1, \mathbf{k}_2) : \mathbf{k}_1 = \mathbf{0}, \mathbf{k}_2 = \mathbf{i}\}$, $I_1(2, 3, \mathbf{k}_1, i_1) = I_1(2, 3, \mathbf{0}, 0) = \{u_\mathbf{j} : u_\mathbf{j} = 0, 2 \leq |\mathbf{j}| \leq 3\}$ and $I_1(2, 4, \mathbf{k}_2, i_2) = I_1(2, 4, \mathbf{i}, 1) = \{u_\mathbf{j} : \sum_{|\mathbf{j}|=2}^{4} \mathbf{j} u_\mathbf{j} = \mathbf{i}, \sum_{|\mathbf{j}|=2}^{4} u_\mathbf{j} = 1\} = \{u_\mathbf{j} : u_\mathbf{i} = 1, u_\mathbf{j} = 0, \text{ for } \mathbf{j} \neq \mathbf{i}, 2 \leq |\mathbf{j}| \leq 4\}$. Also, it can be checked that, for $|\mathbf{i}| = 3$, $F_{\mathbf{i},2} = \delta_\mathbf{i} + \sum_{|\mathbf{j}|=1} \mathcal{E}_{\mathbf{i}+\mathbf{j}} \langle \mathbf{H}_0^\mathbf{j} \rangle$. Recall $H_0 = \mathbf{I}^{-1} \boldsymbol{\Delta}_{\mathbf{1},\mathbf{0}}$. These give, for $|\mathbf{i}| = 3$,

$$N_{\mathbf{i},2} = \frac{F_{\mathbf{i},2}}{\mathbf{i}!} = \frac{1}{\mathbf{i}!}\left( \delta_\mathbf{i} + \sum_{|\mathbf{j}|=1} \mathcal{E}_{\mathbf{i}+\mathbf{j}} \langle (\mathbf{I}^{-1}\boldsymbol{\Delta}_\mathbf{0})^\mathbf{j} \rangle \right).$$

For $|\mathbf{i}| = 5$, $I_2(2, \mathbf{i})$ is empty, so $N_{\mathbf{i},2} = 0$. (In [17], the set $I_2(2,5)$ is also empty; but there $N_{5,2} \neq 0$ and we regard this as a mistake.) Similarly, for $|\mathbf{i}| = 1$ and $|\mathbf{j}| = 3$ with $\mathbf{j} > \mathbf{i}$,

$$N_{\mathbf{j}-\mathbf{i},1} = F_{\mathbf{j}-\mathbf{i},1}/(\mathbf{j} - \mathbf{i})! = \frac{1}{(\mathbf{j} - \mathbf{i})!}\left( \delta_{\mathbf{j}-\mathbf{i}} + \sum_{|\mathbf{l}|=1} \mathcal{E}_{\mathbf{j}-\mathbf{i}+\mathbf{l}} \langle (\mathbf{I}^{-1}\boldsymbol{\Delta}_\mathbf{0})^\mathbf{l} \rangle \right).$$

From these we get the expression for $\mathbf{G}_2$.

(iv) By the above results, the first $d_1$ components of $\mathbf{T}_1$ is $\mathbf{t}_1 = \mathbf{I}^{11}\boldsymbol{\rho}_0^1 + \mathbf{q}_1$, $\mathbf{q}_1 = \mathbf{m}_{\mathbf{0},1}$, which is the $d_1$-dimensional version of $\mathbf{M}_{\mathbf{0},1}$. So $\mathbf{m}_{\mathbf{0},1} = (\{\boldsymbol{\sigma}(\mathbf{a}_1\mathbf{I}^{11})\})^{-1} \times \sum_{|\mathbf{i}|=3} \boldsymbol{\Psi}_\mathbf{i}^1 \mathcal{E}_\mathbf{i}^1 / \mathbf{i}!$, and $\boldsymbol{\Psi}_\mathbf{i}^1$ and $\mathcal{E}_\mathbf{i}^1$ are the corresponding $d_1$-dimensional versions of their counterparts. $\square$

PROOF OF THE PROPOSITION. (i) Let $\tilde{\boldsymbol{\Delta}}_\mathbf{i}$ and $\tilde{\boldsymbol{\Delta}}_\mathbf{0}$ be the weak limits of $\boldsymbol{\Delta}_\mathbf{i}$ and $\boldsymbol{\Delta}_\mathbf{0}$. It is easy to see that $\tilde{\boldsymbol{\Delta}}_\mathbf{i} \sim N(\mathbf{0}, \mathbf{J}_\mathbf{i})$, with $\mathbf{J}_\mathbf{i} = E_{\boldsymbol{\theta}_0}[\mathbf{L}_\mathbf{i}(\mathbf{x}_1|\boldsymbol{\theta}_0)\mathbf{L}_\mathbf{i}'(\mathbf{x}_1|\boldsymbol{\theta}_0)] - \mathbf{E}_\mathbf{i}\mathbf{E}_\mathbf{i}'$, $\tilde{\boldsymbol{\Delta}}_\mathbf{0} \sim N(\mathbf{0}, \mathbf{I}^{-1})$, and $\tilde{\boldsymbol{\Delta}}_\mathbf{i}$ and $\tilde{\boldsymbol{\Delta}}_\mathbf{0}$ are jointly normal with covariance matrix, for $\mathbf{i} = \mathbf{e}_j$, $E_{\boldsymbol{\theta}_0}[\mathbf{L}_\mathbf{i}(\mathbf{x}_1|\boldsymbol{\theta}_0)\mathbf{L}_\mathbf{0}'(\mathbf{x}_1|\boldsymbol{\theta}_0)] - \mathbf{E}_\mathbf{i}\mathbf{E}_\mathbf{0}' = E_{\boldsymbol{\theta}_0}[\mathbf{L}_\mathbf{i}(\mathbf{x}_1|\boldsymbol{\theta}_0)\mathbf{L}_\mathbf{0}'(\mathbf{x}_1|\boldsymbol{\theta}_0)] =$

44   A. YUAN

$\mathbf{D}_j$. Thus $\tilde{\boldsymbol{\Delta}}_{\mathbf{i}}|\tilde{\boldsymbol{\Delta}}_{\mathbf{0}} \sim N(\mathbf{D}_j\mathbf{I}\tilde{\boldsymbol{\Delta}}_{\mathbf{0}}, \mathbf{J} - \mathbf{D}_j\mathbf{I}\mathbf{D}'_j)$, and $\langle(\mathbf{I}^{-1}\tilde{\boldsymbol{\Delta}}_{\mathbf{0}})^{\mathbf{i}}\rangle = {}_j\mathbf{I}^{-1}\tilde{\boldsymbol{\Delta}}_{\mathbf{0}} = \tilde{\boldsymbol{\Delta}}'_{\mathbf{0}}\mathbf{I}_j^{-1}$, so

$$\begin{aligned}
E_{\boldsymbol{\theta}_0}(\tilde{\boldsymbol{\Delta}}_{\mathbf{i}}\langle(\mathbf{I}^{-1}\tilde{\boldsymbol{\Delta}}_{\mathbf{0}})^{\mathbf{i}}\rangle) &= E_{\boldsymbol{\theta}_0}[E_{\boldsymbol{\theta}_0}(\tilde{\boldsymbol{\Delta}}_{\mathbf{i}}|\tilde{\boldsymbol{\Delta}}_{\mathbf{0}})\langle(\mathbf{I}^{-1}\tilde{\boldsymbol{\Delta}}_{\mathbf{0}})^{\mathbf{i}}\rangle] \\
&= E_{\boldsymbol{\theta}_0}[\mathbf{D}_j\mathbf{I}\tilde{\boldsymbol{\Delta}}_{\mathbf{0}}\langle(\mathbf{I}^{-1}\tilde{\boldsymbol{\Delta}}_{\mathbf{0}})^{\mathbf{i}}\rangle] = E_{\boldsymbol{\theta}_0}[\mathbf{D}_j\mathbf{I}\tilde{\boldsymbol{\Delta}}_{\mathbf{0}j}\mathbf{I}^{-1}\tilde{\boldsymbol{\Delta}}_{\mathbf{0}}] \\
&= \mathbf{D}_j\mathbf{I}E_{\boldsymbol{\theta}_0}[\tilde{\boldsymbol{\Delta}}_{\mathbf{0}}\tilde{\boldsymbol{\Delta}}'_{\mathbf{0}}]\mathbf{I}_j^{-1} = \mathbf{D}_j\mathbf{I}\mathbf{I}^{-1}\mathbf{I}_j^{-1} = \mathbf{D}_j\mathbf{I}_j^{-1}.
\end{aligned}$$

Similarly, for $|\mathbf{i}| = 2$, $\mathbf{i} = \mathbf{e}_i + \mathbf{e}_j$ for some $(i,j)$. Thus,

$$\begin{aligned}
E_{\boldsymbol{\theta}_0}\langle(\mathbf{I}^{-1}\tilde{\boldsymbol{\Delta}}_{\mathbf{0}})^{\mathbf{i}}\rangle &= E_{\boldsymbol{\theta}_0}(({}_i\mathbf{I}^{-1}\tilde{\boldsymbol{\Delta}}_{\mathbf{0}})({}_j\mathbf{I}^{-1}\tilde{\boldsymbol{\Delta}}_{\mathbf{0}})') \\
&= {}_i\mathbf{I}^{-1}E_{\boldsymbol{\theta}_0}(\tilde{\boldsymbol{\Delta}}_{\mathbf{0}}\tilde{\boldsymbol{\Delta}}'_{\mathbf{0}})\mathbf{I}_j^{-1} = {}_i\mathbf{I}^{-1}\mathbf{I}^{-1}\mathbf{I}_j^{-1},
\end{aligned}$$

and now the result follows using Fact (i) and taking the corresponding summations.

(ii) Note $E_{\boldsymbol{\theta}_0}(\tilde{\mathbf{G}}_1) = E_{\boldsymbol{\theta}_0}(\tilde{\mathbf{H}}_1^\circ) + \mathbf{I}^{-1}\boldsymbol{\rho}_0 + \mathbf{M}_{\mathbf{0},1}$, and the result follows.

(iii) The proof is similar and is omitted. $\square$
**Acknowledgments.** I am grateful to the reviewers for their helpful suggestions and comments, especially to Guang Cheng for the topic on objective prior.

HOWARD UNIVERSITY
2216 SIXTH STREET, N.W., SUITE 206
WASHINGTON, DISTRICT OF COLUMBIA 20059
USA
E-MAIL: ayuan@howard.edu